\PassOptionsToPackage{unicode}{hyperref}
\PassOptionsToPackage{hyphens}{url}
\documentclass[
  11pt,
]{article}
\usepackage{amsmath,amssymb}
\usepackage{iftex}
\ifPDFTeX
  \usepackage[T1]{fontenc}
  \usepackage[utf8]{inputenc}
  \usepackage{textcomp} 
\else 
  \usepackage{unicode-math} 
  \defaultfontfeatures{Scale=MatchLowercase}
  \defaultfontfeatures[\rmfamily]{Ligatures=TeX,Scale=1}
\fi
\usepackage{lmodern}
\ifPDFTeX\else
\fi
\IfFileExists{upquote.sty}{\usepackage{upquote}}{}
\IfFileExists{microtype.sty}{
  \usepackage[]{microtype}
  \UseMicrotypeSet[protrusion]{basicmath} 
}{}
\makeatletter
\@ifundefined{KOMAClassName}{
  \IfFileExists{parskip.sty}{%
    \usepackage{parskip}
  }{
    \setlength{\parindent}{0pt}
    \setlength{\parskip}{6pt plus 2pt minus 1pt}}
}{
  \KOMAoptions{parskip=half}}
\makeatother
\usepackage{xcolor}
\usepackage[margin=1in]{geometry}
\usepackage{longtable,booktabs,array}
\usepackage{calc} 
\usepackage{etoolbox}
\makeatletter
\patchcmd\longtable{\par}{\if@noskipsec\mbox{}\fi\par}{}{}
\makeatother
\IfFileExists{footnotehyper.sty}{\usepackage{footnotehyper}}{\usepackage{footnote}}
\makesavenoteenv{longtable}
\providecommand{\tightlist}{%
  \setlength{\itemsep}{0pt}\setlength{\parskip}{0pt}}
\newlength{\cslhangindent}
\newlength{\csllabelwidth}
\newlength{\cslentryspacingunit} 
\newenvironment{CSLReferences}[2] 
 {
  \setlength{\parindent}{0pt}
  \ifodd #1
  \let\oldpar\par
  \def\par{\hangindent=\cslhangindent\oldpar}
  \fi
  \setlength{\parskip}{#2\cslentryspacingunit}
 }%
 {}
\usepackage{calc}

\newcommand{\CSLLeftMargin}[1]{\parbox[t]{\csllabelwidth}{#1}}
\newcommand{\CSLRightInline}[1]{\parbox[t]{\linewidth - \csllabelwidth}{#1}\break}

\usepackage{tikz}
\usetikzlibrary{arrows.meta}

\ifLuaTeX
  \usepackage{selnolig}  
\fi
\IfFileExists{bookmark.sty}{\usepackage{bookmark}}{\usepackage{hyperref}}
\IfFileExists{xurl.sty}{\usepackage{xurl}}{} 
\hypersetup{
  pdftitle={Maximal Hamiltonicity of realization graphs of degree sequences},
  pdfauthor={Jeffrey S. Baggett; Department of Mathematics and Statistics, University of Wisconsin--La Crosse},
  hidelinks,
  pdfcreator={LaTeX via pandoc}}

\title{Maximal Hamiltonicity of realization graphs of degree sequences}
\author{Jeffrey S. Baggett \and Department of Mathematics and
Statistics, University of Wisconsin--La Crosse}
\date{}

\begin{document}
\maketitle
\begin{abstract}
We prove that the realization graph of every graphical degree sequence
is maximally Hamiltonian: it is Hamilton-laceable when bipartite on more
than one vertex, and Hamilton-connected otherwise. This answers Problem
P59 of Mütze's survey of combinatorial Gray codes {[}24{]}, and the
Hamiltonicity question recorded as open by Barrus {[}5{]}, in the
strongest form either admits. The argument is an induction on the number
of ground vertices, cutting the realization graph at a single ground
vertex into fibers and the quotient they lie over. The proof is
formalized in Lean 4 and checked by its kernel, with seven results cited
from the literature and nothing else assumed.

Its engine is a classification. The realizable neighborhoods of a ground
vertex form a \textbf{shifted} family --- one closed under replacing an
element by a smaller one --- and the quotient is the Johnson graph of
that family. Such a Johnson graph can fail to be Hamilton-connected, and
we determine exactly when: the failures are one explicit family of
examples, the \textbf{Y-families}, and each of them fails between a
single pair of its members. A shifted family with a greatest member
never fails, and those families are exactly the shifted matroids, where
the conclusion already follows from the theorem of Naddef and
Pulleyblank on the graphs of 0/1-polytopes {[}25,26{]}. The obstruction
lives entirely outside the matroid case, which is why it has not been
met before.
\end{abstract}

\hypertarget{introduction}{%
\section{1. Introduction}\label{introduction}}

Take a graphical degree sequence and form the graph whose vertices are
its labeled realizations, with two joined when a single 2-switch carries
one to the other. This \textbf{realization graph} is connected --- Ryser
proved it in the matrix form {[}31{]} and Fulkerson, Hoffman and
McAndrew for degree sequences {[}16{]} --- and the question of whether
it is Hamiltonian has been open in general since Brualdi raised it for
bipartite degree sequences --- equivalently, for the \textbf{interchange
graph} of a class of zero-one matrices with fixed margins, whose
vertices are the matrices of the class and whose edges are single
interchanges {[}8, Problem 3.7{]}. The split case is the image of that
problem under the correspondence of Arikati and Peled {[}3{]}. Barrus
surveys the realization-graph viewpoint and records the general question
as open; Mütze's survey {[}24{]} lists it as Problem P59, with the split
case as Problem P60. The question belongs to combinatorial generation as
much as to graph theory: a Hamilton path in a flip graph is a Gray code
for the objects it lists, one local move at a time, which is why it
appears in a survey of Gray codes at all.

Three partial results bound what was known. Li and Zhang {[}20{]} proved
the interchange graph Hamiltonian when every column sum is one. Arikati
and Peled {[}3{]} proved the realization graph Hamiltonian when the
degree sequence has majorization gap one, and Barrus {[}5{]} proved it
Hamiltonian whenever it is triangle-free; Section 3 strengthens the
third of these to maximal Hamiltonicity.

We prove the strongest form the question admits. A realization graph can
be bipartite, and a bipartite graph on more than two vertices is never
Hamilton-connected, so the property to aim at is maximal Hamiltonicity:
Hamilton-laceable when bipartite on more than one vertex,
Hamilton-connected otherwise.

\textbf{Theorem 1.1.} \emph{The realization graph of every graphical
degree sequence is maximally Hamiltonian.}

\hypertarget{what-this-paper-continues}{%
\subsection{1.1 What this paper
continues}\label{what-this-paper-continues}}

This paper is the companion announced in our earlier work on interchange
graphs, which proved the same statement for the split case, Brualdi's
Problem P60. That paper identified the obstruction to extending its own
method, and this paper removes it.

The obstruction is worth stating at once, since it shapes everything
here. The proof in the split case turns on the quotient at a chosen line
being the basis graph of a matroid, so that the theorem of Naddef and
Pulleyblank {[}25,26{]} applies. For a general degree sequence the
corresponding quotient, the family of neighborhoods of a fixed ground
vertex that occur in some realization, is not a matroid basis family;
there are degree sequences on eleven ground vertices at which no ground
vertex has a matroid basis family. The classical theory of matroid basis
graphs therefore says nothing about these quotients, and new structure
theory is required.

Section 7 supplies it. The quotient family is always \textbf{shifted},
meaning closed under replacing an element by a smaller one, and we
determine exactly when the Johnson graph of a shifted family fails to be
Hamilton-connected. That classification is the mathematical content of
this paper, and we state it here because it is of independent interest
and because the failures, not the successes, are what Section 8 has to
work around.

\textbf{Theorem 1.2.} \emph{Let \(F\) be a shifted family of
\(k\)-subsets of \([n]\) and let \(A ≠ B\) be members of \(F\). Then the
Johnson graph of \(F\) has a Hamilton \(A\)--\(B\) path unless \(F\) is
a Y-family and \(\{A,B\}\) is its universal pair, in which case it has
none.}

Section 7 settles both terms, but they are quick to say. A
\textbf{Y-family} is named by two of its members: it is the set of
\(k\)-sets Gale-below \([k−1] ∪ \{r\}\) or below \([k+1] ∖ \{s\}\) ---
the first \(k−1\) numbers with one far element adjoined, and the first
\(k+1\) numbers with one small element removed. Its Johnson graph is
then two cliques glued along a shared pair, and that pair is
\textbf{universal} in the literal sense, its two members joined to every
other. Deleting them leaves the two cliques disconnected, so no Hamilton
path can run between them. In the smallest Y-family,
\(F = \{12,13,14,23\}\), that pair is \(\{12,13\}\) and
\(J(F) = K_4 − e\) with \(14\)--\(23\) missing. The same \(F\) shows the
geometry parting company with the combinatorics:
\(\operatorname{conv}(F)\), the convex hull of the indicator vectors of
its members, is a 3-simplex whose skeleton is \(K_4\) and is therefore
Hamilton-connected, while \(J(F)\) is not. The polytope skeleton
acquires edges the Johnson graph does not have, and it acquires them at
every exception --- though not only there, as Section 7.9 shows with a
shifted family that is not a Y-family.

Two features of Theorem 1.2 matter downstream. The exception is
completely classified, so a proof that needs to avoid it can decide
whether it is present. And every Y-family has two incomparable
Gale-maximal members, so a family with a greatest member is never a
Y-family. The families with a greatest member are exactly the basis
families of shifted matroids. So the exception cannot occur in the case
the existing literature studies, which is why it has not been met
before.

\hypertarget{the-shape-of-the-argument}{%
\subsection{1.2 The shape of the
argument}\label{the-shape-of-the-argument}}

\textbf{How the proof works.} A vertex of \(G(d)\) is a labeled graph
with degree function \(d\), and an edge is a single 2-switch, so a
Hamilton path of \(G(d)\) is a walk through every realization, one
switch at a time, from a prescribed start to a prescribed finish. The
proof is an induction on the number of ground vertices, and three of its
four branches are short: a bipartite \(G(d)\) is settled by
classification, a decomposable \(d\) inherits the conclusion from its
factors, and three realizations make \(G(d)\) a triangle.

The fourth branch is the argument. Choose a ground vertex, the
\textbf{pivot}, and sort the realizations by what that vertex is joined
to; each class is a \textbf{fiber}, and a fiber is again a realization
graph on a ground with one fewer vertex, so the induction hypothesis
applies inside it. Contracting the fibers leaves a \textbf{quotient},
which Section 5 identifies exactly: its vertices are the neighborhoods
the pivot actually takes, and two of them are adjacent when they differ
by exchanging a single element. The Hamilton path is then built in one
pass. Walk the quotient, visiting every fiber once; span each fiber as
you pass through it, by the induction hypothesis; and cross from one
fiber to the next along an edge of \(G(d)\) that the interface supplies.

Two obstacles remain, and Sections 6 to 9 are what clear them. A fiber
is entered where the arriving crossing lands and left where the
departing crossing needs, so both of its ends are dictated from outside
and the induction has to supply a path between two prescribed vertices.
That is why the statement proved is maximal Hamiltonicity rather than
traceability: a weaker hypothesis on the fibers would not compose. And
at one fiber both ends are forced at once, by the two chains built
inward from the prescribed start and finish; if that fiber is bipartite
and its two forced ends fall in the same color class, no spanning walk
between them exists. So one fiber along the route is required to be
non-bipartite, and we call it the \textbf{buffer}. Section 6 produces a
pivot carrying one, Section 7 supplies the quotient walk and classifies
the single configuration in which it fails, and Section 8 shows the
pivot can be chosen so that configuration does not arise.

\textbf{Section by section.} Section 2 fixes notation and records what a
2-switch does to a realization. Sections 3 and 4 then reduce Theorem 1.1
to a single case, in which the degree sequence is active and
Tyshkevich-indecomposable and its realization graph is non-bipartite and
larger than a triangle. The bipartite case is settled by classification,
and both reductions are proved directly.

Section 5 cuts the realization graph at a ground vertex and proves that
the quotient is exactly the Johnson graph of a shifted family, together
with the four interface facts the later sections consume. Section 6
produces a ground vertex that separates the two prescribed realizations
and has a non-bipartite fiber; in the split case this is the
corresponding lemma of the earlier paper, and what Section 6 adds is the
extension to every degree sequence. The \textbf{prescribed} pair is
simply the pair the Hamilton path is required to join: maximal
Hamiltonicity asks for one between \emph{every} two realizations, so the
argument fixes an arbitrary pair at the outset and everything downstream
is relative to it.

Section 7 is the engine, and the longest part of the paper. It settles
exactly when the Johnson graph of a shifted family has a Hamilton path
between two prescribed members: always, except when the family is one of
an explicit three-parameter list and the pair is its two Gale-least
members. That is Theorem 1.2. It is proved from the foundations alone,
it is the part with no counterpart in the split case, and everything
after it consumes the classification rather than reopening it.

Sections 8 and 9 supply what the assembly needs beyond those. Section 8
shows the ground vertex can be chosen so that the prescribed pair is not
the exception of Theorem 1.2, which is a strictly stronger demand than
Section 6 meets and does not follow from it. Section 9 proves the two
interface extensions that carry a spanning walk from one fiber into the
next, one of them avoiding a single prescribed target. Section 10
assembles a Hamilton path in one pass along a Hamilton path of the
quotient, and validates it step by step.

The proof is complete and it is machine-checked. Every numbered result
of Sections 3 through 10 is proved in our Lean development, and, except
where the text says otherwise, the argument given here is the argument
the kernel checks. What the development assumes rather than proves is
seven results from the literature, listed in Section 10.4 and matched
there against the axiom trace. Section 11.1 gives the accounting.

Computation appears in this paper for the things a proof does not
settle: that the exception of Theorem 7.7 actually occurs, that the
routes we say are unavailable are unavailable, and that our formal
definitions describe realization graphs and not some other object.
Section 11 reports those and nothing else.

\hypertarget{independent-work}{%
\subsection{1.3 Independent work}\label{independent-work}}

Independently and in parallel with our earlier paper {[}4{]}, posted to
arXiv on 16 July 2026, Hladík and Fink {[}18{]} proved that the
realization graph of every degree sequence has a Hamilton path from any
prescribed vertex. They announced the result at the Czech and Slovak
student conference SVOČ in May 2026 and posted their preprint later that
July, and their paper records the two efforts as independent.

Their argument passes through the claim that the family of neighborhoods
of a fixed vertex has a greatest element, which carries their Lemma 3.2
and through it their Theorem 3.3. That claim does not hold in general:
Section 5.5 exhibits a five-vertex degree sequence whose family has two
maximal elements. The gap is in the route rather than the destination.
\textbf{Their theorem is true} --- Corollary 10.2 recovers it from ours
in two lines --- and what the greatest-element claim concealed is a
genuine exception, which Section 7 isolates and the corrected statement
there accounts for. We raised the point with the authors, who have
confirmed it and are preparing a revision. What this paper proves is the
stronger conclusion: maximal Hamiltonicity controls both endpoints where
a Hamilton path from a prescribed start controls one.

\hypertarget{conventions}{%
\section{2. Conventions}\label{conventions}}

All graphs are finite and simple. A \textbf{degree function} on a vertex
set \(V\) of size \(n\) is a map \(d : V → \{0, 1, …, n−1\}\); it is
\textbf{graphical} when some graph on \(V\) realizes it, and such a
graph is a \textbf{realization}. Whether \(d\) is graphical is decided
by \(d\) alone, through the Erdős--Gallai inequalities {[}13{]};
Sections 5 and 8 use them, and nothing before that does.

\begin{figure}[tbp]
\centering
\begin{tikzpicture}[x=1cm,y=1cm,
  minidot/.style={circle,fill=black!72,inner sep=0.9pt},
  miniedge/.style={draw=black!52,line width=0.42pt},
  gedge/.style={draw=black!30,line width=0.7pt},
  ann/.style={font=\scriptsize,text=black!62}]
\draw[gedge] (1.929e-16,3.15) -- (2.728,-1.575);
\draw[gedge] (1.929e-16,3.15) -- (-2.728,-1.575);
\draw[gedge] (1.929e-16,3.15) -- (-2.728,1.575);
\draw[gedge] (1.929e-16,3.15) -- (2.728,1.575);
\draw[gedge] (2.728,-1.575) -- (-2.728,-1.575);
\draw[gedge] (2.728,-1.575) -- (-5.786e-16,-3.15);
\draw[gedge] (2.728,-1.575) -- (2.728,1.575);
\draw[gedge] (-2.728,-1.575) -- (-5.786e-16,-3.15);
\draw[gedge] (-2.728,-1.575) -- (-2.728,1.575);
\draw[gedge] (-5.786e-16,-3.15) -- (-2.728,1.575);
\draw[gedge] (-5.786e-16,-3.15) -- (2.728,1.575);
\draw[gedge] (-2.728,1.575) -- (2.728,1.575);
\fill[white] (1.929e-16,3.15) circle (0.72);
\draw[miniedge] (2.247e-16,3.67) -- (-0.4945,3.311);
\draw[miniedge] (2.247e-16,3.67) -- (-0.3056,2.729);
\draw[miniedge] (2.247e-16,3.67) -- (0.3056,2.729);
\draw[miniedge] (-0.4945,3.311) -- (-0.3056,2.729);
\draw[miniedge] (0.3056,2.729) -- (0.4945,3.311);
\node[minidot] at (2.247e-16,3.67) {};
\node[minidot] at (-0.4945,3.311) {};
\node[minidot] at (-0.3056,2.729) {};
\node[minidot] at (0.3056,2.729) {};
\node[minidot] at (0.4945,3.311) {};
\node[ann] at (2.578e-16,4.21) {$G_{1}$};
\fill[white] (2.728,-1.575) circle (0.72);
\draw[miniedge] (2.728,-1.055) -- (2.233,-1.414);
\draw[miniedge] (2.728,-1.055) -- (2.422,-1.996);
\draw[miniedge] (2.728,-1.055) -- (3.034,-1.996);
\draw[miniedge] (2.233,-1.414) -- (3.034,-1.996);
\draw[miniedge] (2.422,-1.996) -- (3.223,-1.414);
\node[minidot] at (2.728,-1.055) {};
\node[minidot] at (2.233,-1.414) {};
\node[minidot] at (2.422,-1.996) {};
\node[minidot] at (3.034,-1.996) {};
\node[minidot] at (3.223,-1.414) {};
\node[ann] at (3.646,-2.105) {$G_{2}$};
\fill[white] (-2.728,-1.575) circle (0.72);
\draw[miniedge] (-2.728,-1.055) -- (-3.223,-1.414);
\draw[miniedge] (-2.728,-1.055) -- (-3.034,-1.996);
\draw[miniedge] (-2.728,-1.055) -- (-2.422,-1.996);
\draw[miniedge] (-3.223,-1.414) -- (-2.233,-1.414);
\draw[miniedge] (-3.034,-1.996) -- (-2.422,-1.996);
\node[minidot] at (-2.728,-1.055) {};
\node[minidot] at (-3.223,-1.414) {};
\node[minidot] at (-3.034,-1.996) {};
\node[minidot] at (-2.422,-1.996) {};
\node[minidot] at (-2.233,-1.414) {};
\node[ann] at (-3.646,-2.105) {$G_{3}$};
\fill[white] (-5.786e-16,-3.15) circle (0.72);
\draw[miniedge] (-5.468e-16,-2.63) -- (-0.4945,-2.989);
\draw[miniedge] (-5.468e-16,-2.63) -- (-0.3056,-3.571);
\draw[miniedge] (-5.468e-16,-2.63) -- (0.4945,-2.989);
\draw[miniedge] (-0.4945,-2.989) -- (0.3056,-3.571);
\draw[miniedge] (-0.3056,-3.571) -- (0.3056,-3.571);
\node[minidot] at (-5.468e-16,-2.63) {};
\node[minidot] at (-0.4945,-2.989) {};
\node[minidot] at (-0.3056,-3.571) {};
\node[minidot] at (0.3056,-3.571) {};
\node[minidot] at (0.4945,-2.989) {};
\node[ann] at (-7.734e-16,-4.21) {$G_{4}$};
\fill[white] (-2.728,1.575) circle (0.72);
\draw[miniedge] (-2.728,2.095) -- (-3.223,1.736);
\draw[miniedge] (-2.728,2.095) -- (-2.422,1.154);
\draw[miniedge] (-2.728,2.095) -- (-2.233,1.736);
\draw[miniedge] (-3.223,1.736) -- (-3.034,1.154);
\draw[miniedge] (-3.034,1.154) -- (-2.422,1.154);
\node[minidot] at (-2.728,2.095) {};
\node[minidot] at (-3.223,1.736) {};
\node[minidot] at (-3.034,1.154) {};
\node[minidot] at (-2.422,1.154) {};
\node[minidot] at (-2.233,1.736) {};
\node[ann] at (-3.646,2.105) {$G_{5}$};
\fill[white] (2.728,1.575) circle (0.72);
\draw[miniedge] (2.728,2.095) -- (2.422,1.154);
\draw[miniedge] (2.728,2.095) -- (3.034,1.154);
\draw[miniedge] (2.728,2.095) -- (3.223,1.736);
\draw[miniedge] (2.233,1.736) -- (2.422,1.154);
\draw[miniedge] (2.233,1.736) -- (3.034,1.154);
\node[minidot] at (2.728,2.095) {};
\node[minidot] at (2.233,1.736) {};
\node[minidot] at (2.422,1.154) {};
\node[minidot] at (3.034,1.154) {};
\node[minidot] at (3.223,1.736) {};
\node[ann] at (3.646,2.105) {$G_{6}$};
\node[ann] at (5.55,1.46) {the five ground};
\node[ann] at (5.55,1.22) {vertices, fixed};
\node[ann] at (5.55,-1.3) {$u$ with degree $d(u)$};
\node[minidot] at (5.55,0.52) {};
\node[ann] at (5.55,0.9256) {$1\,(3)$};
\node[minidot] at (5.055,0.1607) {};
\node[ann] at (4.67,0.286) {$2\,(2)$};
\node[minidot] at (5.244,-0.4207) {};
\node[ann] at (5.006,-0.7488) {$3\,(2)$};
\node[minidot] at (5.856,-0.4207) {};
\node[ann] at (6.094,-0.7488) {$4\,(2)$};
\node[minidot] at (6.045,0.1607) {};
\node[ann] at (6.43,0.286) {$5\,(1)$};
\end{tikzpicture}
\caption{The realization graph $G(d)$ for $d = (3,2,2,2,1)$, drawn whole. Each of its six vertices is a realization, shown as the graph it is, with the five ground vertices in the same five positions every time and keyed at the right by number and degree. The six are named $G_1$ to $G_6$, and Figure~\ref{fig:fibers} sorts these same six, under those same names, into the fibers of a pivot. Two are joined when their edge sets differ in exactly four places, which is a single $2$-switch. The three missing joins are the three long diagonals, so $G(d)$ is $K_6$ minus a perfect matching, and it is Hamilton-connected. Section~5.1 takes the same sequence as its running example and cuts this graph into fibers.}
\label{fig:realizationgraph}
\end{figure}
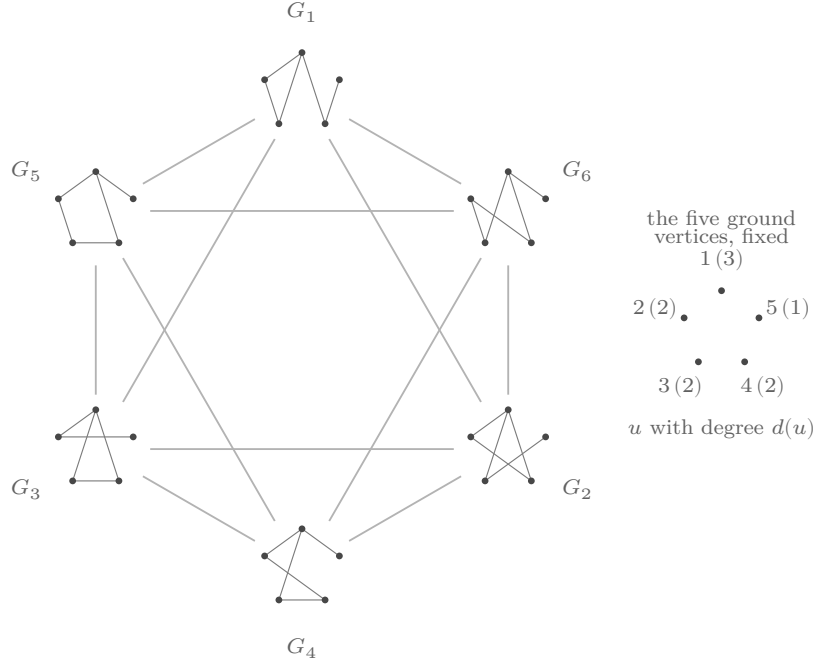

Two bracket conventions run through the paper and are not the same.
\([n]\) is always the set \(\{1, …, n\}\), and \([t]\) its prefix of
length \(t\). A bracket around a \emph{statement} is an indicator:
\([t ∈ S]\) is \(1\) when \(t ∈ S\) and \(0\) otherwise. The first
encloses a number, the second a proposition, and they never occur in the
same position. We call the elements of \(V\) the \textbf{ground
vertices}; Section 7 calls the same set \emph{the ground} and writes it
\([n]\). \emph{Ground order \(n\)} means \(|V| = n\), and the induction
of Section 4 is on it. Degree sequences are labeled throughout: two
realizations that differ by a relabeling are distinct vertices of
\(G(d)\). A \textbf{2-switch} replaces edges \(ab\) and \(cd\), on four
distinct vertices with \(ad\) and \(cb\) absent, by \(ad\) and \(cb\);
those four vertices are its \textbf{support}, and any four ground
vertices taken together are a \textbf{quartet}. Two realizations are
adjacent in \(G(d)\) exactly when their edge sets have symmetric
difference of size four. Figure \ref{fig:realizationgraph} draws one,
small enough to take in at once; Section 5.1 returns to the same degree
sequence and cuts that graph into fibers.

\textbf{Alternating structure, which is where 2-switches come from.}
Given two realizations \(G\) and \(H\), color the edges of
\(E(G) △ E(H)\) by which of the two they belong to. A closed walk in the
symmetric difference is \textbf{alternating} when consecutive edges have
different colors. Alternation is forced by the degrees: at any vertex
the edges split into shared, \(G\)-only and \(H\)-only, and since the
shared ones count toward both degrees, equal degrees leave equally many
\(G\)-only and \(H\)-only edges at every vertex. That balance is exactly
the condition for the symmetric difference to decompose into closed
alternating \textbf{trails}, and a 2-switch is the smallest case: an
alternating four-cycle, swapped by exchanging the colors on it.

\textbf{Trails, not cycles, and the distinction is not pedantic.} A
trail may revisit a vertex; a cycle may not. An alternating cycle has
even length, since the colors alternate around it, so an odd cycle never
alternates on its own --- but two odd cycles sharing a vertex, traversed
as a figure-eight, do. The symmetric difference of two realizations
therefore need not decompose into alternating cycles even though it
always decomposes into alternating trails. Section 1.3 records an
independent argument that passes through this point, and the
Acknowledgments give its authors' own diagnosis: the decomposition into
edge-disjoint alternating cycles is available when the underlying graph
is bipartite, and in general the obstruction is exactly the figure-eight
above.

A graph is \textbf{maximally Hamiltonian} when it is Hamilton-laceable
if bipartite, meaning a Hamilton path joins every pair of vertices in
opposite color classes, and Hamilton-connected otherwise; in the
bipartite case we require the proper 2-coloring to be onto, so that both
color classes are nonempty. The graphs \(K_1\) and \(K_2\) are maximally
Hamiltonian by convention. \textbf{A maximally Hamiltonian bipartite
graph on more than one vertex is balanced}, since a Hamilton path
between opposite classes alternates between them, so the two classes
differ in size by at most the parity of the order and, the path having
both endpoints in opposite classes, are equal. We use this in Section 9
and record it here because it is a property of the definition rather
than of realization graphs. Realization graphs are connected, a
consequence of a theorem of Fulkerson, Hoffman and McAndrew {[}16{]} as
Barrus records {[}5{]}, and \(G(d)\) has one vertex exactly when \(d\)
is the degree function of a threshold graph. A \textbf{spanning walk} of
a graph from \(x\) to \(y\) is a Hamilton path with those two endpoints;
we say \emph{walk} rather than \emph{path} when it is to be one segment
of a longer path, which is how Sections 9 and 10 use them. It never
repeats a vertex: the word marks the role the object plays, not a weaker
object.

For a ground vertex \(v\) and a set \(S ⊆ V ∖ \{v\}\) we write \(f_S\)
for the \textbf{residual degree function}, which lowers by one the
degree of each member of \(S\) and removes \(v\). Throughout, \(S\)
ranges over the sets that occur as a neighborhood of the pivot, so it
depends on \(v\) even though the letter does not record it; the pivot is
fixed wherever \(S\) appears, and Section 5.1 fixes it for the whole of
that section. Section 5 develops the fibers over \(v\), the family of
realizable neighborhoods, and the interface between adjacent fibers, and
fixes the notation \(Φ_S\), \(I_v\), \(w_{ab}\) and \(δ\) used from
Section 6 onward.

We use \(G\) and \(H\) for realizations, \(A\) and \(B\) for their
projections into a quotient, and \(X\) and \(Y\) for members of a set
family. Wherever a prescribed pair of realizations is in play ---
Sections 5 through 10 --- \(A\) and \(B\) are those projections and
nothing else. The other uses all sit before any quotient is taken or are
quoted: Section 3.3 writes \(B\) for a partial product, Section 4.2
writes \(A\) and \(B\) for the clique and independent sides of a
splitted factor, and \textbf{C4} and Section 12 keep the letters their
sources use. Section 7 is written in the language of set families alone
and does not refer to degree sequences; the two vocabularies meet only
through Corollary 5.5.

A graph with a proper 2-coloring whose two color classes have equal size
is \textbf{balanced}. Such a graph is \textbf{paired two-disjoint-path
coverable} when, for any two prescribed pairs of terminals that are
opposite-colored and together involve four distinct vertices, it has two
vertex-disjoint paths joining the prescribed pairs and covering every
vertex. For balanced bipartite graphs of order at least four this is
strictly stronger than Hamilton-laceability, and Sections 3 and 4 need
the stronger form: laceability alone does not survive a Cartesian
product, and the paired property does.

This paper cites results from outside it in the ordinary way, at the
point of use. Five of them are labeled \textbf{C1} through \textbf{C5}
--- \textbf{C} for \emph{cited} --- and each is stated in full where it
is used, \textbf{C1} to \textbf{C3} in Section 3.2, \textbf{C4} in
Section 4.2 and \textbf{C5} in Section 6.4, because they are quoted from
developments accompanying this one and a reader should not have to fetch
them. Nothing else is labeled.

The arguments below stand as mathematics, independently of any
formalization: every proof is written to be read and checked on the
page. Section 11 separately reports what a machine has verified, and
which results the formal development charges as assumptions rather than
proving. A reader with no interest in that may skip it without losing a
step of the proof.

Tyshkevich's canonical decomposition {[}34{]} is cited in Section 4 for
a structural reason explained there, and the argument does not depend on
it.

\hypertarget{the-bipartite-case}{%
\section{3. The bipartite case}\label{the-bipartite-case}}

The bipartite case closes by classification, and it is the only case
that does. Barrus's theorem {[}5{]} describes every triangle-free
realization graph as a Cartesian product of graphs from a short fixed
list, so the result follows from a product argument and one base case.
Everything after this section is the machinery for the non-bipartite
case, where no classification is available and the graphs must be taken
apart in place.

\hypertarget{the-classification}{%
\subsection{3.1 The classification}\label{the-classification}}

By Barrus's classification, for a graphical \(d\) the conditions
``\(G(d)\) is triangle-free'' and ``\(G(d)\) is bipartite'' are
equivalent, and the bipartite realization graphs are exactly the
Cartesian products

\[
G(d) ≅ T_{k_1} □ ⋯ □ T_{k_r} □ E,
\]

where each \(CT_k\) is a transposition graph --- permutations of
\([k]\), adjacent when they differ by one transposition --- and \(E\) is
either absent or a single copy of the crown graph \(K_{6,6} − 6K_2\).
Barrus's \(CT_k\) is the complete transposition graph \(CT_k\) of the
companion paper {[}4{]}, and his printed identifications \(T_2 ≅ K_2\)
and \(T_3 ≅ K_{3,3}\) agree with ours.

We use this equivalence again in Section 6, where a fiber is
non-bipartite exactly when it contains a triangle.

\hypertarget{the-results-we-quote}{%
\subsection{3.2 The results we quote}\label{the-results-we-quote}}

Five results from outside this paper are used in Sections 3, 4 and 6,
each stated in the form in which it is used rather than merely named
there. They are labeled \textbf{C1} through \textbf{C5} throughout ---
\textbf{C} for \emph{cited} --- so that a result carrying a letter is
never one of ours. The three this section needs, \textbf{C1} to
\textbf{C3}, are stated here; \textbf{C4} is stated in Section 4.2 where
the product lift is applied and \textbf{C5} in Section 6.4 where the
buffer is taken, so that neither is a toll on the way to Theorem 3.1.
Section 10.4 lists all five in one place.

Two definitions are needed first. A \textbf{weld} of graphs
\(H_1, …, H_ℓ\) of equal order is their disjoint union together with a
perfect matching between every pair \(H_i, H_j\). A
\textbf{transposition-like graph of rank one} is a Hamilton-connected
graph, or a graph carrying a proper surjective 2-coloring for which it
is Hamilton-laceable; of \textbf{rank \(r ≥ 2\)}, a weld of \(ℓ ≥ r\)
transposition-like graphs of rank \(r − 1\) of equal order. Unfolding
the recursion presents a rank-\(r\) graph as a welding tower.

\begin{figure}[tbp]
\centering
\begin{tikzpicture}[x=1cm,y=1cm,
  dot/.style={circle,fill=black!78,inner sep=1.5pt},
  leafe/.style={draw=black!85,line width=1.5pt},
  matche/.style={draw=black!42,line width=0.7pt},
  ann/.style={font=\scriptsize,text=black!62},
  box/.style={draw=black!35,rounded corners=2pt,dashed,line width=0.6pt}]
\draw[matche] (5,0.55) -- (2.5,0.55);
\draw[matche] (5,0.55) -- (0,-0.55);
\draw[matche] (2.5,0.55) -- (0,0.55);
\draw[matche] (5,-0.55) -- (0,0.55);
\draw[matche] (5,-0.55) -- (2.5,-0.55);
\draw[matche] (2.5,-0.55) -- (0,-0.55);
\draw[leafe] (5,0.55) -- (5,-0.55);
\draw[leafe] (2.5,0.55) -- (2.5,-0.55);
\draw[leafe] (0,0.55) -- (0,-0.55);
\node[dot] at (5,0.55) {};
\node[ann] at (5,0.85) {123};
\node[dot] at (2.5,0.55) {};
\node[ann] at (2.5,0.85) {132};
\node[dot] at (5,-0.55) {};
\node[ann] at (5,-0.85) {213};
\node[dot] at (0,0.55) {};
\node[ann] at (0,0.85) {231};
\node[dot] at (2.5,-0.55) {};
\node[ann] at (2.5,-0.85) {312};
\node[dot] at (0,-0.55) {};
\node[ann] at (0,-0.85) {321};
\draw[box] (-0.42,-1.0) rectangle (0.42,1.0);
\node[ann] at (0,-1.52) {leaf 1};
\draw[box] (2.08,-1.0) rectangle (2.92,1.0);
\node[ann] at (2.5,-1.52) {leaf 2};
\draw[box] (4.58,-1.0) rectangle (5.42,1.0);
\node[ann] at (5,-1.52) {leaf 3};
\node[ann] at (2.5,-2.05) {the bottom layer: the three rank-one leaves};
\node[ann,anchor=west] at (6.0,0.55) {bold: the leaves' own edges};
\node[ann,anchor=west] at (6.0,0.15) {light: a perfect matching};
\node[ann,anchor=west] at (6.0,-0.2) {between every pair of leaves};
\end{tikzpicture}
\caption{A welding tower of rank two, drawn on $CT_3$. Its six vertices are the permutations of $\{1,2,3\}$, two joined when they differ by a transposition. Grouping them by their last symbol splits the graph into three parts of equal order, each inducing a single edge; those are the \emph{rank-one leaves}, the bottom layer of the tower, and they are what \textbf{C1} imposes its condition on. Between every pair of leaves the remaining edges form a perfect matching, which is what makes the graph a \emph{weld} of the three rather than their disjoint union. A rank-$r$ graph is the same picture with each leaf replaced by a graph of rank $r-1$, unfolded until the bottom layer is reached. In the proof of Theorem~3.1 the tower is $CT_a\square B$, all of whose rank-one leaves are copies of $B$, and \textbf{C1} asks of $B$ only that it be a single vertex or of even order, which is why the crown can be carried at all.}
\label{fig:weld}
\end{figure}
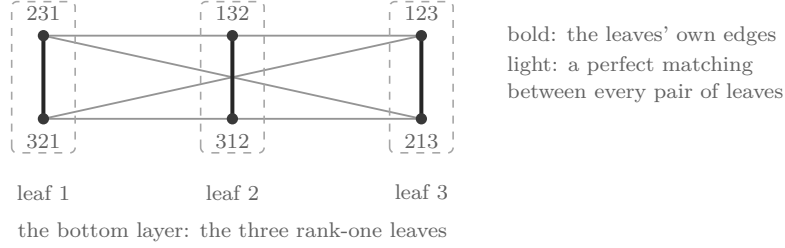

Its bottom layer consists of rank-one graphs, and those are its
\textbf{rank-one leaves} --- the term \textbf{C1} and Section 3.3 both
turn on, so it is worth fixing here. Figure \ref{fig:weld} draws a tower
of rank two with its leaves marked.

The leaf condition in \textbf{C1} is imposed over one such tower, not
over some tower. And Section 2's paired two-disjoint-path coverability
is the case \(k = 2\) of the following. For \(k ≥ 1\), a
\textbf{\(k\)-demand} in a bipartite graph is a family of \(k\) terminal
pairs whose \(2k\) terminals are distinct and each of which joins the
two color classes; a \textbf{paired \(k\)-disjoint path cover} for it is
a family of \(k\) pairwise vertex-disjoint paths, the \(i\)th joining
the \(i\)th pair, whose vertex sets together exhaust the graph.

\textbf{C1 (welding).} \emph{Let \(n ≥ 2\) and let \(H\) be a bipartite
transposition-like graph of rank \(n\) in which every rank-one leaf of
the welding is a single vertex or has even order. Then \(H\) admits a
paired \((n−1)\)-disjoint path cover for every \((n−1)\)-demand.}

\textbf{C2 (descent).} \emph{Let \(H\) be balanced bipartite with
\(2k ≤ |V(H)|\). If \(H\) admits a paired \(k\)-disjoint path cover for
every \(k\)-demand, then for every \(ℓ\) with \(1 ≤ ℓ ≤ k\) it admits a
paired \(ℓ\)-disjoint path cover for every \(ℓ\)-demand.}

\textbf{C1} and \textbf{C2} are Theorem 1.5 and Proposition 1.1(c) of
Coleman, Fischberg, Gong, Harrington and Wong {[}10{]}, quoted here as
the companion paper quotes them {[}4, Theorems 7.2 and 7.3{]}. The order
condition \(2k ≤ |V(H)|\) in \textbf{C2} says that \(k\)-subsets of each
color class exist. It is implicit in the source, whose proof extends an
\(ℓ\)-demand upward to a \(k\)-demand, and we state it because without
it the hypothesis can hold vacuously while the conclusion fails. It
holds at every use below.

\textbf{C3 (products of transposition graphs).} \emph{Every Cartesian
product of complete transposition graphs \(CT_a\) with \(a ≥ 2\), of
order at least four, is paired two-disjoint-path coverable.}

\textbf{C3} is the companion paper's {[}4{]} Theorem 7.1. Two further
results of theirs, \textbf{C4} the product lift and \textbf{C5} the
buffer lemma, are quoted where they are used --- \textbf{C4} in Section
4.2 before Corollary 4.3, \textbf{C5} in Section 6.4 --- since neither
is needed to reach Theorem 3.1 and neither is quoted exactly as printed.

\textbf{None of C1 through C5 is charged as an assumption by the formal
development.} Each is proved there: \textbf{C1} and \textbf{C2} were
formalized from the foundations rather than taken from the literature,
and \textbf{C3}, \textbf{C4} and \textbf{C5} are theorems of the
companion's development. Section 10.4 sets out the seven results this
paper does take on trust, and none of them is on this list.

\hypertarget{the-statement}{%
\subsection{3.3 The statement}\label{the-statement}}

Where this leaves us. Barrus's classification has reduced the whole
bipartite branch to Cartesian products of complete transposition graphs
together with at most one crown, and \textbf{C1} to \textbf{C3} with the
prism lemma upgrade those products from Hamiltonian to
Hamilton-laceable. The theorem below is that upgrade, stated for the
realization graph.

\textbf{Theorem 3.1.} \emph{If \(G(d)\) is triangle-free then \(G(d)\)
is maximally Hamiltonian. A triangle-free realization graph is
bipartite, so the assertion is that \(G(d)\) is Hamilton-laceable.}

\textbf{Lemma 3.2.} \emph{The crown graph \(K_{6,6} − 6K_2\) is paired
two-disjoint-path coverable, and hence Hamilton-laceable.}

Lemma 3.2 concerns a single graph on twelve vertices. All 3,600
admissible demands were certified exhaustively, and the certificate is
part of the verification suite described in Section 11. It is
machine-checked in the kernel, with no compiled-evaluation certificate.

\emph{Proof of Theorem 3.1.} Take the factorization supplied by the
classification. Every factor is a single vertex or a balanced bipartite
graph of even order: \(CT_1 = K_1\); for \(k ≥ 2\) the graph \(CT_k\)
has \(k!\) vertices and is bipartite by the parity of permutations; and
the crown has twelve, is balanced, and is paired two-disjoint-path
coverable by Lemma 3.2. A Cartesian product of balanced bipartite graphs
is balanced bipartite, so every partial product below inherits that
much. There are two shapes to treat, according to whether a crown
occurs.

If no crown occurs, the graph is a product of complete transposition
graphs. Every nontrivial factor has even order, so the product has order
one, order two, or order at least four. In the first two cases it is
\(K_1\) or \(K_2\), maximally Hamiltonian by convention; in the last,
\textbf{C3} gives the paired property directly and \textbf{C2} at
\(ℓ = 1\) turns it into Hamilton-laceability.

If a crown occurs, multiply the factors onto the crown one at a time,
maintaining that the partial product is balanced bipartite, of even
order at least four, and paired two-disjoint-path coverable. The crown
itself satisfies all three, by Lemma 3.2 and its twelve vertices. Let
\(B\) be the partial product so far and \(CT_a\) the next factor.

If \(a ≥ 3\), then \(CT_a □ B\) is a bipartite transposition-like graph
of rank \(a\) all of whose rank-one leaves are copies of \(B\): it is
the weld of the \(a\) cosets of \(CT_a\), each carrying a copy of
\(CT_{a−1} □ B\), matched by the transpositions moving the last symbol,
and the tower bottoms out at \(CT_1 □ B = B\). The leaf is of even
order, so \textbf{C1} applies and yields a paired \((a−1)\)-disjoint
path cover; \textbf{C2} descends it to \(k = 2\). This is the step that
needs a leaf with no structure beyond even order and laceability, and it
is the only reason the crown can be carried at all.

If \(a = 2\), then \(CT_2 = K_2\) has rank two and \textbf{C1} does not
reach it. Multiplying by it is taking the \textbf{prism} over \(B\), and
Lemma A.1 supplies exactly that case: the prism over a balanced
bipartite paired two-disjoint-path coverable graph on at least four
vertices is again one. Together the two cases carry every factor, and
\textbf{C2} at \(ℓ = 1\) gives Hamilton-laceability of the whole.
\(\blacksquare\)

Theorem 3.1 upgrades Barrus's corollary that every triangle-free
realization graph is Hamiltonian, from bare Hamiltonicity to maximal
Hamiltonicity. Two remarks. His corollary is printed without a proviso
for the one-vertex realization graph, and the convention that \(K_1\) is
maximally Hamiltonian is ours, fixed in Section 2. And the upgrade is
available because \textbf{C1} admits an arbitrary even-order
Hamilton-laceable rank-one leaf; the exceptional factor enters as one
such leaf, and nothing about it beyond Lemma 3.2 is used.

\hypertarget{what-section-4-takes-from-this}{%
\subsection{3.4 What Section 4 takes from
this}\label{what-section-4-takes-from-this}}

Corollary 4.3 propagates maximal Hamiltonicity along a Tyshkevich
composition, and for bipartite factors it needs more than the conclusion
of Theorem 3.1. It needs the \textbf{paired two-disjoint-path} property
itself, which is what the product argument above actually establishes
for every bipartite factor. We state that separately, because the
distinction is easy to lose and the product argument does not go through
with Hamilton-laceability alone.

\textbf{Corollary 3.3.} \emph{Every bipartite realization graph of order
at least two is paired two-disjoint-path coverable.}

\emph{Proof.} The proof of Theorem 3.1 establishes the paired property
in every shape it treats and then descends to Hamilton-laceability; what
that proof does not do is collect the three cases, since its conclusion
is the descended one. We collect them. By the classification, \(G(d)\)
is a Cartesian product of complete transposition graphs and at most one
crown, and a nontrivial factor has even order, so the product has order
one, order two, or order at least four.

\emph{Order two.} Then \(G(d) = K_2\), and the property holds vacuously:
a \(2\)-demand needs four distinct terminals and there are two vertices.
This is the same vacuity \textbf{C4} rests on for a \(K_2\) factor,
recorded in Section 4.2, and it is worth naming as vacuity rather than
passing it off as a case that was checked.

\emph{No crown, order at least four.} \textbf{C3} gives the paired
property directly, for every product of complete transposition graphs of
that order. The descent by \textbf{C2} at \(ℓ = 1\) in the proof of
Theorem 3.1 is a consequence of the paired property, not a substitute
for it.

\emph{A crown occurs.} The induction in the proof of Theorem 3.1 carries
``balanced bipartite, of even order at least four, and paired
two-disjoint-path coverable'' as an invariant of every partial product:
it holds of the crown by Lemma 3.2, is preserved by \textbf{C1} and
\textbf{C2} when the next factor is \(CT_a\) with \(a ≥ 3\), and by
Lemma A.1 when it is \(CT_2 = K_2\). The last partial product is
\(G(d)\). \(\blacksquare\)

\hypertarget{reductions-and-the-shape-of-the-induction}{%
\section{4. Reductions, and the shape of the
induction}\label{reductions-and-the-shape-of-the-induction}}

Theorem 3.1 settles the bipartite case outright. This section removes
three further cases, so that the rest of the paper may assume a degree
function of one specific kind. Everything here is short, and two of the
three steps are quoted from the literature.

\hypertarget{inactive-vertices}{%
\subsection{4.1 Inactive vertices}\label{inactive-vertices}}

A ground vertex is \textbf{active} when at least two sets occur as its
neighborhood across the realizations of \(d\), and \textbf{inactive}
otherwise: an inactive vertex has the same neighborhood in every
realization. We also call \textbf{\(d\) itself active} when every one of
its ground vertices is. Activity is \textbf{not} a hypothesis of
Sections 6, 8 or 10 --- Section 4.3 explains why it need not be --- but
the notion is used there, in the form of an active line of a matrix.

\textbf{Theorem 4.1.} \emph{Let \(d\) be graphical. If \(v\) is inactive
for \(d\), then \(G(d) ≅ G(d^*)\), where \(d^*\) is the degree function
obtained by deleting \(v\) and lowering the degree of each of its
neighbors by one.}

\emph{Proof.} Because \(v\) is inactive there is a single set \(S\) with
\(N_G(v) = S\) for \textbf{every} realization \(G\) of \(d\). A 2-switch
acts on four distinct vertices, so one using an edge at \(v\) would
leave \(v\) with one neighbor exchanged for a different one, changing
\(N_G(v)\). No 2-switch of \(G(d)\) therefore touches \(v\), and every
switch lies inside \(V ∖ \{v\}\). Deleting \(v\) therefore sends
realizations of \(d\) to realizations of \(d^*\) bijectively --- the
inverse adds \(v\) back with neighborhood \(S\) --- and carries
2-switches to 2-switches in both directions. That is an isomorphism of
realization graphs. \(\blacksquare\)

The theorem is not a branch of the induction of Section 10, which never
needs to dispatch on activity; it is recorded because it is the
one-fiber instance of Lemma 5.1 and because the notion it turns on is
used later.

Activity has an equivalent local description, which we do not use: a
vertex is active exactly when it lies in an induced \(2K_2\), \(P_4\) or
\(C_4\) of some realization. Schvöllner and Pastine {[}32{]} work with
that form, and their switch-degree formula quantifies it; we record the
agreement in Section 6.2 and take the neighborhood form as the
definition, since it is what the later sections consume.

\hypertarget{tyshkevich-composition}{%
\subsection{4.2 Tyshkevich composition}\label{tyshkevich-composition}}

The composition takes two ingredients of different kinds. The left one,
\(d_1\), is a \textbf{split} sequence together with a chosen partition
of its ground into a clique side \(A\) and an independent side \(B\).
The right one, \(d_2\), is any degree sequence at all, on a ground
\(Y\). To form \(d_1 ∘ d_2\), take a realization of each and join
\textbf{every vertex of \(A\) to every vertex of \(Y\)}: the clique side
goes complete to the whole of the other factor, and \(B\) is left
untouched. In degrees, a vertex of \(A\) gains \(|Y|\), a vertex of
\(Y\) gains \(|A|\), and a vertex of \(B\) is unchanged, so the
operation is visible in the sequence itself.

Every degree sequence has a canonical decomposition into indecomposable
parts under this operation {[}34{]}. What the proof below uses is not
the canonicity but the fact underneath it: whether \(d\) decomposes, and
how, is a property of the sequence rather than of any one realization.

To make the decomposable branch inductive we need more than that the
sequence splits: we need the splitting to separate the realizations and
their 2-switches into independent coordinates, so that a walk in each
factor can be run without disturbing the other. The next theorem turns
that separation into a Cartesian product, which is what makes the
factors usable as an induction hypothesis.

\textbf{Theorem 4.2.} \emph{If \(d = d_1 ∘ d_2\) is a nontrivial
Tyshkevich composition, then \(G(d) ≅ G(d_1) □ G(d_2)\).}

\emph{Proof.} Write \(A\) and \(B\) for the clique and independent sides
of the splitted factor and \(Y\) for the ground of \(d_2\). The cross
adjacencies are the same in \textbf{every} realization of \(d\), not
merely in the one the composition is built from: that is the content of
Tyshkevich's decomposition theorem {[}34{]}, which makes the canonical
decomposition a property of the degree sequence rather than of a chosen
realization. It can also be seen by counting, and the count is worth
doing in the aggregate rather than one vertex at a time, because a
vertex-by-vertex reading has to presuppose the very allocation it is
trying to establish. Write \(a = |A|\) and \(y = |Y|\). Summing the
degrees of the composition over \(A\) and over \(B\), the \(A\)--\(B\)
adjacencies appear once in each sum and cancel:

\[
Σ_{u ∈ A} d(u) − Σ_{u ∈ B} d(u) = a(a − 1) + a y.
\]

Now take \textbf{any} realization \(H\) of that same labeled degree
function. The same cancellation leaves

\[
\begin{gathered}
Σ_{u ∈ A} d(u) − Σ_{u ∈ B} d(u) \\
= 2|E_H(A)| + |E_H(A, Y)| − 2|E_H(B)| − |E_H(B, Y)| \\
≤ a(a − 1) + a y,
\end{gathered}
\]

since \(|E_H(A)| ≤ a(a−1)/2\) and \(|E_H(A,Y)| ≤ a y\) and the two
subtracted terms are nonnegative. The degree identity forces equality,
and equality in a sum of nonnegative deficits forces it in each: \(A\)
is a clique, \(A\)--\(Y\) is complete, \(B\) is independent, and
\(B\)--\(Y\) is empty, in every realization of \(d\). Setting \(Y = ∅\)
gives the fact Section 6.4 uses, that a split partition in one
realization is a split partition in every realization of the same degree
function.

Given that, a 2-switch alters four adjacencies among two disjoint pairs,
so a switch using a cross adjacency would change one, and no switch
crosses. Every switch therefore lies inside one part. Conversely a
switch inside either factor lifts: it leaves the cross adjacencies and
the other factor untouched, so it is a switch of \(d\). Realizations of
\(d\) thus correspond to pairs of realizations of \(d_1\) and \(d_2\),
and adjacency is ``equal in one coordinate and switch-adjacent in the
other'' --- the Cartesian product. \(\blacksquare\)

\textbf{C4 (the product lift).} \emph{Let \(A\) and \(B\) be graphs
whose Cartesian product \(A □ B\) is non-bipartite. Suppose each of
\(A\) and \(B\) is either non-bipartite and Hamilton-connected, or
carries a proper surjective 2-coloring for which it is both
Hamilton-laceable and paired two-disjoint-path coverable. Then \(A □ B\)
is Hamilton-connected.}

The companion's Section 4 proves \textbf{C4} in three cases and numbers
only the middle one. Both factors non-bipartite: both are then
Hamilton-connected, and no parity device is needed. One factor bipartite
of order at least four: this is its Proposition 4.1, the substantive
case, and the paired property is what it consumes --- a bipartite factor
cannot be Hamilton-connected, so the product's Hamilton-connectedness
has to be manufactured from a factor that cannot have it, by doubling a
layer and covering it with two disjoint paths. One factor \(K_2\): the
product is then a prism over the other factor and is treated directly. A
factor of order two satisfies the paired hypothesis vacuously, there
being no four distinct terminals, which is why the statement above needs
no case distinction.

\textbf{Corollary 4.3 (component closure).} \emph{Let \(d = d_1 ∘ d_2\)
be a nontrivial composition. If each of \(d_1\), \(d_2\) has a maximally
Hamiltonian realization graph, then so does \(G(d)\).}

\emph{Proof.} If either factor has a single realization its realization
graph is \(K_1\), the product is the other factor's realization graph,
and there is nothing to prove; so assume both have order at least two.
By Theorem 4.2, \(G(d) ≅ G(d_1) □ G(d_2)\).

If both factors are bipartite then so is the product: give \((x, y)\)
the color \(c_1(x) + c_2(y)\) reduced modulo two, where \(c_i\) is the
coloring of the \(i\)th factor; a product edge moves exactly one
coordinate, so it changes exactly one summand and hence the sum. So
\(G(d)\) is bipartite, hence triangle-free, and \textbf{Theorem 3.1
gives the conclusion outright}. This is the one case \textbf{C4} does
not reach, since it assumes its product non-bipartite; Section 3 already
covers it, and there is no need to build the product path by hand.

Otherwise some factor is non-bipartite. A layer of the product is a copy
of that factor, so \(G(d)\) inherits its odd cycle and is non-bipartite,
and \textbf{C4} applies to the two factors. What \textbf{C4} asks of a
non-bipartite factor is Hamilton-connectedness, which is what maximal
Hamiltonicity gives there. What it asks of a bipartite factor is a
proper surjective 2-coloring, Hamilton-laceability for it, and the
paired two-disjoint-path property. Realization graphs are connected and
the factor has order at least two, so its bipartition is proper and
surjective; laceability is again maximal Hamiltonicity; and the paired
property is Corollary 3.3. Hence \(G(d)\) is Hamilton-connected.
\(\blacksquare\)

When one factor is bipartite and the other is not, it is the paired
property and not laceability that crosses the product --- Section 2
records why --- so there the corollary consumes Corollary 3.3. When both
factors are bipartite it consumes Theorem 3.1 directly instead. It is
stated for a single composition into two factors, which is all Section
10.1 asks of it; iterating it needs nothing further, since the product
is again the realization graph of the composed sequence and Corollary
3.3 applies to that as it stands.

Corollary 4.3 propagates a conclusion; it supplies none for an
indecomposable factor. That is what the rest of the paper is for.

\hypertarget{the-induction}{%
\subsection{4.3 The induction}\label{the-induction}}

Let \(d\) be graphical with at least two realizations. One of four
branches applies: \(G(d)\) is bipartite, and Theorem 3.1 gives the
conclusion; \(d\) is Tyshkevich-decomposable, and Corollary 4.3
propagates from the factors; \(G(d)\) is the \(K_3\) base graph, which
is Hamilton-connected directly; or \(d\) is indecomposable and
non-bipartite with \(G(d)\) not \(K_3\).

\textbf{Activity is not among the branches, and Section 4.1 is not one
of them.} Under the standing hypotheses of the last branch it comes for
free: the only indecomposable graph carrying an inactive vertex is the
one-vertex graph {[}32{]}, whose realization graph is a single vertex
and therefore bipartite, so it has already been taken by the first
branch. Theorem 4.1 remains true and is worth having --- it is the
one-fiber instance of Lemma 5.1 --- but the induction never needs to
dispatch on it.

The second branch replaces \(d\) by strictly smaller degree functions,
so an induction on the number of ground vertices reaches the fourth.
Sections 5 to 10 treat it, and we call it the \textbf{main case}.
Section 10.1 carries the dispatch out as the case split of the induction
that proves Theorem 10.1, which is the only place it is performed.

We record what the main case gives, since every later section assumes
it. Every ground vertex is active, which is to say it has at least two
realizable neighborhoods, so every quotient has at least two vertices.
\(d\) is Tyshkevich-indecomposable. \(G(d)\) contains an odd cycle, and
it has more than three vertices.

\hypertarget{the-fibers-and-the-quotient}{%
\section{5. The fibers and the
quotient}\label{the-fibers-and-the-quotient}}

Fix a graphical degree function \(d\) and a ground vertex \(v\), which
we call the \textbf{pivot}: it is the vertex whose neighborhood
everything below is indexed by, and it stays fixed for the whole
section.

Once \(v\) is fixed, every 2-switch is of one of two kinds, and the
structure of this section is what follows from that. A 2-switch that
does not use \(v\) leaves \(N(v)\) alone, so if the realizations are
grouped by their value of \(N(v)\), those moves stay inside a group: the
groups are the \textbf{fibers}, and the switches internal to one are its
edges. A 2-switch that does use \(v\) replaces exactly one neighbor of
\(v\) by one other, so it carries a realization to a group whose label
differs in a single element. Those are the moves between fibers, and
recording them is what the \textbf{quotient} does. Nothing else can
happen, so \(G(d)\) is the fibers together with the edges of the second
kind, and the rest of the section is this picture made precise.

What it sets up for the rest of the proof: the fibers over \(v\), the
family of neighborhoods that occur, the quotient graph on that family,
and the \textbf{interface} between two adjacent fibers --- the edges of
\(G(d)\) that run from one fiber to the next, which is what the one-pass
construction crosses on. Everything here is proved from first principles
or from classical results about degree sequences.

\textbf{Notation, in one place.} Each item is defined again where it is
first used; this is for looking back, not for reading through.

Every letter is chosen to be read rather than memorized, and the gloss
says how.

\begin{longtable}[]{@{}
  >{\raggedright\arraybackslash}p{(\columnwidth - 2\tabcolsep) * \real{0.2000}}
  >{\raggedright\arraybackslash}p{(\columnwidth - 2\tabcolsep) * \real{0.8000}}@{}}
\toprule\noalign{}
\endhead
\bottomrule\noalign{}
\endlastfoot
\(v\) & the \textbf{pivot}, fixed for the whole section; the
\textbf{v}ertex everything is indexed by \\
\(Φ_S\) & the \textbf{fiber} over \(S\), the realizations of \(d\) with
\(N(v) = S\); \textbf{Φ} for fiber \\
\(I_v\) & the family of \(S\) for which \(Φ_S\) is nonempty: the
\textbf{i}ndex set of the fibers, one \(S\) per fiber \\
\(f_S\) & the \textbf{residual} degree function \(d(t) − [t ∈ S]\) on
the ground; lowercase \(d\), being what is left of \(d\) \\
\(F = G − v\) & the \textbf{ground realization} of \(G ∈ Φ_S\), which
realizes \(f_S\); uppercase for a graph, as \(G\) is, with the pivot
gone \\
\(S\) & the \textbf{s}ource neighborhood, the fiber a step starts in \\
\(T = S − b + a\) & the \textbf{t}arget, differing from \(S\) in one
element: \(b\) is the neighbor \(v\) gives up, \(a\) the one it takes
on \\
\(δ = f(a) − f(b)\) & the degree difference across that interface \\
a \textbf{witness} & a ground vertex \(c\) adjacent to \(a\) but not to
\(b\), the vertex a connector swaps at \\
\(W_{ab}\), \(w_{ab}\) & the set \(N_F(a) ∖ (N_F(b) ∪ \{b\})\) of
\textbf{w}itnesses, from \(a\) toward \(b\), and its size \\
\(e_u\) & a unit at \(u\): the degree function that is \(1\) at \(u\)
and \(0\) elsewhere, so \(f − e_a + e_b\) lowers \(a\) by one and raises
\(b\) \\
\end{longtable}

Two conventions are worth stating with them. All of \(W\), \(w\) and the
witness condition are read in the \textbf{ground} realization \(F\),
never in \(G\). And \(W_{ba}\), \(w_{ba}\) mean the same formulas with
\(a\) and \(b\) interchanged, which is a different quantity in the same
graph, not the same quantity read backwards.

\hypertarget{fibers}{%
\subsection{5.1 Fibers}\label{fibers}}

Every realization of \(d\) determines the neighborhood of the pivot, so
\(G ↦ N_G(v)\) is a map from the realizations of \(d\) to the subsets of
\(V ∖ \{v\}\), and this section is about its fibers. For
\(S ⊆ V ∖ \{v\}\) write

\[
Φ_S = \{ G : G \text{ realizes } d \text{ and } N_G(v) = S \},
\]

the \textbf{fiber} over \(S\) --- \(Φ\) for fiber. In words: a fiber is
all the realizations in which \(v\) has one particular set of neighbors.
Its members are graphs, not sets; the set \(S\) is only the label they
share. The fibers partition the realizations of \(d\), one class for
each neighborhood that occurs. The \textbf{residual degree function} is

\[
f_S(t) = d(t) − [t ∈ S], \quad \text{ for } t ≠ v,
\]

which is what remains after deleting \(v\) and paying for the edges it
used.

Inside a fiber the pivot is dead weight. Its neighbors are the same in
every member, so by the classification above every internal edge of the
fiber is a move of the first kind and none of them touches \(v\).
Deleting it therefore throws away nothing the fiber records, and what is
left is the same problem on one fewer vertex.

\textbf{Lemma 5.1.} \emph{Deleting \(v\) is a bijection from \(Φ_S\)
onto the set of realizations of \(f_S\), and it carries the 2-switches
inside \(Φ_S\) to the 2-switches of that smaller realization graph. So
\(Φ_S ≅ G(f_S)\) as graphs, and \(f_S\) has one fewer ground vertex than
\(d\).}

A fiber is named by its set of realizations and treated as a graph
throughout: \(Φ_S\) carries the subgraph of \(G(d)\) induced on that
set, and that is the graph Lemma 5.1 identifies.

\emph{Proof.} A realization \(G ∈ Φ_S\) has \(N_G(v) = S\), so \(G − v\)
realizes \(f_S\), and conversely adjoining a vertex joined exactly to
\(S\) inverts the map; the two constructions are mutually inverse, so
deletion is a bijection. The 2-switches inside \(Φ_S\) are exactly the
moves of the first kind, those not using \(v\), so deletion leaves all
four of their toggled edges intact and each survives as a 2-switch of
\(f_S\); conversely a 2-switch of \(f_S\) lifts by re-adjoining \(v\).
\(\blacksquare\)

Lemma 5.1 is what makes the induction of Section 10 well founded: every
fiber is a strictly smaller instance of the same problem.

Call \(S\) \textbf{realizable} when \(Φ_S\) is nonempty, and let \(I_v\)
denote the family of realizable neighborhoods. The next statement
identifies that family without constructing a single realization.

\textbf{Lemma 5.2.} \emph{\(S\) is realizable if and only if
\(|S| = d(v)\) and \(f_S\) is graphical.}

\emph{Proof.} Both directions are elementary. Deleting \(v\) from a
realization proves one; adding a new vertex joined to \(S\) to a
realization of \(f_S\) proves the other. \(\blacksquare\) We use Lemma
5.2 constantly, since it turns a question about realizations into a
question about degree sequences.

\textbf{A running example, carried through the paper.} Take
\(d = (3,2,2,2,1)\) on \(\{1,…,5\}\) and pivot at \(v = 2\), the
sequence Section 5.5 returns to and the one drawn in Figure 2 of Hladík
and Fink {[}18{]}. It is the smallest sequence exhibiting what this
section is about, and there is no other of its order. There are six
realizations, and the neighborhoods that occur at \(v\) are

\[
I_v = \{ \{1,3\}, \{1,4\}, \{1,5\}, \{3,4\} \},
\]

To see where that comes from, ask each realization who \(v\) is joined
to in it, and group by the answer. All six, writing \(ij\) for the edge
\(\{i,j\}\) and setting the pivot's own edges in bold:

\begin{longtable}[]{@{}lll@{}}
\toprule\noalign{}
& edges of the realization & \(N(2)\) \\
\midrule\noalign{}
\endhead
\bottomrule\noalign{}
\endlastfoot
\(G_1\) & \textbf{12}, 13, 14, \textbf{23}, 45 & \(\{1,3\}\) \\
\(G_5\) & \textbf{12}, 14, 15, \textbf{23}, 34 & \(\{1,3\}\) \\
\(G_2\) & \textbf{12}, 13, 14, \textbf{24}, 35 & \(\{1,4\}\) \\
\(G_4\) & \textbf{12}, 13, 15, \textbf{24}, 34 & \(\{1,4\}\) \\
\(G_3\) & \textbf{12}, 13, 14, \textbf{25}, 34 & \(\{1,5\}\) \\
\(G_6\) & 13, 14, 15, \textbf{23}, \textbf{24} & \(\{3,4\}\) \\
\end{longtable}

\(G_1\) and \(G_5\) agree nowhere else and agree at the pivot, so they
lie in one fiber; \(G_2\) does not. And \(G_1\) and \(G_5\) differ by
the 2-switch on the alternating four-cycle \(1-3-4-5\), which removes
\(13\) and \(45\) and adds \(34\) and \(15\): a move of the first kind,
leaving \(v\) untouched, which is why the two sit in one fiber and why
Lemma 5.1 can delete \(v\) without losing that edge. Grouping all six
this way gives four fibers: two of size two, over \(\{1,3\}\) and
\(\{1,4\}\), and two singletons, over \(\{1,5\}\) and \(\{3,4\}\). Any
two of these neighborhoods differ in a single element except that last
pair, which differs in all four. Figure \ref{fig:realizationgraph}
already drew \(G(d)\) whole; Figure \ref{fig:fibers} cuts the same six
realizations into these four fibers, once Section 5.4 supplies the last
piece to label.

The example is worth keeping in view, because one small family carries
most of the paper. \(I_v\) is the smallest \textbf{Y-family} of Section
7.5, so this quotient is the exception of Theorem 7.7: \(K_4\) with the
edge between \(\{1,5\}\) and \(\{3,4\}\) missing, and no Hamilton path
between \(\{1,3\}\) and \(\{1,4\}\). Those last two are the pair the
theorem forbids, and the obstruction is plain in so small a graph: a
Hamilton path joining them would have to visit \(\{1,5\}\) and
\(\{3,4\}\) one after the other, and that is the one pair this quotient
leaves unjoined. Section 5.5 reads the same family the other way, as a
quotient with two maximal members and hence no greatest one. And the
failure is confined to the quotient: \(G(d)\) here is \(K_6 − 3K_2\),
which is Hamilton-connected. That a quotient can fail while its
realization graph does not is the reason Sections 6 and 8 exist.

\hypertarget{the-family-of-realizable-neighborhoods-is-shifted}{%
\subsection{5.2 The family of realizable neighborhoods is
shifted}\label{the-family-of-realizable-neighborhoods-is-shifted}}

Order the ground vertices other than \(v\) as
\(u_1 ≻ u_2 ≻ ⋯ ≻ u_{n−1}\) by non-increasing degree, with ties broken
by label. Read a subset of \(V ∖ \{v\}\) as a subset of \([n−1]\)
through this order. The Gale order and the notion of a shifted family
are those of Section 7, and we state them here because everything from
this subsection on turns on them. For \(k\)-sets
\(X = \{x_1 < ⋯ < x_k\}\) and \(Y = \{y_1 < ⋯ < y_k\}\), write \(X ≼ Y\)
when \(x_i ≤ y_i\) for every \(i\); that is the \textbf{Gale order}. A
family of \(k\)-sets is \textbf{shifted} when it is closed under a
single downward step of that order --- replacing a member of a set by a
smaller element not already in it keeps you inside the family. Single
steps generate the order, so a shifted family is the same thing as a
nonempty down-set of \(≼\), and Section 7 uses both descriptions. The
same objects are called \emph{Gale ideals} and \emph{left-compressed} or
\emph{stable} families elsewhere; Section 7.9 gives seven names for the
special case in which the family is the basis family of a matroid.

\textbf{Theorem 5.3.} \emph{\(I_v\) is a shifted family of
\(d(v)\)-subsets of \([n−1]\).}

\emph{Proof.} Let \(S ∈ I_v\), let \(j ∈ S\) and \(i ∉ S\) with \(i\)
earlier than \(j\) in the order, and put \(S' = S − j + i\). If no such
pair exists, which is the case exactly when \(S\) is Gale-least, the
condition holds vacuously at \(S\) and there is nothing to prove.
Earlier in the order means \(d(i) ≥ d(j)\). Compare the residuals. They
agree away from \(i\) and \(j\), and

\[
\begin{gathered}
f_S(i) = d(i), \quad f_S(j) = d(j) − 1, \\
f_{S'}(i) = d(i) − 1, \quad f_{S'}(j) = d(j).
\end{gathered}
\]

So \(f_{S'}\) is obtained from \(f_S\) by moving one unit from the entry
at \(i\) to the entry at \(j\), and that transfer runs from a strictly
larger entry to a strictly smaller one, since
\(f_S(i) = d(i) ≥ d(j) > d(j) − 1 = f_S(j)\).

A transfer of this shape preserves graphicality, by one edge swap. Let
\(G\) realize \(f_S\). Since \(\deg_G(i) > \deg_G(j)\), some vertex
\(w\) outside \(\{i, j\}\) is adjacent to \(i\) and not to \(j\):
otherwise every neighbor of \(i\) other than \(j\) would also be a
neighbor of \(j\), and as \(i\) and \(j\) are either adjacent to each
other or not, both degrees would be reduced by the same amount on
discarding that edge, giving \(\deg_G(i) ≤ \deg_G(j)\). Deleting \(iw\)
and inserting \(jw\) is therefore legal, lowers the degree at \(i\) by
one, raises it at \(j\) by one, and leaves every other degree unchanged,
so the result realizes \(f_{S'}\).

Since \(f_S\) is graphical by Lemma 5.2, so is \(f_{S'}\), and Lemma 5.2
returns \(S' ∈ I_v\). \(\blacksquare\)

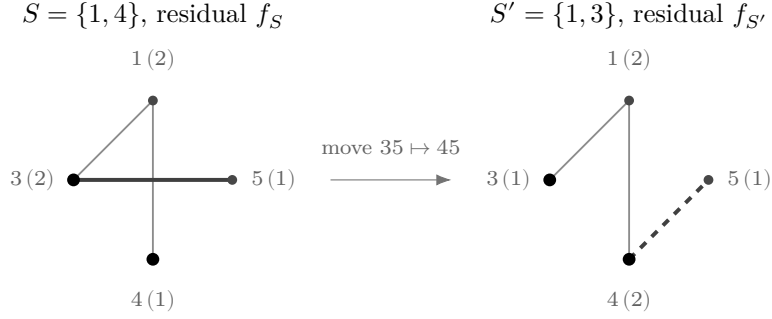
\begin{figure}[tbp]
\centering
\begin{tikzpicture}[x=1cm,y=1cm,
  dot/.style={circle,fill=black!72,inner sep=1.3pt},
  hot/.style={circle,fill=black,inner sep=1.7pt},
  gedge/.style={draw=black!45,line width=0.7pt},
  del/.style={draw=black!75,line width=1.5pt},
  add/.style={draw=black!75,line width=1.5pt,dashed},
  ann/.style={font=\scriptsize,text=black!62},
  ttl/.style={font=\small}]
\node[ttl] at (0,2.18) {$S=\{1,4\}$, residual $f_S$};
\draw[gedge] (6.429e-17,1.05) -- (-1.05,1.286e-16);
\draw[gedge] (6.429e-17,1.05) -- (-1.929e-16,-1.05);
\draw[del] (-1.05,1.286e-16) -- (1.05,-2.572e-16);
\node[dot] at (6.429e-17,1.05) {};
\node[ann] at (9.773e-17,1.596) {$1\,(2)$};
\node[hot] at (-1.05,1.286e-16) {};
\node[ann] at (-1.596,1.955e-16) {$3\,(2)$};
\node[hot] at (-1.929e-16,-1.05) {};
\node[ann] at (-2.932e-16,-1.596) {$4\,(1)$};
\node[dot] at (1.05,-2.572e-16) {};
\node[ann] at (1.596,-3.909e-16) {$5\,(1)$};
\node[ttl] at (6.3,2.18) {$S'=\{1,3\}$, residual $f_{S'}$};
\draw[gedge] (6.3,1.05) -- (5.25,1.286e-16);
\draw[gedge] (6.3,1.05) -- (6.3,-1.05);
\draw[add] (6.3,-1.05) -- (7.35,-2.572e-16);
\node[dot] at (6.3,1.05) {};
\node[ann] at (6.3,1.596) {$1\,(2)$};
\node[hot] at (5.25,1.286e-16) {};
\node[ann] at (4.704,1.955e-16) {$3\,(1)$};
\node[hot] at (6.3,-1.05) {};
\node[ann] at (6.3,-1.596) {$4\,(2)$};
\node[dot] at (7.35,-2.572e-16) {};
\node[ann] at (7.896,-3.909e-16) {$5\,(1)$};
\draw[-{Latex[length=2mm]},black!55] (2.35,0) -- (3.95,0);
\node[ann,align=center] at (3.15,0.44) {move $35\mapsto 45$};
\end{tikzpicture}
\caption{One downward step of the Gale order, on the running example $d=(3,2,2,2,1)$ with pivot $v=2$, drawn on the ground $V\setminus\{v\}$. Left: a realization of the residual $f_S$ for $S=\{1,4\}$. The step replaces $j=4$ by the earlier $i=3$, and each vertex is labelled with its residual degree. What is easy to read past is that the set step and the residual step go in OPPOSITE directions: $4$ enters the residual and $3$ leaves it, so one unit moves off the entry at $i=3$, which is the larger, onto the entry at $j=4$. Right: the same realization after the single edge move. The witness $w=5$ is adjacent to $i$ and not to $j$, so deleting $i w$ (bold) and inserting $j w$ (dashed) is legal, lowers the degree at $i$ by one and raises it at $j$ by one, and realizes $f_{S'}$. Theorem~5.3 is this move, made once for every downward step.}
\label{fig:shift}
\end{figure}

Figure \ref{fig:shift} carries out one such step on the running example.
The direction is the part worth pausing on: the step replaces \(j\) by
the earlier \(i\) in the \emph{set}, and moves one unit off the entry at
\(i\) onto the entry at \(j\) in the \emph{residual}.

The argument uses no property of degree sequences beyond Lemma 5.2 and
that one swap. In particular it does not invoke the classical fact, due
to Ruch and Gutman {[}28{]}, that graphical sequences of a given sum
form a down-set in the \textbf{majorization order}, where one sequence
majorizes another when its sorted partial sums dominate the other's term
by term: a single unit moved between two coordinates is the one case of
that fact which can be exhibited directly, and it is the only case
Theorem 5.3 needs.

\hypertarget{the-quotient}{%
\subsection{5.3 The quotient}\label{the-quotient}}

The \textbf{quotient at \(v\)} is the graph on \(I_v\) in which \(S\)
and \(T\) are adjacent when some member of \(Φ_S\) is 2-switch adjacent
to some member of \(Φ_T\). Adjacent fibers differ in one element, which
is the second kind of move from the opening of the section: in a
2-switch that uses \(v\), the pivot lies in exactly one of the two
removed edges and one of the two added ones, since each pair is
disjoint, so the switch replaces a single neighbor of \(v\) by a single
other. The next theorem says the converse holds, so the quotient is not
merely contained in the Johnson graph of \(I_v\) but equal to it.

\textbf{What the theorem does.} Suppose the two fibers are joined at
all, and ask what such an edge must look like. It is a 2-switch of
\(G(d)\) carrying a realization with \(N(v) = S\) to one with
\(N(v) = T\), so it has to delete \(vb\) and add \(va\), and by the
classification at the head of the section its other two edges lie in the
ground and swap the same partner from \(a\) to \(b\). That pins the
shape completely and leaves one free choice, the ground vertex being
swapped. The theorem counts the legal choices, and the count comes out
in terms of the two ground degrees \(f(a)\) and \(f(b)\) in the source
realization. These edges are the \textbf{connectors} of Section 5.4,
counted there in this same way.

The two clauses below look asymmetric, and the reason is worth having in
advance: they count in the same way but in different degree functions.
The source side is read in \(f = f_S\); the target side is read in
\(f_T = f − e_a + e_b\), which lowers \(a\) by one and raises \(b\) by
one. So a difference read at the target exceeds the negated source
difference by exactly two, which is where every \(+2\) and every
\(2 − δ\) in this subsection and the next comes from. There is nothing
special about the constant; it is the one unit of degree that the step
moves, counted at both ends.

\textbf{Theorem 5.4 (quotient adjacency).} \emph{Let \(S ∈ I_v\) and
\(T = S − b + a ∈ I_v\). Then some member of \(Φ_S\) is adjacent to some
member of \(Φ_T\). Writing \(f = f_S\), if \(f(a) > f(b)\) then every
member of \(Φ_S\) has such an edge, with at least \(f(a) − f(b)\)
choices; if \(f(a) ≤ f(b) + 1\) then every member of \(Φ_T\) does, with
at least \(f(b) − f(a) + 2\) choices.}

\emph{Proof.} Write \(f = f_S\) and \(g = f_T\), both graphical. They
agree except that \(g(a) = f(a) − 1\) and \(g(b) = f(b) + 1\). The
\(Φ_S\) side supplies a cross edge when some member \(G ∈ Φ_S\) admits
the 2-switch that removes \(vb\) and \(ca\) and adds \(va\) and \(cb\),
for some \(c ∉ \{v,a,b\}\) with \(ca ∈ E(G)\) and \(cb ∉ E(G)\); the
edge is that 2-switch, joining \(G\) to its image. Writing
\(F = G − v\), which realizes \(f\) and carries every edge at \(a\)
because \(a ∉ S\), the condition reads: \(c ∈ N_F(a) ∖ \{b\}\) with
\(cb ∉ E(F)\).

Suppose \(f(a) > f(b)\) and let \(F\) realize \(f\). If no witness
existed, every \(c ∈ N_F(a) ∖ \{b\}\) would lie in \(N_F(b) ∖ \{a\}\);
comparing cardinalities gives \(f(a) − [ab ∈ F] ≤ f(b) − [ab ∈ F]\),
against \(f(a) > f(b)\). Discarding the shared term, the witnesses
number at least \(f(a) − f(b)\), and distinct witnesses give distinct
results.

Suppose instead \(f(a) ≤ f(b) + 1\). Then
\(g(b) = f(b) + 1 ≥ f(a) > f(a) − 1 = g(a)\), and the same count applies
in any realization of \(g\) with the roles of \(a\) and \(b\) exchanged,
giving at least \(g(b) − g(a) = f(b) − f(a) + 2\) witnesses whose
switches land in \(Φ_S\). The two cases overlap at \(f(a) = f(b) + 1\)
and together cover all values. \(\blacksquare\)

\textbf{Corollary 5.5.} \emph{The quotient at \(v\) is the Johnson graph
\(J(I_v)\) of Section 7.}

\emph{Proof.} Theorem 5.3 makes \(I_v\) a shifted family, so \(J(I_v)\)
is defined. Adjacent fibers have neighborhoods differing in one element,
so the quotient is a subgraph of \(J(I_v)\); Theorem 5.4 gives the
converse, that every such pair of neighborhoods is realized by a
2-switch between the fibers. \(\blacksquare\)

Corollary 5.5 is what licenses Sections 7 and 8 to speak about the
quotient. Every statement proved there about \(J(F)\) for a shifted
family \(F\) applies verbatim to the quotient at any ground vertex, with
\(F = I_v\).

Theorem 5.3 and Corollary 5.5 are what let Section 7 speak only of
shifted families: \(I_v\) goes in and a Hamilton path between two of its
members comes out, with the degree sequence left behind.

\hypertarget{the-interface}{%
\subsection{5.4 The interface}\label{the-interface}}

\textbf{What this subsection is for.} Section 10 walks a path in the
quotient and has to get from one fiber to the next at every step. The
edges available for that are the \textbf{connectors}, and this
subsection is entirely about counting them. One identity does the work:
the number of connectors out of a realization, minus the same count with
\(a\) and \(b\) interchanged, equals a difference of two degrees. That
pins the \textbf{difference} at every realization of the fiber, and it
is free. It does not pin either count on its own, and converting the
difference into an actual number is what the rest of the subsection
buys, at the price of assuming the source fiber is triangle-free. From
that single fact come the two things Section 9 uses --- a lower bound on
how many connectors each realization has, and, when the fiber is
bipartite, the exact statement that the connector graph expands by a
constant factor. The rest is the case analysis that establishes the
identity's consequences when the source fiber is triangle-free.

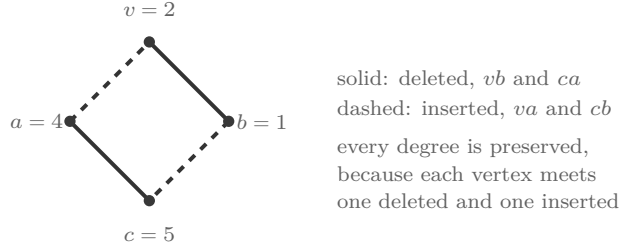
\begin{figure}[tbp]
\centering
\begin{tikzpicture}[x=1cm,y=1cm,
  dot/.style={circle,fill=black!78,inner sep=1.5pt},
  del/.style={draw=black!80,line width=1.5pt},
  add/.style={draw=black!80,line width=1.5pt,dashed},
  ann/.style={font=\scriptsize,text=black!62},
  ttl/.style={font=\small}]
\node[ttl] at (0,2.15) {the four toggled edges alone};
\draw[del] (0,1.05) -- (1.05,0);
\draw[add] (1.05,0) -- (0,-1.05);
\draw[del] (0,-1.05) -- (-1.05,0);
\draw[add] (-1.05,0) -- (0,1.05);
\node[dot] at (0,1.05) {};
\node[ann] at (0,1.491) {$v=2$};
\node[dot] at (1.05,0) {};
\node[ann] at (1.491,0) {$b=1$};
\node[dot] at (0,-1.05) {};
\node[ann] at (0,-1.491) {$c=5$};
\node[dot] at (-1.05,0) {};
\node[ann] at (-1.491,0) {$a=4$};
\node[ann,anchor=west] at (2.35,0.55) {solid: deleted, $vb$ and $ca$};
\node[ann,anchor=west] at (2.35,0.15) {dashed: inserted, $va$ and $cb$};
\node[ann,anchor=west] at (2.35,-0.35) {every degree is preserved,};
\node[ann,anchor=west] at (2.35,-0.7) {because each vertex meets};
\node[ann,anchor=west] at (2.35,-1.05) {one deleted and one inserted};
\end{tikzpicture}
\caption{A connector, with only the edges it changes drawn. On the running example $d=(3,2,2,2,1)$ at the pivot $v=2$, the neighborhood $S=\{1,3\}$ is carried to $T=S-b+a=\{3,4\}$ by the $2$-switch $\{v1,54\}\mapsto\{v4,51\}$, with witness $c=5$. The four toggled edges are exactly the edges of the four-cycle $v$--$b$--$c$--$a$--$v$, and they alternate around it: the pivot drops $b$ and picks up $a$ while the witness does the reverse. Alternation is why the move is a legal $2$-switch and why every degree is preserved, since each of the four vertices meets one deleted edge and one inserted edge. Figure~\ref{fig:fibers} shows a connector in place among the realizations; this is the same kind of move with everything else removed.}
\label{fig:connector}
\end{figure}

\begin{figure}[!ht]
\centering
\begin{tikzpicture}[x=1cm,y=1cm,
  miniedge/.style={black!55,line width=0.6pt},
  swedge/.style={black!85,line width=1.4pt},
  minidot/.style={circle,fill=black!70,inner sep=0pt,minimum size=1.7mm},
  pivotdot/.style={circle,draw=black!85,fill=white,line width=0.8pt,inner sep=0pt,minimum size=2.6mm},
  fibrebox/.style={rounded corners=2.5pt,draw=black!35,densely dashed},
  qlab/.style={font=\scriptsize,align=center,anchor=south},
  ann/.style={font=\scriptsize,text=black!70},
  cross/.style={black!85,line width=1.3pt},
]
\draw[fibrebox] (-0.82,-2.47) rectangle (0.82,0.82);
\node[qlab] at (0,0.92) {$N(v)=\{1,3\}$\\[-1pt]{\tiny fiber of 2}};
\draw[fibrebox] (1.815,-2.47) rectangle (3.455,0.82);
\node[qlab] at (2.635,0.92) {$N(v)=\{1,4\}$\\[-1pt]{\tiny fiber of 2}};
\draw[fibrebox] (4.45,-0.82) rectangle (6.09,0.82);
\node[qlab] at (5.27,0.92) {$N(v)=\{1,5\}$\\[-1pt]{\tiny fiber of 1}};
\draw[fibrebox] (7.085,-0.82) rectangle (8.725,0.82);
\node[qlab] at (7.905,0.92) {$N(v)=\{3,4\}$\\[-1pt]{\tiny fiber of 1}};
\node[ann,anchor=south] at (3.952,3.3) {the quotient $J(I_v)$, one node per fiber};
\draw[miniedge,line width=0.9pt] (3.202,2.77) -- (4.702,2.77);
\draw[miniedge,line width=0.9pt] (3.202,2.77) -- (2.002,2.05);
\draw[miniedge,line width=0.9pt] (3.202,2.77) -- (5.902,2.05);
\draw[miniedge,line width=0.9pt] (4.702,2.77) -- (2.002,2.05);
\draw[miniedge,line width=0.9pt] (4.702,2.77) -- (5.902,2.05);
\draw[black!45,densely dotted,line width=0.9pt] (2.002,2.05) -- (5.902,2.05);
\node[ann,font=\tiny,anchor=north] at (3.952,1.99) {no edge};
\node[minidot] at (3.202,2.77) {};
\node[ann,font=\tiny,anchor=south] at (3.202,2.87) {$\{1,3\}$};
\node[minidot] at (4.702,2.77) {};
\node[ann,font=\tiny,anchor=south] at (4.702,2.87) {$\{1,4\}$};
\node[minidot] at (2.002,2.05) {};
\node[ann,font=\tiny,anchor=south] at (2.002,2.15) {$\{1,5\}$};
\node[minidot] at (5.902,2.05) {};
\node[ann,font=\tiny,anchor=south] at (5.902,2.15) {$\{3,4\}$};
\draw[miniedge] (3.184e-17,0.52) -- (-0.4945,0.1607);
\draw[miniedge] (3.184e-17,0.52) -- (-0.3056,-0.4207);
\draw[miniedge] (3.184e-17,0.52) -- (0.3056,-0.4207);
\draw[swedge] (-0.4945,0.1607) -- (-0.3056,-0.4207);
\draw[swedge] (0.3056,-0.4207) -- (0.4945,0.1607);
\node[minidot] at (3.184e-17,0.52) {};
\node[pivotdot] at (-0.4945,0.1607) {};
\node[minidot] at (-0.3056,-0.4207) {};
\node[minidot] at (0.3056,-0.4207) {};
\node[minidot] at (0.4945,0.1607) {};
\node[ann,font=\tiny] at (0,-0.66) {$G_{1}$};
\draw[miniedge] (2.635,0.52) -- (2.14,0.1607);
\draw[miniedge] (2.635,0.52) -- (2.329,-0.4207);
\draw[miniedge] (2.635,0.52) -- (2.941,-0.4207);
\draw[swedge] (2.14,0.1607) -- (2.941,-0.4207);
\draw[swedge] (2.329,-0.4207) -- (3.13,0.1607);
\node[minidot] at (2.635,0.52) {};
\node[pivotdot] at (2.14,0.1607) {};
\node[minidot] at (2.329,-0.4207) {};
\node[minidot] at (2.941,-0.4207) {};
\node[minidot] at (3.13,0.1607) {};
\node[ann,font=\tiny] at (2.635,-0.66) {$G_{2}$};
\draw[miniedge] (5.27,0.52) -- (4.775,0.1607);
\draw[miniedge] (5.27,0.52) -- (4.964,-0.4207);
\draw[miniedge] (5.27,0.52) -- (5.576,-0.4207);
\draw[miniedge] (4.775,0.1607) -- (5.765,0.1607);
\draw[miniedge] (4.964,-0.4207) -- (5.576,-0.4207);
\node[minidot] at (5.27,0.52) {};
\node[pivotdot] at (4.775,0.1607) {};
\node[minidot] at (4.964,-0.4207) {};
\node[minidot] at (5.576,-0.4207) {};
\node[minidot] at (5.765,0.1607) {};
\node[ann,font=\tiny] at (5.27,-0.66) {$G_{3}$};
\draw[miniedge] (2.635,-1.13) -- (2.14,-1.489);
\draw[miniedge] (2.635,-1.13) -- (2.329,-2.071);
\draw[miniedge] (2.635,-1.13) -- (3.13,-1.489);
\draw[miniedge] (2.14,-1.489) -- (2.941,-2.071);
\draw[miniedge] (2.329,-2.071) -- (2.941,-2.071);
\node[minidot] at (2.635,-1.13) {};
\node[pivotdot] at (2.14,-1.489) {};
\node[minidot] at (2.329,-2.071) {};
\node[minidot] at (2.941,-2.071) {};
\node[minidot] at (3.13,-1.489) {};
\node[ann,font=\tiny] at (2.635,-2.31) {$G_{4}$};
\draw[miniedge] (3.184e-17,-1.13) -- (-0.4945,-1.489);
\draw[miniedge] (3.184e-17,-1.13) -- (0.3056,-2.071);
\draw[miniedge] (3.184e-17,-1.13) -- (0.4945,-1.489);
\draw[miniedge] (-0.4945,-1.489) -- (-0.3056,-2.071);
\draw[miniedge] (-0.3056,-2.071) -- (0.3056,-2.071);
\node[minidot] at (3.184e-17,-1.13) {};
\node[pivotdot] at (-0.4945,-1.489) {};
\node[minidot] at (-0.3056,-2.071) {};
\node[minidot] at (0.3056,-2.071) {};
\node[minidot] at (0.4945,-1.489) {};
\node[ann,font=\tiny] at (0,-2.31) {$G_{5}$};
\draw[miniedge] (7.905,0.52) -- (7.599,-0.4207);
\draw[miniedge] (7.905,0.52) -- (8.211,-0.4207);
\draw[miniedge] (7.905,0.52) -- (8.4,0.1607);
\draw[miniedge] (7.41,0.1607) -- (7.599,-0.4207);
\draw[miniedge] (7.41,0.1607) -- (8.211,-0.4207);
\node[minidot] at (7.905,0.52) {};
\node[pivotdot] at (7.41,0.1607) {};
\node[minidot] at (7.599,-0.4207) {};
\node[minidot] at (8.211,-0.4207) {};
\node[minidot] at (8.4,0.1607) {};
\node[ann,font=\tiny] at (7.905,-0.66) {$G_{6}$};
\draw[cross] (0.62,0) -- (2.015,0);
\node[ann,anchor=south] at (1.317,0.14) {a connector};
\node[ann,anchor=north west] at (-0.9,-3) {$\circ$ the pivot $v$; the two heavy edges in each of the joined realizations are the 2-switch};
\node[ann,anchor=south] at (-2.62,0.515) {the ground vertices};
\node[minidot] at (-2.62,-0.305) {};
\node[ann] at (-2.62,0.111) {$1$};
\node[pivotdot] at (-3.115,-0.6643) {};
\node[ann] at (-3.51,-0.5358) {$2$};
\node[minidot] at (-2.926,-1.246) {};
\node[ann] at (-3.17,-1.582) {$3$};
\node[minidot] at (-2.314,-1.246) {};
\node[ann] at (-2.07,-1.582) {$4$};
\node[minidot] at (-2.125,-0.6643) {};
\node[ann] at (-1.73,-0.5358) {$5$};
\node[ann,anchor=north] at (-2.62,-2.125) {the pivot is $v=2$};
\end{tikzpicture}
\caption{The quotient at a ground vertex, on $d = (3,2,2,2,1)$ with pivot $v$ the open circle. Its 6 realizations divide into 4 fibers by the neighborhood of $v$, and those neighborhoods are $I_v = \{\{1,3\}, \{1,4\}, \{1,5\}, \{3,4\}\}$, a shifted family; the key at the left numbers the ground vertices, so each set can be read off the drawing, and the six realizations carry the names $G_1$ to $G_6$ they are given in Figure~\ref{fig:realizationgraph}, which draws the same six as one graph. Two fibers are adjacent in the quotient exactly when their neighborhoods differ in one element, which is Corollary 5.5: the quotient is the Johnson graph of $I_v$. The heavy edge joins two realizations in adjacent fibers, and it is a single 2-switch: the four heavy edges inside those two realizations are the ones it toggles. Such an edge is a \emph{connector}, and Section 5.4 counts them: the signed identity of Theorem 5.6 holds at every realization of the source fiber here, with both counts taken on the ground.}
\label{fig:fibers}
\end{figure}
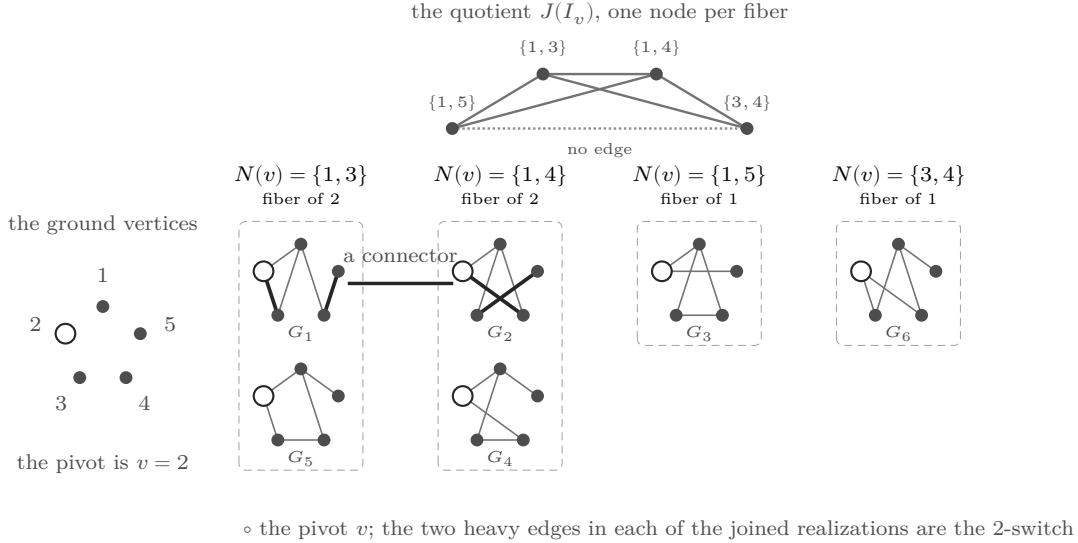

Fix adjacent \(S, T = S − b + a\) in \(I_v\) and write \(f = f_S\) and
\(δ = f(a) − f(b)\). A \textbf{connector} from \(G ∈ Φ_S\) is the
2-switch \(\{vb, ca\} → \{va, cb\}\) for a witness \(c ∉ \{v,a,b\}\)
with \(ca ∈ E(G)\) and \(cb ∉ E(G)\); its four toggled edges form the
alternating four-cycle \(v-b-c-a-v\), drawn on its own in Figure
\ref{fig:connector}. Every connector is an edge of \(G(d)\).

Write \(F = G − v\) for the ground realization of \(f\), and put

\[
w_{ab}(G) = |N_F(a) ∖ (N_F(b) ∪ \{b\})|,
\]

the number of connectors from \(G\). Distinct witnesses toggle distinct
ground edges and each connector determines its witness, so \(G\) has
exactly \(w_{ab}(G)\) distinct connector partners. The \textbf{connector
graph} of the interface is the graph between \(Φ_S\) and \(Φ_T\) in
which \(G ∈ Φ_S\) and \(H ∈ Φ_T\) are joined when some connector carries
one to the other; \(w_{ab}(G)\) is the degree of \(G\) in it. It is
bipartite by construction, the two sides being the two fibers, which are
disjoint because a realization has one neighborhood at \(v\); nothing is
being claimed there. What is not free is the degrees, and Theorem 5.7 is
where they are settled.

\textbf{The orientation is in the pair, not in the edges.} Writing
\(T = S − b + a\) names \(b\) as the neighbor of \(v\) that leaves and
\(a\) as the one that arrives, and that is the only sense in which a
connector goes anywhere: \(w_{ab}(G)\) and \(w_{ba}(H)\) at a partner
\(H ∈ Φ_T\) count the same edges from opposite ends. What is not an edge
count is \(w_{ba}(G)\), the mirror quantity in \(G\)'s own ground graph,
present because Theorem 5.6 states its identity at a single realization.

\textbf{Both counts are taken on the ground, and this matters.} Since
\(b ∈ S\) and \(a ∉ S\), the pivot lies in \(N_G(b)\) and not in
\(N_G(a)\). It is therefore never a connector witness, so reading
\(w_{ab}\) in \(G\) would do no harm; but \(w_{ba}\) read in \(G\) would
count the pivot, which shifts the identity below by one. Every later
use, here and in Section 9, needs the ground form.

\textbf{Theorem 5.6 (signed identity).}
\emph{\(w_{ab}(G) − w_{ba}(G) = f(a) − f(b) = δ\) for every \(G\). In
the target orientation,
\(w_{ba}(H) − w_{ab}(H) = f_T(b) − f_T(a) = 2 − δ\) for every
\(H ∈ Φ_T\).}

\emph{Proof.} Count \(N_F(a) ∖ (N_F(b) ∪ \{b\})\) and
\(N_F(b) ∖ (N_F(a) ∪ \{a\})\) against \(f(a)\) and \(f(b)\). The shared
neighbors and the edge \(ab\) contribute equally to both and cancel. The
second statement is the first applied to \(f_T = f − e_a + e_b\), where
\(e_u\) is the degree function that is \(1\) at \(u\) and \(0\)
elsewhere, so that subtracting it lowers one entry by one.
\(\blacksquare\)

The \(2\) in that second statement is not an artifact of the
bookkeeping. The step moves a single unit of degree from \(a\) to \(b\),
and a difference registers a unit crossing it twice, once on each side.
The end of this subsection exhibits the same \(2\) without degrees at
all, as one witness leaving \(W_{ab}\) and entering \(W_{ba}\).

Figure \ref{fig:fibers} carries this section's apparatus out on one
degree sequence: the fibers of Section 5.1, the shifted family and
Johnson identification of Sections 5.2 and 5.3, a connector, and the
identity above verified at every member of the source fiber.

\textbf{What ``triangle-free'' is a property of.} The next results
hypothesize that a fiber is triangle-free, and the object meant is the
\textbf{fiber as a graph}: \(Φ_S\) carries the subgraph of \(G(d)\)
induced on its realizations, and triangle-free says that subgraph has no
three mutually adjacent realizations, that is no three realizations
pairwise one 2-switch apart. It says nothing about the realizations
themselves, which are graphs of their own and may contain as many
triangles as they like. By Lemma 5.1 the condition can be read
equivalently on \(Φ_S\) or on the smaller realization graph \(G(f_S)\),
and the statements below use whichever is closer to hand.

\textbf{And on a fiber, triangle-free and bipartite are the same
condition.} A fiber is itself a realization graph, \(Φ_S ≅ G(f_S)\) with
\(f_S\) graphical, so Barrus's classification applies to it and Section
3.1's equivalence holds: \(Φ_S\) is triangle-free exactly when it is
bipartite. Theorem 5.7 below is stated with the triangle-free form
because its proof produces triangles, and Theorem 5.8 with the bipartite
form because its proof uses the color classes, but \textbf{the two
theorems assume the same thing} and 5.8 asks for nothing beyond 5.7.

\textbf{What the two of them are for.} Section 10 leaves each fiber by a
connector, and when the fiber is bipartite it can only be traversed by a
Hamilton-laceable walk, whose endpoints are constrained by color. So the
assembly needs \textbf{endpoint control}: it must be able to leave from
a realization of the color the walk requires, not merely from some
realization. Theorem 5.7 supplies that by giving every realization the
same number of connectors, so no color class can be the poor one, and
Theorem 5.8 turns it into the exact count \(|N(Λ)| = k|Λ|\) that Section
9 spends. Neither is needed when the fiber is non-bipartite, because a
Hamilton-connected fiber can be entered and left anywhere.

\textbf{Theorem 5.7 (cross-degree).} \emph{Let \(S\) and \(T\) both be
realizable. If \(Φ_S\) is triangle-free then \(δ ≥ 0\), and every member
of \(Φ_S\) has exactly \(k = \max(1, δ)\) connectors, where
\(δ = f(a) − f(b)\).}

\textbf{The realizability hypothesis}, that \(S\) and \(T\) both lie in
\(I_v\), is the standing assumption of this section: it says each of the
two neighborhoods is actually attained, so \(Φ_S\) and \(Φ_T\) are both
nonempty, which by Lemma 5.2 is the same as \(f_S\) and \(f_T\) both
being graphical. It is not decorative. Without it take
\(f = (3,1,1,1)\), the star: it has a single realization, so \(G(f)\) is
one vertex and vacuously triangle-free, and two of its three degree-one
vertices, taken as \(a\) and \(b\), give \(δ = 0\) and \(w_{ab} = 0\),
against \(\max(1,0) = 1\). What excludes the example is that
\(f − e_a + e_b = (3,0,2,1)\) is not graphical, which is exactly \(T\)
failing to be realizable.

It is convenient to drop the pivot first. Since \(Φ_S ≅ G(f)\) by Lemma
5.1, and \(S − b + a\) is realizable exactly when \(f − e_a + e_b\) is
graphical by Lemma 5.2, the pivot carries no information, and the
statement to prove is: \emph{if \(f\) is graphical with \(G(f)\)
triangle-free and \(f' = f − e_a + e_b\) is graphical, then every
realization \(F\) of \(f\) has \(w_{ab}(F) = \max(1, f(a) − f(b))\).}

Write \(W_{ab} = N_F(a) ∖ (N_F(b) ∪ \{b\})\) and
\(W_{ba} = N_F(b) ∖ (N_F(a) ∪ \{a\})\), so \(w_{ab} = |W_{ab}|\) and
\(w_{ba} = |W_{ba}|\): the count is the size of the set it is named for.

The next three results are internal to the proof of Theorem 5.7 and are
used nowhere else in the paper. They are stated separately because each
is a distinct step and each is machine-checked on its own, not because
anything later cites them.

Two 2-switches out of \(F\) land on an edge of \(G(f)\) exactly when
their symmetric differences share two edges, and a shared edge has the
same status in both, since presence is measured in the same \(F\).

\textbf{Lemma 5.7a.} \emph{If \(G(f)\) is triangle-free then
\(|W_{ba}| ≥ 1\) implies \(|W_{ab}| ≤ 1\), for every realization \(F\)
and every \(a ≠ b\).}

\emph{Proof.} Suppose \(|W_{ba}| ≥ 1\) and \(|W_{ab}| ≥ 2\). Pick
\(x ∈ W_{ba}\) and distinct \(t, t' ∈ W_{ab}\). Then \(F\) contains
\(xb\), \(at\) and \(at'\), and omits \(xa\), \(bt\) and \(bt'\). The
two 2-switches

\[
σ : xb, at \quad ↦ \quad xa, bt \quad τ : xb, at' \quad ↦ \quad xa, bt'
\]

are both legal, and \(F_σ △ F_τ = \{at, bt, at', bt'\}\), which is an
alternating four-cycle on \(\{a, t, b, t'\}\) and so a legal 2-switch.
The five vertices are distinct: \(a ≠ b\); \(x ∉ \{a,b\}\); and
\(t, t' ∈ N_F(a)\) while \(x ∉ N_F(a)\), so \(t, t' ≠ x\), and
\(t, t' ≠ b\) by definition of \(W_{ab}\). So \(F\), \(F_σ\) and \(F_τ\)
are pairwise adjacent --- a triangle in \(G(f)\). \(\blacksquare\)

The three regimes of the identity are used repeatedly below, and Theorem
5.7's proof derives each where it needs it.

Say \(a\) and \(b\) are \textbf{twins} in \(F\) when
\(N_F(a) ∖ \{b\} = N_F(b) ∖ \{a\}\); by the identity this is exactly the
case \(w_{ab}(F) = w_{ba}(F) = 0\).

\textbf{Lemma 5.7b.} \emph{Let \(G(f)\) be triangle-free and
\(f(a) = f(b)\). No 2-switch joins a realization in which \(a, b\) are
twins to one in which they are not.}

\emph{Proof.} Let \(τ\) be the transposition exchanging \(a\) and \(b\).
Since \(f(a) = f(b)\), relabeling by \(τ\) carries realizations of \(f\)
to realizations of \(f\) and 2-switches to 2-switches, so it is an
automorphism of \(G(f)\).

Suppose \(F\) and \(F'\) are adjacent with \(a, b\) twins in \(F\) and
not in \(F'\). Twinhood in \(F\) says exactly that \(τ\) fixes \(F\), so
applying \(τ\) to that adjacency gives \(F\) adjacent to \(F'^τ\).

Since \(f(a) = f(b)\), Theorem 5.6 gives \(w_{ab}(F') = w_{ba}(F')\),
and Lemma 5.7a bounds the common value by one; it is not zero, because
\(a, b\) are not twins in \(F'\). So exactly one vertex \(x\) is
adjacent to \(a\) and not to \(b\), and exactly one \(y\) is adjacent to
\(b\) and not to \(a\). Then \(F'\) and \(F'^τ\) differ in precisely the
four edges \(ax, by, bx, ay\), an alternating four-cycle on \(a-x-b-y\),
so \(F'^τ ≠ F'\) and the two are adjacent.

Now \(F\), \(F'\) and \(F'^τ\) are pairwise adjacent and distinct --- a
triangle in \(G(f)\), contrary to hypothesis. \(\blacksquare\)

\textbf{Corollary 5.7c (twin constancy).} \emph{Let \(G(f)\) be
triangle-free and \(f(a) = f(b)\). Then either every realization of
\(f\) has \(a, b\) twins, or none does.}

\emph{Proof.} By Lemma 5.7b the twin indicator is constant on each
connected component of \(G(f)\), and \(G(f)\) is connected (§2).
\(\blacksquare\)

\emph{Proof of Theorem 5.7.} Write \(x = f(a) − f(b)\), so Theorem 5.6
gives \(w_{ab} − w_{ba} = x\). Lemma 5.7a, applied in both directions,
leaves three regimes.

If \(x ≥ 1\) and \(w_{ba} ≥ 1\), then Lemma 5.7a gives \(w_{ab} ≤ 1\)
while the identity gives \(w_{ab} = w_{ba} + x ≥ 2\). So \(w_{ba} = 0\)
and \(w_{ab} = x = \max(1, x)\).

If \(x ≤ −1\) and \(w_{ab} ≥ 1\), then symmetrically \(w_{ba} ≤ 1\)
while \(w_{ba} = w_{ab} + |x| ≥ 2\). So \(w_{ab} = 0\), and the claimed
formula would be false. \textbf{This regime is empty}, which is the
content of \(δ ≥ 0\), and Theorem 5.4 is what empties it. Its first
clause says that some member of \(Φ_S\) is adjacent to some member of
\(Φ_T\), and under the identification of Lemma 5.1 that is a realization
\(F_0\) of \(f\) with \(w_{ab}(F_0) ≥ 1\). Lemma 5.7a gives
\(w_{ba}(F_0) ≤ 1\), and therefore
\(x = w_{ab}(F_0) − w_{ba}(F_0) ≥ 0\). The appeal is not circular:
Theorem 5.4 is proved before this one, from graphicality and a count,
and does not use it.

If \(x = 0\) then \(w_{ab} = w_{ba}\), and Lemma 5.7a forces the common
value into \(\{0, 1\}\). It is \(0\) exactly when
\(N_F(a) ∖ \{b\} = N_F(b) ∖ \{a\}\), that is when \(a\) and \(b\) are
twins in \(F\). This is where triangle-freeness is used a second time.
By Corollary 5.7c the twin case is all-or-nothing: either every
realization of \(f\) has \(a, b\) twins, or none does. The first
alternative is impossible: the realization \(F_0\) supplied by Theorem
5.4 above has \(w_{ab}(F_0) ≥ 1\), so \(a\) and \(b\) are not twins in
it. Hence no realization is twin, every \(F\) has \(w_{ab}(F) = 1\), and
\(1 = \max(1, 0)\). \(\blacksquare\)

The star that opened this discussion is exactly the first alternative
surviving. In \(f = (3,1,1,1)\) the two leaves are twins, both having
the center as their only neighbor, so \(w_{ab} = 0\); \(G(f)\) is a
single vertex, so twin constancy holds but says nothing; and no non-twin
realization can be built, because \(f' = (3,0,2,1)\) is not graphical.
Realizability of \(T\) enters the proof only as the hypothesis
\(T ∈ I_v\) that Theorem 5.4 needs, and producing \(F_0\) is what it
buys.

The mechanism is the \textbf{vicinal preorder} on the ground --- \(u\)
below \(u'\) when \(N(u) ∖ \{u'\} ⊆ N(u') ∖ \{u\}\), the standard order
on a graph's vertices by neighborhood containment {[}21, Basic
Terminology, p.~3{]}, and \(w_{ua}(F) = 0\) says exactly that --- which
in a triangle-free \(G(f)\) is total except that equal-degree vertices
may sit at distance one. No numerical proximity to a threshold sequence
is involved; the slack ranges of triangle-free and has-triangle fibers
overlap, so Erdős--Gallai slack does not discriminate.

The hypothesis of the next theorem is about the \textbf{fiber}, not the
interface. The connector graph is bipartite for free, as just noted;
what is assumed here is that \(Φ_S\) is bipartite as a graph in its own
right, that is that the 2-switch graph on the realizations of \(f_S\)
has no odd cycle. That can fail, and where it holds it is what forces
laceability rather than Hamilton-connectedness later.

\textbf{Theorem 5.8 (exact color expansion).} \emph{If the fiber \(Φ_S\)
is bipartite as a graph, with its realizations split into color classes
\(C_0\) and \(C_1\), then, with \(k\) as in Theorem 5.7, every
\(Λ ⊆ C_i\) satisfies \(|N(Λ)| = k|Λ|\) in the connector graph.}

\emph{Proof.} Theorem 5.7 gives every member of \(Λ\) exactly \(k\)
connectors, so \(|N(Λ)| = k|Λ|\) requires only that no two members of
\(Λ\) share a target. If two did, they would differ by the 2-switch on
\(\{a,b,c,c'\}\) and so be adjacent, and adjacent realizations lie in
opposite color classes. \(\blacksquare\)

A warning about the sign, because it is easy to over-read Theorem 5.7.
The inequality \(δ ≥ 0\) is asserted only for a triangle-free source. It
is false in general: reversing the orientation of an interface sends
\(δ\) to \(2 − δ\), so a negative value in one direction is a value
above two in the other, and negative values occur at 8,876 of the 88,518
interfaces through ground order seven. No argument in this paper uses
\(δ ≥ 0\) outside the triangle-free case.

\textbf{The difference is constant across the interface, and one vertex
changing sides explains it.} A connector from \(G\) with witness \(c\)
deletes \(ca\) and adds \(cb\), so on the ground \(c\) stops being a
neighbor of \(a\) and becomes one of \(b\). It therefore \textbf{leaves}
\(W_{ab}\) and \textbf{enters} \(W_{ba}\), and nothing else moves:

\[
w_{ab}(H) = w_{ab}(G) − 1, \quad w_{ba}(H) = w_{ba}(G) + 1.
\]

One element crossing between two sets changes their difference by two,
which is the whole of the \(2\) in Theorem 5.6: reading the identity at
\(H\) rather than at \(G\) reverses the roles, worth \(−δ\), and moves
one witness across, worth \(+2\). Taking the two displacements together
gives the same value at every edge of the connector graph,

\[
w_{ab}(G) − w_{ba}(H) = δ − 1,
\]

which Section 9.3 uses at \(δ = 1\). It is machine-checked from the
foundations alone.

\hypertarget{the-quotient-is-not-a-matroid}{%
\subsection{5.5 The quotient is not a
matroid}\label{the-quotient-is-not-a-matroid}}

\textbf{Why that is worth a subsection.} Nothing about a shifted family
suggests it should be the basis family of a matroid, so on its own the
title would be answering a question nobody asked. It is worth saying
because in the companion paper the corresponding quotients \emph{are}
matroids, which is what lets that paper reach for the classical
Hamiltonicity theory of matroid base graphs. This subsection is where
that route stops being available, and Section 7 is what replaces it.

Two things can fail here, and they are worth separating because they are
quantified differently. At a \textbf{fixed} ground vertex the family
\(I_v\) need not have a Gale-greatest member, and five vertices are
enough. At \textbf{every} ground vertex at once the family can fail to
be a matroid basis family, and that needs eleven.

For the first, take \(d = (3,2,2,2,1)\) on \(\{1,…,5\}\) and the pivot
\(v = 2\), so \(k = d(v) = 2\) and the ground of the family is
\(\{1,3,4,5\}\). Of the six two-subsets, \(\{3,5\}\) and \(\{4,5\}\)
have residual \((3,2,1,0)\), which fails the first Erdős--Gallai
inequality since \(3 > \min(2,1)+\min(1,1)+\min(0,1) = 2\), and the
other four are graphical. So by Lemma 5.2

\[
I_v = \{ \{1,3\}, \{1,4\}, \{1,5\}, \{3,4\} \}.
\]

This family is shifted, as Theorem 5.3 requires. But \(\{1,5\}\) and
\(\{3,4\}\) are Gale-incomparable, since \(1 < 3\) while \(5 > 4\), and
each is maximal. \textbf{\(I_v\) has no greatest element.} Every
argument that selects a top neighborhood at a fixed pivot fails here,
and the failure is not a large-order phenomenon.

This family is not merely an example. Relabel its ground \(\{1,3,4,5\}\)
as \(\{1,2,3,4\}\) in order and it becomes \(\{12, 13, 14, 23\}\), which
is \(F(2; 4, 1)\), the smallest Y-family of Section 7.5. So the smallest
obstruction to a greatest element and the smallest obstruction to
Hamilton-connectivity of the quotient are the same family, and the two
obstructions sit at opposite ends of it: the incomparable maximal
members that defeat the greatest-element argument are \(\{1,5\}\) and
\(\{3,4\}\), while the universal pair carrying Theorem 7.7's exception
is the two Gale-\textbf{least} members \(\{1,3\}\) and \(\{1,4\}\). The
maximal pair is traversable --- relabeled, \(\{1,4\}\) and \(\{2,3\}\)
are joined by the Hamilton path \(14, 12, 13, 23\) --- and it is the
least pair that is not.

The second failure is the one that rules out the matroid route entirely
rather than at a chosen pivot.

For a regular degree sequence \(I_v\) is the complete family of
\(k\)-sets and the quotient is the whole Johnson graph; when the degrees
take only two values, one more than the other, it is still a matroid
basis family. Neither regime persists. Eleven ground vertices are enough
for it to fail outright: exactly six degree sequences at that order have
no ground vertex whose realizable family is a matroid, and eleven is the
smallest order at which any does.

This is why Section 7 proves a theorem about shifted families rather
than citing one about matroid basis graphs. It is also why the exception
of Theorem 7.7 exists: by Corollary 7.9 a family with a Gale-greatest
member is never a Y-family, and among shifted families the matroid ones
are exactly those.

\hypertarget{the-separator-buffer-theorem}{%
\section{6. The separator-buffer
theorem}\label{the-separator-buffer-theorem}}

The induction of Section 10 walks a Hamilton path of the quotient at a
single ground vertex, and needs one fiber along that walk to be
Hamilton-connected rather than merely Hamilton-laceable. That fiber is
the \textbf{buffer}, and this section produces a ground vertex carrying
one.

\textbf{Why one fiber must be better than the rest.} The construction
builds two chains inward, one from each of the two prescribed
realizations \(G_0\) and \(G_1\), and they have to meet. At every fiber
except the meeting point the entry is forced but the exit is still free,
so the exit can be chosen in the color class the fiber needs and a
laceable fiber is enough. At the meeting fiber both ends are forced at
once, one by each chain, and nothing controls their colors. A bipartite
fiber is only Hamilton-laceable, so two same-colored ends there have no
spanning walk at all and the construction fails at the one fiber it
cannot reroute. A non-bipartite fiber is Hamilton-connected, which is
indifferent to color. Section 10.2 gives this again in place, with the
two chains in front of the reader.

\textbf{This is the price of prescribing both endpoints.} A Hamilton
path from one prescribed realization to an unspecified second one needs
no buffer: there is a single chain, every fiber keeps a free exit, and
no fiber ever has both ends dictated. The whole of this section is what
the stronger conclusion costs.

\textbf{A reader can stop after Theorem 6.1.} It is the only thing the
construction takes from this section: Sections 8 and 10.1 to 10.3 cite
it and its hypotheses and cite nothing else here, though Section 10.4
discusses this section's own trust dependencies, and the subsections
below are the case analysis that proves it. Nothing later depends on
how.

\textbf{Tyshkevich-indecomposable} is the condition that Section 4.2 has
nothing left to factor, and it is this paper's analogue of the companion
paper's \emph{prime}, since Theorem 4.2 turns a Tyshkevich composition
into a Cartesian product of realization graphs.

\textbf{Activity is not assumed, and does not need to be.} A degree
function is \textbf{active} when every ground vertex has at least two
realizable neighborhoods, as defined in Section 4.1. Earlier drafts
carried it here as a second hypothesis. It is redundant: the only
indecomposable graph with an inactive vertex is the one-vertex graph
{[}32{]}, whose realization graph is a single vertex and so is
bipartite, which this theorem's hypotheses exclude. So indecomposability
and non-bipartiteness already deliver every active line the arguments
below use, and Section 10.1 needs no separate inactive-vertex branch. In
the margin setting the two notions are linked by a theorem of Brualdi
and Manber; here one simply implies the other. Throughout, \(d\) is a
labeled graphical degree function, \(G\) and \(H\) denote realizations
of \(d\), and we write \(N_G(v)\) for the neighborhood of \(v\) in
\(G\).

In the split case this is \textbf{C5}, the companion paper's buffer
lemma, whose statement is the same one line for line with \emph{line} in
place of \emph{ground vertex}. What this section adds is the extension
to every graphical degree sequence; the split branch below invokes
\textbf{C5} rather than reproving it.

\textbf{Theorem 6.1 (separator-buffer).} \emph{Let \(d\) be
Tyshkevich-indecomposable, with \(G(d)\) non-bipartite and not the
\(K_3\) base graph. For any two distinct realizations \(G_0\) and
\(G_1\) there is a ground vertex \(v\) with \(N_{G_0}(v) ≠ N_{G_1}(v)\)
such that some \(v\)-fiber is non-bipartite.}

A triangle of \(G(d)\) is three realizations pairwise joined by a
2-switch, and each of those switches has a \textbf{support}, the four
ground vertices its two exchanged pairs lie on; the \textbf{ground
support} of the triangle is the union of its three edges' supports, that
is the set of every ground vertex any of the three switches touches. It
is a set of ground vertices, not of realizations, and Section 6.2 shows
it always has size four or five.

By the classification of Section 3, a fiber is non-bipartite exactly
when it contains a triangle, so the conclusion says that some triangle
of \(G(d)\) avoids \(v\) in its ground support. We use the two readings
interchangeably.

\hypertarget{reduction-to-a-common-core}{%
\subsection{6.1 Reduction to a common
core}\label{reduction-to-a-common-core}}

Put \(D = D(G_0, G_1) = \{v : N_{G_0}(v) ≠ N_{G_1}(v)\}\), the vertex
support of \(E(G_0) △ E(G_1)\); we write \(D\) bare wherever the pair is
fixed. The two descriptions agree because \(v\) has a different
neighborhood in the two realizations exactly when some edge at \(v\)
differs, so \(D\) is the set of vertices meeting an edge that changes.
Note \(G_0\) and \(G_1\) are arbitrary realizations here and need not be
adjacent, so the changing edges are in general many more than one
2-switch's four. These and only these vertices separate the pair. The
symmetric difference is balanced at every ground vertex, in the sense of
Section 2: \(G_0\) and \(G_1\) have the same degrees, so each vertex
meets as many edges of one as of the other. It therefore decomposes into
alternating closed trails, and the shortest of those has four edges, so
\(|D| ≥ 4\).

Let \(K(d)\) be the intersection of the ground supports of all triangles
of \(G(d)\), the \textbf{common core} this subsection is named for. A
pair fails Theorem 6.1 exactly when \(D(G_0, G_1) ⊆ K(d)\), which is not
a definition but a consequence of the two readings of \(K(d)\) recorded
immediately below. Every triangle support has size four or five, by the
classification below, so \(|D| ≥ 6\) already rules out failure and only
\(|D| ∈ \{4, 5\}\) remains.

Two readings of \(K(d)\) are worth having side by side, because the
second is the one a reader will ask for. A triangle whose ground support
avoids \(v\) lies inside a single \(v\)-fiber, since none of its three
switches touches \(v\); and a fiber is non-bipartite exactly when it
contains a triangle. So \(K(d)\) is equally the set of ground vertices
all of whose fibers are bipartite, and Theorem 6.1 says that no
differing set lies inside it. Since every \(|D| ≥ 4\), it would be
enough to know that \(K(d)\) has at most three members. The rest of the
section is about \(K(d)\) under that name: Section 6.3 rules out failure
when \(G(d)\) has a triangle of the first or third type below, and
Section 6.4 when it has neither.

\hypertarget{the-triangle-classification}{%
\subsection{6.2 The triangle
classification}\label{the-triangle-classification}}

\textbf{What is being classified.} Section 6.1 reduced the problem to
knowing which ground vertices a triangle of \(G(d)\) can touch, so this
subsection determines the possible shapes of a triangle and reads off
its ground support in each. Everything is measured from one chosen
corner \(G\): the two other corners are each one 2-switch away from it,
and \emph{removed} and \emph{added} below always mean relative to that
\(G\), so an edge is removed when it lies in \(G\) and not in the corner
named, and added when the reverse. The outcome is that the support has
size four or five and never more, which is what Section 6.1 consumes.

Let \(G, H_1, H_2\) be a triangle of \(G(d)\) and put
\(Δ_i = E(G) △ E(H_i)\). All three pairs are 2-switch adjacent, so
\(|Δ_1| = |Δ_2| = |Δ_1 △ Δ_2| = 4\), which forces \(|Δ_1 ∩ Δ_2| = 2\).
The two common edges are removed, mixed, or added, and these three cases
are exhaustive.

\textbf{Type I.} Both common edges are removed. The four support
vertices induce \(2K_2\) in \(G\), and the two switches install the
other two perfect matchings. The support has size four.

\textbf{Type II.} One common edge is removed and one added. They share
exactly one endpoint, and after relabeling the switches are
\(\{xy, zt\} → \{xz, yt\}\) and \(\{xy, zt'\} → \{xz, yt'\}\) on five
distinct vertices. The support has size five.

\textbf{Type III.} Both common edges are added. They are disjoint and
form a perfect matching on four vertices; each switch removes one of the
other two matchings, and the switches are distinct, so both occur. Then
\(G\) contains their union, an induced \(C_4\), and omits the common
matching. The support has size four, and this is the exact
complement-dual of Type I.

We call Types I and III \textbf{rich}. Schvöllner and Pastine's degree
formula {[}32{]} \(\deg(G) = 2q_G(2K_2) + 2q_G(C_4) + q_G(P_4)\), in
which \(q_G(H)\) counts the induced copies of \(H\) in the realization
\(G\) and \(\deg(G)\) is its degree in \(G(d)\), gives the same
trichotomy from the switch-degree side: an induced \(2K_2\) or \(C_4\)
supplies two switches and hence a rich triangle, while an induced
\(P_4\) supplies only one.

\hypertarget{the-rich-four-core-lemma}{%
\subsection{6.3 The rich four-core
lemma}\label{the-rich-four-core-lemma}}

\textbf{The idea.} A pair fails only if \emph{every} triangle of
\(G(d)\) uses all of the four or five vertices where the two
realizations differ. That is a strong demand, and this subsection shows
it is too strong to satisfy. Take a vertex \(w\) outside the differing
set \(D\) and ask how it attaches to \(D\): most answers immediately
build a triangle that misses part of \(D\), and the survivors leave
\(D\) so rigidly attached that the realization splits as a Tyshkevich
composition, which the hypothesis forbids. What is left is a degree
sequence with no outside vertex at all, and that is a base case. The
bookkeeping is the table of which attachments survive.

Suppose a failing pair has a rich triangle. Its support has four
vertices and contains \(D\), so \(D\) is exactly that support.
Complementation exchanges Types I and III and changes neither the
realization graph, nor activity, nor Tyshkevich decomposability, nor
triangle supports, so we may analyze Type I.

Write \(D = \{0,1,2,3\}\). Across the three triangle vertices the
induced graphs on \(D\) are the three perfect matchings \(01|23\),
\(02|13\) and \(03|12\), and every adjacency from outside \(D\) is
unchanged.

Let \(w\) lie outside \(D\). If some quartet containing \(w\) and three
vertices of \(D\) induced \(2K_2\) or \(C_4\) in any of the three
realizations, it would supply a rich triangle whose support does not
contain \(D\), contradicting the common-core condition. That leaves
exactly three possibilities: \(N_D(w) = ∅\), \(N_D(w) = D\), or
\(N_D(w) = D − \{j\}\) for a single \(j\).

Sizes one and two are excluded by one witness. Suppose
\(1 ≤ |N_D(w)| ≤ 2\) and take the matching that pairs the members of
\(N_D(w)\) with each other --- any of the three, if there is only one
member. Its other pair \(\{k, l\}\) is then disjoint from \(N_D(w)\), so
for \(i ∈ N_D(w)\) the quartet \(\{w, i, k, l\}\) has the edges \(wi\)
and \(kl\) and no others: \(wk\) and \(wl\) are absent because
\(k, l ∉ N_D(w)\), and \(ik\) and \(il\) are absent because the matching
pairs \(i\) elsewhere. That is an induced \(2K_2\), and its support
misses the fourth vertex of \(D\).

Name these outside types for the attachment that defines them: \(Z_∅\)
when \(N_D(w) = ∅\), \(Z_D\) when \(N_D(w) = D\), and \(Z_j\) when
\(N_D(w) = D − \{j\}\). Substituting two outside types, the status of
the edge between them, and the three matchings into Types I, II and III
gives a complete table of permitted statuses.

\begin{longtable}[]{@{}llll@{}}
\toprule\noalign{}
& \(Z_∅\) & \(Z_D\) & \(Z_j\) \\
\midrule\noalign{}
\endhead
\bottomrule\noalign{}
\endlastfoot
\(Z_∅\) & \(0\) & \(0\) or \(1\) & \(0\) \\
\(Z_D\) & & \(1\) & \(1\) \\
\(Z_i\) & & & none permitted \\
\end{longtable}

Each exclusion is checkable on the single realization \(G[D] = 01\|23\),
and we give the witness for every one of them rather than assert them,
writing \(u\) and \(w\) for the two outside vertices. The
\(Z_∅\)--\(Z_∅\) case is a Type I support \(\{0, 1, u, w\}\): \(01\) and
\(uw\) are edges and neither \(u\) nor \(w\) has any neighbor in \(D\),
so the quartet induces \(2K_2\). The rest are Type II, in the labeling
\(\{xy, zt\} → \{xz, yt\}\) and \(\{xy, zt'\} → \{xz, yt'\}\) of Section
6.2, which needs \(xy, zt, zt' ∈ E\) and \(xz, yt, yt' ∉ E\) on five
distinct vertices:

\begin{longtable}[]{@{}lllllll@{}}
\toprule\noalign{}
excluded cell & \(x\) & \(y\) & \(z\) & \(t\) & \(t'\) & the support
misses \\
\midrule\noalign{}
\endhead
\bottomrule\noalign{}
\endlastfoot
\(Z_D\)--\(Z_D\), non-edge & \(u\) & \(0\) & \(w\) & \(2\) & \(3\) &
\(1\) \\
\(Z_∅\)--\(Z_0\), edge & \(0\) & \(1\) & \(w\) & \(2\) & \(u\) &
\(3\) \\
\(Z_D\)--\(Z_0\), non-edge & \(2\) & \(w\) & \(1\) & \(0\) & \(u\) &
\(3\) \\
\(Z_0\)--\(Z_0\), either status & \(1\) & \(0\) & \(2\) & \(u\) & \(w\)
& \(3\) \\
\(Z_0\)--\(Z_1\), \(uw\) a non-edge & \(u\) & \(1\) & \(w\) & \(2\) &
\(3\) & \(0\) \\
\(Z_0\)--\(Z_1\), \(uw\) an edge & \(1\) & \(0\) & \(w\) & \(u\) & \(3\)
& \(2\) \\
\end{longtable}

\textbf{The table is a reference, not a reading assignment.} Every entry
is read off the type definitions and the matching \(01\|23\) by the same
routine, so we work one row here and leave the rest to the audit cited
below. Take the \(Z_∅\)--\(Z_0\) row: \(01\) is a matching edge, \(w2\)
and \(wu\) are edges because \(w\) is \(Z_0\) and the cell assumes the
\(Z_∅\)--\(Z_0\) edge, while \(0w\) is absent because \(w\) is \(Z_0\),
\(12\) is absent in the matching, and \(1u\) is absent because \(u\) is
\(Z_∅\). The \(Z_0\)--\(Z_0\) row never consults the status of \(uw\)
--- its six edge conditions are \(10 ∈ E\), \(2u, 2w ∈ E\) and
\(12, 0u, 0w ∉ E\) --- so one witness excludes both statuses there.
\textbf{The \(Z_0\)--\(Z_1\) cell needs two}, because a Type II
labelling requires \(xz ∉ E\), and taking \(x = u\), \(z = w\) makes
that condition \(uw ∉ E\). The second row supplies the edge case, and it
is the mirror image of the difficulty: its removed pair is
\(\{10, uw\}\), so it consumes the very edge the first row cannot have.
Every displayed support misses a vertex of \(D\), so each excluded cell
contradicts failure.

These witnesses were also checked by enumeration: for each cell, the
configuration was built explicitly, the whole realization graph of its
degree sequence enumerated, and a triangle whose ground support misses a
vertex of \(D\) found in every case
(\texttt{compute/sec63\_witnesses.py}).

The table says the \(Z_∅\)-vertices form an independent set, the
\(Z_D\)-vertices form a clique, every \(Z_∅\) is anticomplete and every
\(Z_D\) complete to \(D\) and to the at most one \(Z_j\), and the
\(Z_∅\)--\(Z_D\) edges are arbitrary. Suppose \(Z_∅ ∪ Z_D\) is nonempty,
and write \(S\) for the split graph induced on it, with \(Z_D\) its
clique side and \(Z_∅\) its independent side --- \textbf{clique side
first}, as Section 4.2 fixes for a splitted factor. The cross pattern
the table describes is exactly the composition's: the clique side
complete to \(D ∪ \{Z_j\}\), the independent side anticomplete to it. So
the realization is \(S ∘ G[D ∪ \{Z_j\}]\). Composing two graphs this way
composes their degree sequences, so \(d\) itself is a nontrivial
Tyshkevich composition, contrary to indecomposability. (Section 4.2
defines \(∘\) on sequences; the promotion to graphs is legitimate in
this direction, which is the only one used here, and Theorem 4.2's proof
is where the two readings are reconciled.)

If a lone \(Z_0\) remains, take \(G[D] = 01|23\) and
\(N_D(Z_0) = \{1,2,3\}\). The switch \(\{01, 3Z_0\} → \{0Z_0, 13\}\)
produces a realization inducing the \(C_4\) with edges
\(13, 1Z_0, 23, 2Z_0\) on \(\{1,2,3,Z_0\}\), whose Type III triangle
misses \(0\); the common core fails again.

So an active indecomposable failure with a rich triangle has no outside
vertex at all. It is \(d = (1^4)\) or, after complementation,
\(d = (2^4)\), and \(G(d) = K_3\). Both are base cases.

\hypertarget{the-remaining-branch}{%
\subsection{6.4 The remaining branch}\label{the-remaining-branch}}

\textbf{C5 (buffer existence).} \emph{Let \(G(R,S)\) be the interchange
graph of a nonempty class of zero-one matrices with fixed margins, with
at least three rows and at least three columns, every line active, no
invariant position, and \(G(R,S)\) non-bipartite. Then any two distinct
matrices of the class differ on some line one of whose fibers is
non-bipartite.}

The companion states \textbf{C5} inside a section carrying two further
standing hypotheses: that the class is invariant-free, and that
\(|V(G(R,S))| > 6\). We have replaced the first by the absence of an
invariant position, which is the same condition. We have dropped the
second on the companion's own authority: the paragraph opening its
Section 5 says that the order bound \emph{``selects the §5 branch rather
than the §6 branch; it is not used in the local proofs of Lemmas
5.8--5.10.''} The bound is a routing condition in their induction, not a
hypothesis of the lemma, and Section 6.4 is not running their induction
--- it wants the conclusion only. The accompanying formalization agrees:
the buffer lemma is formalized with the hypotheses stated above and no
order bound. We record the discrepancy rather than let it pass, because
a referee checking \textbf{C5} against the printed statement would find
the bound in the ambient regime.

Suppose no rich triangle occurs. Then no realization contains an induced
\(2K_2\) or \(C_4\), since either quartet would generate a Type I or
Type III triangle, so every realization is \textbf{pseudo-split}: free
of induced \(2K_2\) and of induced \(C_4\).

Let \(G\) be a realization that is not split. By the forbidden-subgraph
characterization of split graphs {[}15{]}, it contains an induced
\(2K_2\), \(C_4\) or \(C_5\); the first two are excluded, so it contains
an induced \(C_5\). \textbf{That characterization is not reproved here
and it is not assumed either}: it is machine-checked in the accompanying
formal development, which is why Section 10.4 lists it among the two
facts proved there rather than among the seven taken on trust. A reader
who wants it on the page has Földes and Hammer. In a
\(\{2K_2, C_4\}\)-free graph every vertex outside that cycle is complete
or anticomplete to it, since a mixed cyclic neighborhood gives a
forbidden quartet on that vertex and three cycle vertices. The complete
vertices form a clique, because two nonadjacent ones together with two
nonadjacent cycle vertices induce \(C_4\); the anticomplete vertices
form an independent set, because an edge there together with a cycle
edge induces \(2K_2\). Hence \(G = S ∘ C_5\) unless \(G = C_5\), where
\(S\) is the split graph on the vertices outside the cycle with the
complete ones as its clique side; and indecomposability forces
\(G = C_5\). But \(d(C_5) = (2^5)\) has realization graph
\(K_{6,6} − 6K_2\), which is bipartite, contrary to hypothesis. So \(G\)
is split.

For a split degree sequence the whole of \(G(d)\) can be read on a
matrix, and the rest of this branch works there. Two facts license the
translation, and both are \textbf{proved in the accompanying formal
development rather than on this page}, which is why Section 10.4 lists
them apart from the seven taken on trust: the characterization used just
above, and \textbf{every realization of \(d\) is split at one and the
same partition of the ground}. A reader checking the development will
find them as \texttt{split\_or\_induced\_fiveCycle} and
\texttt{splitAt\_all\_realizations}. The second is what makes the
translation single-valued: the partition is fixed once for the whole
class, not chosen afresh for each realization.

Fix that partition, let \(R\) and \(S\) be the degree vectors it induces
on the clique side and the independent side, and let \(A(R,S)\) be the
class of zero-one matrices with those margins. In every realization the
clique side is complete and the independent side empty, so all the
variation sits in the cross entries, and reading each realization by its
cross entries is a graph isomorphism from \(G(d)\) onto the interchange
graph \(G(R,S)\). Under it:

\begin{longtable}[]{@{}
  >{\raggedright\arraybackslash}p{(\columnwidth - 2\tabcolsep) * \real{0.5000}}
  >{\raggedright\arraybackslash}p{(\columnwidth - 2\tabcolsep) * \real{0.5000}}@{}}
\toprule\noalign{}
\begin{minipage}[b]{\linewidth}\raggedright
in \(G(d)\)
\end{minipage} & \begin{minipage}[b]{\linewidth}\raggedright
in \(A(R,S)\)
\end{minipage} \\
\midrule\noalign{}
\endhead
\bottomrule\noalign{}
\endlastfoot
a realization & a matrix \\
a ground vertex, used as a pivot & a \emph{line}, that is, a row or a
column \\
the fiber at that pivot & the matrices agreeing on that line \\
a 2-switch & an interchange \\
\end{longtable}

Two more terms are needed, and they are the ones \textbf{C5} asks for. A
line is \emph{active} when it is not the same in every matrix of the
class, and a position is \emph{invariant} when its entry is the same in
every matrix. Both are properties of the class, not of any one matrix.

Two facts transfer to the matrix. \textbf{No line of the class is
constant}: a constant row would make its clique vertex's whole
neighborhood constant, since the rest of that neighborhood is the
clique, and that vertex would then be inactive --- which
indecomposability forbids here, by the argument at the head of this
section. \textbf{Indecomposability also rules out an invariant
position}, and this is the step that needs a theorem rather than
bookkeeping. If the class had one, then by a result of Ryser, in the
form given by Brualdi {[}9, Theorem 3.4.1{]}, every matrix of the class
would decompose with an all-ones block in one corner and an all-zeros
block in the opposite one. Reading that block back on the graph gives a
nontrivial Tyshkevich composition, so an indecomposable sequence leaves
no invariant position.

The point is worth stating separately because activity alone does not
give it. There are classes in which every row and every column varies
while some single entry does not, and the smallest is a split sequence
on eight vertices.

When the class has at least three active rows and three active columns,
every hypothesis of \textbf{C5} is now in hand: the class is nonempty
because \(d\) is graphical, it has three rows and three columns by the
case assumption, every line is active and it has no invariant position,
both by the paragraph beginning ``Two facts transfer to the matrix'',
and its interchange graph is \(G(d)\), which Theorem 6.1 assumes
non-bipartite. So \textbf{C5} produces, for the two prescribed matrices,
a separating line with a non-bipartite fiber.

The remaining case is that some side has fewer than three active lines,
and there it always has \emph{exactly} two. Two distinct matrices of the
class have the same margins, so their difference is a nonzero integer
matrix whose every row sum and every column sum is zero. A row of such a
matrix carrying one nonzero entry must carry a second of the opposite
sign, and likewise a column; so from any nonzero entry, its row forces a
second varying column and its column a second varying row. Theorem 6.1
is given two distinct realizations, so the class does have two distinct
matrices, and neither side can have fewer than two active lines.

Say the two are rows, and let \(f\) be the number of active opposite
lines. This \(f\) is a count, unrelated to the residual degree function
\(f_S\) of Sections 5 and 9; the collision is local to this branch,
which never mentions a residual.

\emph{The class is the Johnson graph \(J(f,k)\).} Each active column has
exactly one of its two entries equal to one: its sum is fixed by the
margins, and a column of a two-row matrix whose sum is \(0\) or \(2\) is
constant, hence inactive. So a matrix of the class is named by the set
of columns placing their one in the first row, a set whose size is that
row's sum, say \(k\); and every \(k\)-subset occurs, since any
assignment of \(k\) of the \(f\) columns to the first row meets both
margins. An interchange on a two-row matrix carries a \(2 × 2\)
submatrix \([[1,0],[0,1]]\) to \([[0,1],[1,0]]\), moving one column out
of the named set and one in, so adjacency is symmetric difference two
and the class is \(J(f,k)\). Both rows being active forces
\(1 ≤ k ≤ f−1\), and in particular \(f ≥ 2\).

\emph{Triangles in a Johnson graph.} For \(1 ≤ j ≤ m−1\), the graph
\(J(m,j)\) contains a triangle if and only if \(m ≥ 3\). If \(j = 1\) it
is \(K_m\). If \(j ≥ 2\), take a \((j+1)\)-subset \(Y\), which exists
because \(j+1 ≤ m\), and three distinct \(x, y, z ∈ Y\), which exist
because \(j+1 ≥ 3\); then \(Y ∖ \{x\}\), \(Y ∖ \{y\}\) and \(Y ∖ \{z\}\)
are \(j\)-subsets pairwise at symmetric difference two. Conversely
\(m ≤ 2\) leaves only \(J(2,1) = K_2\).

\emph{The pivot is an opposite line.} Fixing one of the two active lines
fixes the whole named set and so fixes the matrix, leaving a singleton
fiber, which carries no triangle. A differing opposite line does exist,
since the prescribed realizations name distinct \(k\)-sets and any
element of their symmetric difference is such a line. At an opposite
line the two fibers are the matrices whose named set contains it and
those whose does not, namely \(J(f−1,k−1)\) and \(J(f−1,k)\).

\emph{If \(f ≥ 4\),} at least one of those two is nondegenerate:
\(J(f−1,k)\) is unless \(k = f−1\), and \(J(f−1,k−1)\) is unless
\(k = 1\), and those two exceptions coincide only at \(f = 2\), so at
most one of them holds. The nondegenerate fiber is a Johnson graph on
\(f − 1 ≥ 3\) ground elements, so by the criterion above it carries a
triangle, and the line is the required pivot.

\emph{If \(f = 3\),} then \(k ∈ \{1,2\}\) and the class is \(J(3,1)\) or
\(J(3,2)\), each of which is \(K_3\), so \(G(d) = K_3\) and the
hypothesis excludes it. Every line is active, so the matrix is exactly
\(2 × 3\) with all column sums one, and the split graph it codes has two
clique vertices, of degrees \(k+1\) and \(4−k\), and three independent
vertices of degree one: the chair sequence \((3,2,1,1,1)\), or its
complement according to which side is taken as the clique.

\emph{If \(f = 2\),} then \(k = 1\) and the class is \(J(2,1) = K_2\),
which is bipartite, contrary to hypothesis. Here the matrix is
\(2 × 2\), so up to labels there is only the one instance,
\(d = (2,2,1,1)\), split at the middle pair.

Since \(f ≥ 2\), the three cases exhaust the possibilities. Every claim
in this branch was recomputed independently over the whole
two-active-row family for \(f ≤ 12\), and the two-active-lines count
over 215,291 interchange classes, with no exception
(\texttt{compute/sec64\_two\_active\_lines\_seal.py}).

This exhausts both branches and proves Theorem 6.1. \(\blacksquare\)

\hypertarget{what-the-theorem-does-not-give}{%
\subsection{6.5 What the theorem does not
give}\label{what-the-theorem-does-not-give}}

Theorem 6.1 selects a pivot by two conditions only: it separates the
prescribed pair, and it carries a non-bipartite fiber. Section 8 shows
that the assembly needs a third condition, and that the third one is not
implied by the first two.

\hypertarget{the-quotient-theorem}{%
\section{7. The quotient theorem}\label{the-quotient-theorem}}

We prove that the Johnson graph of a shifted family is
Hamilton-connected apart from one family of exceptions, which we
classify. Nothing in this section refers to degree sequences.

\textbf{This is much the longest section of the paper, and most of it is
not used anywhere else in it.} That is deliberate. The characterization
is a theorem in its own right, and we prove it in full rather than
proving only the fragment the realization-graph argument consumes. What
that argument uses is two statements: Theorem 7.6, the classification,
including the explicit arm structure; and Theorem 7.7, the main theorem.
\textbf{A reader who wants only the realization graphs can take those
two and go on to Section 8}; everything between here and there is their
proof. Corollary 7.9, the principal case, is why the exception stays
invisible in the setting the literature covers, and the literature
subsection at the end of this section reads it that way; Sections 8 to
10 never invoke it. Sections 11 and 12 also report Theorem 7.11 and
Lemma 7.4, but nothing rests on them.

\textbf{The result, stated here so the reader knows the destination.}
Let \(F\) be a shifted family of \(k\)-sets and let \(A ≠ B\) be two of
its members. Then the Johnson graph \(J(F)\) has a Hamilton \(A\)--\(B\)
path --- \textbf{unless} \(F\) is a \emph{Y-family} and \(\{A, B\}\) is
its pair of Gale-least members, the two smallest in the order that
Section 7.1 fixes, in which case it has none. That is Theorem 7.7.

The exception is small and explicit. The Y-families are a
three-parameter list \(F(k; r, s)\), written out in Section 7.5, so
whether a given family is one is decided by inspection rather than by
search; Corollary 7.8 draws the consequence, that a shifted family which
is not a Y-family has \(J(F)\) Hamilton-connected with no exceptional
pair at all. Corollary 7.9 adds that a family with a Gale-greatest
member is never a Y-family, which is why the exception cannot arise in
the case the existing literature covers.

\textbf{The route, so the reader can see where each piece goes.} Section
7.1 fixes the Gale order and the two smallest sets \(m_1\) and \(m_2\).
Section 7.2 counts the decrements of a set, and that count drives
Sections 7.3 to 7.5: a set far from the bottom has many decrements, so
anything that can go wrong is forced to sit low. Section 7.3 locates the
two smallest members and shows that a family with three or more is not
bipartite, which is why the theorem can ask for Hamilton-connectivity at
all. Section 7.4 exhibits one shape that certainly fails --- the
\textbf{clique-sum}, two members adjacent to everything with two cliques
hanging off them --- and Section 7.5 identifies the clique-sums as
exactly one explicit family, the \textbf{Y-families}; that is Theorem
7.6.

Note the direction, because it is easy to read past. Up to that point we
have shown only that Y-families fail. \textbf{The converse, that nothing
else does, is Theorem 7.7, and it is proved in Section 7.8}, by an
induction whose one hard input is the crossing lemma, Theorem 7.11;
Section 7.7 supplies a duality that halves its case analysis. Section
7.9 then places the result against the literature --- the matroid class
where the exception cannot occur, and the Gray-code results that come
nearest to it.

\hypertarget{shifted-families}{%
\subsection{7.1 Shifted families}\label{shifted-families}}

Fix the ground \([n] = \{1, …, n\}\) with its natural order and an
integer \(k\) with \(1 ≤ k ≤ n\).

The \textbf{Gale order} compares two \(k\)-sets by sorted entries:
writing \(X = \{x_1 < ⋯ < x_k\}\) and \(Y = \{y_1 < ⋯ < y_k\}\), we set
\(X ≼ Y\) when \(x_i ≤ y_i\) for every \(i\). A family \(F\) of
\(k\)-subsets is \textbf{shifted of rank \(k\)} when it is nonempty and
closed under a single decrement: \(X ∈ F\), \(j ∈ X\), \(i < j\),
\(i ∉ X\) imply \(X − j + i ∈ F\). We call \(X − j + i\) a
\textbf{decrement} of \(X\), and write \((i, j]\) and \((i, j)\) for the
integer intervals \(\{i+1, …, j\}\) and \(\{i+1, …, j−1\}\). Throughout,
\(12\) abbreviates the 2-subset \(\{1,2\}\) and \(123\) the 3-subset
\(\{1,2,3\}\); no digit repeats.

The order has a second description and we use the two interchangeably.
Writing \(c_Z(t) = |Z ∩ [t]|\), an element \(x_i\) of \(X\) satisfies
\(x_i ≤ t\) exactly when \(c_X(t) ≥ i\), so for two \(k\)-sets

\[
X ≼ Y \quad \text{ if } \text{ and } \text{ only } \text{ if } \quad c_X(t) ≥ c_Y(t) \quad \text{ for } \text{ every } t ∈ [n],
\]

that is, \(X\) has at least as many elements below every threshold. The
equivalence requires the two sets to have equal size, which holds
throughout. The formal development is written in the threshold form and
proves the equivalence, so a result stated below in either description
is machine-checked in the description it is stated in. Gale domination
is generated by single decrements, so the shifted families are exactly
the nonempty down-sets of \(≼\), and we use both descriptions.

The \textbf{Johnson graph of \(F\)}, written \(J(F)\), has vertex set
\(F\) with \(X\) adjacent to \(Y\) when \(|X △ Y| = 2\) --- that is,
when one is obtained from the other by exchanging a single element for a
single element outside it. It is the subgraph of \(J(n,k)\) induced on
\(F\), and taking \(F\) to be all \(k\)-subsets recovers \(J(n,k)\).

A family is \textbf{principal} when it has a Gale-greatest member \(M\),
so that \(F = \{Z : Z ≼ M\}\) is everything at or below one set. The
principal families are exactly the basis families of shifted matroids,
that is, of matroids whose bases form a shifted family. Corollary 7.9
and the discussion after it identify exactly which families these are.

Everything below is order-theoretic, so all statements transfer along an
order isomorphism of the ground. We use this without comment when
passing to a punctured ground \([n] ∖ \{e\}\).

Write \(m_1 = [k]\) and \(m_2 = [k−1] ∪ \{k+1\}\). The second is a
\(k\)-subset of \([n]\) only when \(k < n\); if \(k = n\) the only
\(k\)-subset is \([n]\) itself, so \(|F| = 1\). Every statement below
that names \(m_2\) therefore either assumes \(|F| ≥ 2\) or, like the
census of Lemma 7.1, presupposes a \(k\)-set other than \(m_1\) and so
forces \(k < n\) on its own.

\hypertarget{the-decrement-census}{%
\subsection{7.2 The decrement census}\label{the-decrement-census}}

Everything in Sections 7.3 to 7.5 turns on \textbf{how many decrements a
\(k\)-set has}, and this subsection counts them. A set with few
decrements must be close to \(m_1\), and the classification is
eventually read off the sets with exactly one or two.

The count depends only on how far \(W\) is from \(m_1 = [k]\), so the
statement is organized by that. When \(W\) misses exactly one element of
\([k]\), the two numbers \(i\) and \(a\) below are simply \emph{which}
element of \([k]\) is missing and \emph{which} element outside it is
present --- they are determined by \(W\), not chosen.

\textbf{Lemma 7.1 (census).} \emph{Let \(W\) be a \(k\)-set with
\(W ≠ m_1\).}

\emph{(a) If \(|W ∖ [k]| = 1\), write \(W = [k] − i + a\) with
\(i ≤ k < a\). The distinct decrements of \(W\) are \(W − a + a'\) for
\(a' ∈ \{i\} ∪ (k, a)\), and \(W − j + i\) for \(j ∈ (i, k]\). So \(W\)
has exactly \((a − k) + (k − i) = a − i\) of them.}

\emph{(b) \(a − i = 1\) occurs only for \(W = m_2\), whose sole
decrement is \(m_1\).}

\emph{(c) For \(k ≥ 2\), \(a − i = 2\) occurs only for
\(W = [k] − (k−1) + (k+1)\) and \(W = [k] − k + (k+2)\); for \(k = 1\)
only for \(W = \{3\}\). In each of these cases the decrement set is
exactly \(\{m_1, m_2\}\).}

\emph{(d) If \(|W ∖ [k]| ≥ 2\) then \(W\) has at least four distinct
decrements.}

\emph{Proof.} Parts (a), (b) and (c) are immediate from the definition
of a decrement --- no family enters them, which is what makes this a
census --- and we prove (d). Let \(a > b > k\) be the two largest
elements of \(W\) outside \([k]\). The decrements at the \(a\)-slot
number \(|[a−1] ∖ W| ≥ (a−1) − (k−1) = a − k ≥ 2\), since \(a ≥ k+2\).
Since \(a\) and \(b\) lie in \(W\) and are both at least \(b\), at most
\(k−2\) elements of \(W\) lie in \([b−1]\), so the decrements at the
\(b\)-slot number \(|[b−1] ∖ W| ≥ (b−1) − (k−2) = b − k + 1 ≥ 2\), since
\(b ≥ k+1\). The two groups are distinct: one deletes \(a\) and retains
\(b\), the other deletes \(b\) and retains \(a\). \(\blacksquare\)

Every part of the census is consumed below, by Lemmas 7.2 and 7.5, by
Corollary 7.3 and by Theorem 7.6, so we state them rather than leaving
the census as a remark.

\hypertarget{the-bottom-of-a-shifted-family}{%
\subsection{7.3 The bottom of a shifted
family}\label{the-bottom-of-a-shifted-family}}

\textbf{Why the bottom is where to look.} A shifted family is nearly
unconstrained at the top: which large sets it contains is exactly what
distinguishes one family from another. At the bottom it has no freedom
at all. Closure downward forces \(m_1 = [k]\) into every family, then
forces \(m_2 = [k−1] ∪ \{k+1\}\) into every family with a second member,
and in that order. So every shifted family in the paper has the same two
smallest elements, and that fixed foot is what the rest of the section
stands on. It is also, as Section 7.5 will show, exactly where the one
exception lives.

\textbf{Lemma 7.2.} \emph{Let \(F\) be shifted of rank \(k\) with
\(|F| ≥ 2\). Then:} \emph{(i) \(m_1 ∈ F\) and \(m_2 ∈ F\);} \emph{(ii)
\(m_2 ≼ X\) for every \(X ∈ F ∖ \{m_1\}\), so \(m_1\) and \(m_2\) are
the two Gale-least members and the second-least is unique;} \emph{(iii)
if \(|F| = 2\) then \(F = \{m_1, m_2\}\) and \(J(F)\) is a single edge;}
\emph{(iv) if \(|F| = 3\) then \(J(F)\) is a triangle.}

\emph{Proof.} (ii) A member \(X ≠ m_1\) has sorted values \(x_i ≥ i\)
for all \(i\), and \(x_k ≥ k+1\), since \(X ⊆ [k]\) would force
\(X = m_1\). As \(m_2\) has sorted values \((1, …, k−1, k+1)\), we get
\(m_2 ≼ X\). (i) Pick any \(X ≠ m_1\), which exists; down-closure with
(ii) gives \(m_2 ∈ F\), and then \(m_1 ∈ F\). (iii) Immediate from (i)
and (ii), with \(|m_1 △ m_2| = 2\). (iv) Let \(F = \{m_1, m_2, W\}\).
Every decrement of \(W\) lies in \(F ∖ \{W\}\), so \(W\) has at most
two; Lemma 7.1(d) forces \(|W ∖ [k]| = 1\), and 7.1(b) excludes
\(a − i = 1\), so \(a − i = 2\). For \(k = 1\), \(W = \{3\}\) and all
three singletons are pairwise adjacent. For \(k ≥ 2\), \(W\) is one of
the two sets of Lemma 7.1(c). Each has \(|W △ m_1| = 2\), and each is
one swap from \(m_2\), since \(([k−2] ∪ \{k, k+1\}) △ m_2 = \{k−1, k\}\)
and \(([k−1] ∪ \{k+2\}) △ m_2 = \{k+1, k+2\}\). \(\blacksquare\)

\textbf{Corollary 7.3.} \emph{Let \(F\) be shifted of rank \(k\). If
\(|F| ≥ 3\) then \(J(F)\) contains a triangle, so \(J(F)\) is not
bipartite.}

Shiftedness is not decoration here. Without it the statement is false:
\(F = \{12, 13, 34\}\) has three members and \(J(F)\) is the path
\(12 – 13 – 34\).

\emph{Proof.} Lemma 7.2(i) puts \(m_1\) and \(m_2\) in \(F\), and they
are distinct, so \(|F| ≥ 3\) leaves \(F ∖ \{m_1, m_2\}\) nonempty;
choose \(W\) Gale-minimal in it. Every decrement of \(W\) lies in \(F\)
and is Gale-below \(W\), so minimality places it in \(\{m_1, m_2\}\),
and \(W\) has at most two decrements. As in Lemma 7.2(iv) this gives
\(a − i = 2\), so \(W\) is adjacent to both \(m_1\) and \(m_2\), which
are adjacent to each other. \(\blacksquare\)

We take a Gale-\emph{minimal} member rather than ``the three least
members'': the Gale order is a partial order, so two \(k\)-sets can be
incomparable and there need be no ``three least'' to speak of; it is not
total above \(m_2\), and \(F = \{123, 124, 125, 134\}\) is shifted with
\(125\) and \(134\) incomparable.

Corollary 7.3 lets the main theorem be stated as Hamilton-connectivity,
since a bipartite graph on three or more vertices is never
Hamilton-connected.

\hypertarget{clique-sums}{%
\subsection{7.4 Clique-sums}\label{clique-sums}}

\textbf{Why a condition on the graph before a condition on the family.}
Theorem 7.7 will say that exactly one thing obstructs a Hamilton path
between two prescribed members, and the way to find that one thing is to
ask first what shape in \(J(F)\) could block such a path at all,
ignoring for the moment that \(J(F)\) came from a family. This
subsection isolates the shape and Lemma 7.4 shows it is fatal. Section
7.5 then translates it back into a condition on \(F\) that can be
checked without drawing the graph.

Let \(G\) be a graph and \(u ≠ v\) two of its vertices. The pair
\((G, \{u,v\})\) is a \textbf{clique-sum} when \(u\) and \(v\) are
adjacent to every other vertex, and the remaining vertices split into
two nonempty parts that are each complete, with no edge between them.
Both halves are needed: without completeness the parts would only be the
components of \(G − \{u,v\}\), which is a weaker condition and not what
the name records. Equivalently \(G ≅ K_2 ∨ (K_p ⊔ K_q)\) for some
\(p, q ≥ 1\), where \(⊔\) is disjoint union and \(∨\) the join, which
adds every edge between its two sides: adding \(u\) and \(v\) to each
part turns it into a clique, so \(G\) is two cliques \(K_{p+2}\) and
\(K_{q+2}\) glued along the single shared edge \(uv\), which is what the
name records. The pair \(\{u,v\}\) is the \textbf{universal pair}, and
it is determined by \(G\) alone, since a member of one clique is
non-adjacent to every member of the other and no other pair is adjacent
to everything.

\textbf{And a family inherits it.} The vertices of \(J(F)\) \emph{are}
the members of \(F\), so a condition on the graph is automatically a
condition on the family: we say \(F\) is a clique-sum when \(J(F)\) is
one, and then its universal pair is a pair of \emph{members} of \(F\).
Nothing else is meant by the phrase. Theorem 7.6 says that for shifted
families this graph condition and the coordinates of Section 7.5 pick
out the same objects. Lemma 7.4 below is pure graph theory: it holds of
any graph with that structure, knows nothing about families, and does
not use shiftedness.

\textbf{Lemma 7.4.} \emph{Let \((G, \{u,v\})\) be a clique-sum. Then
\(G\) has no Hamilton \(u\)--\(v\) path.}

This is the standard obstruction and nothing more: \textbf{if
\(G − \{u,v\}\) is disconnected then \(G\) has no Hamilton \(u\)--\(v\)
path}, since deleting the two ends of such a path leaves a single path,
which cannot cover a disconnected set.

Every pair \emph{other} than \(\{u,v\}\) does have a Hamilton path
between it, but that is not needed here and is not proved here: for a
shifted family it falls out of Theorem 7.7, and we record it as a remark
after Corollary 7.8 rather than proving it twice.

\emph{Proof.} Deleting \(\{u,v\}\) leaves the two cliques with no edge
between them, which is disconnected. The interior of a Hamilton
\(u\)--\(v\) path is a path covering everything else, and a path is
connected, so it cannot cover a disconnected set. \(\blacksquare\)

\textbf{Lemma 7.5.} \emph{Let \(F\) be shifted of rank \(k\) and suppose
\(J(F)\) is a clique-sum. Then its universal pair is \(\{m_1, m_2\}\).}

This is the first place the family is doing work rather than the graph:
\(m_1\) and \(m_2\) are defined by the Gale order, and the proof runs on
decrements.

\emph{Proof.} Let \(U\) be the universal pair and \(P, Q\) the two
parts, both nonempty.

\emph{Step 1: \(m_1 ∈ U\).} Suppose not, say \(m_1 ∈ P\). Every
\(W ∈ Q\) is non-adjacent to \(m_1\), so \(|W ∖ m_1| ≥ 2\). Choose \(W\)
Gale-minimal in \(Q\). Every decrement of \(W\) lies in \(F\), is
Gale-below \(W\), and is \textbf{adjacent} to \(W\); adjacency places it
in \(U ∪ Q\), and minimality excludes \(Q\), so all decrements lie in
\(U\), which has two members. But \(|W ∖ m_1| ≥ 2\) puts \(W\) in the
range of Lemma 7.1(d), which gives at least four decrements.
Contradiction.

\emph{Step 2: \(m_2 ∈ U\).} Suppose not, so \(U = \{m_1, u\}\) with
\(u ≠ m_2\), and \(m_2\) lies in a clique \(R\). Choose \(W\)
Gale-minimal in the other clique, which is nonempty. As in Step 1, all
decrements of \(W\) lie in \(U\). If \(|W ∖ m_1| ≥ 2\) then Lemma 7.1(d)
gives at least four, a contradiction, so \(W = [k] − i + a\). Here
\(a − i = 1\) would make \(W = m_2\), which lies in \(R\); and
\(a − i ≥ 3\) gives at least three decrements. So \(a − i = 2\), and
Lemma 7.1(c) makes the decrement set exactly \(\{m_1, m_2\}\), forcing
\(m_2 ∈ U\). \(\blacksquare\)

\hypertarget{the-classification-1}{%
\subsection{7.5 The classification}\label{the-classification-1}}

\textbf{What is about to happen.} Section 7.4 exhibited a shape that
fails --- the clique-sum --- and that is a condition on the \emph{graph}
\(J(F)\). What one wants is a condition on the \emph{family}, checkable
from \(F\) itself, and this subsection supplies it: the Y-families
\(F(k; r, s)\), three integer parameters and nothing else. Theorem 7.6
is the equivalence between the two.

For \(k ≥ 2\) and parameters \(r\) and \(s\) with \(k+2 ≤ r ≤ n\) and
\(1 ≤ s ≤ k−1\), define

\[
\begin{gathered}
F(k; r, s) = \{m_1, m_2\} ∪ P ∪ Q, \\
P = \{ [k−1] ∪ \{a\} : k+2 ≤ a ≤ r \}, \quad Q = \{ [k+1] ∖ \{i\} : s ≤ i ≤ k−1 \}.
\end{gathered}
\]

\textbf{Said in one line, the family is a down-closure.} \(F(k; r, s)\)
is exactly the set of \(k\)-sets Gale-below one or other of

\[
[k−1] ∪ \{r\} \quad \text{ and } \quad [k+1] ∖ \{s\},
\]

the two arm tops: the short prefix with one far element adjoined, and
the long prefix with one small element removed. Down-closing those two
produces both arms and the stem beneath them, so the whole family is
named by two of its members.

That form also explains the two parameter ranges, which otherwise look
arbitrary. Push either parameter toward the middle and both formulas run
into the \textbf{same} two sets:

\[
[k−1] ∪ \{k\} \quad = [k+1] ∖ \{k+1\} = m_1, \quad [k−1] ∪ \{k+1\} = [k+1] ∖ \{k\} = m_2.
\]

So \(k+2 ≤ r\) and \(s ≤ k−1\) are one condition stated twice --- each
generator must lie strictly beyond the stem, or it \emph{is} the stem,
and the arm it was supposed to top is empty.

\textbf{How many there are.} Because the pair \((r, s)\) determines the
family and distinct pairs give distinct families, the two parameter
ranges count the class outright: the ground \([n]\) carries exactly

\[
(n − k − 1)(k − 1)
\]

Y-families of rank \(k\). Summing over the ranks, and writing
\(m = n − 2\),

\[
Σ_{k=2}^{n−1} (n − k − 1)(k − 1) = Σ_{j=1}^{m} j(m − j) = m(m² − 1)/6 = C(n−1, 3),
\]

so all ranks together contribute \(C(n−1, 3)\) of them. The exception
class is therefore not just explicit but cubic in the ground: a Y-family
is named by two integers, and one can list every Y-family on a given
ground without generating a single shifted family.

The two arms are built in opposite directions from the two prefixes on
either side of \([k]\), and that is the thing to hold on to.
\textbf{\(P\) is grown upward from \([k−1]\)}: take the short prefix and
adjoin one element from beyond \([k+1]\), letting that element run out
to \(r\). \textbf{\(Q\) is cut down from \([k+1]\)}: take the long
prefix and delete one element from inside it, letting the deleted
element run down to \(s\). So \(P\)'s members share a common low core
and each carries one far element, while \(Q\)'s members live entirely
inside \([k+1]\) and each is missing one small element.

We call \(P\) the \textbf{outer arm}, \(Q\) the \textbf{inner arm}, and
\(F(k; r, s)\) a \textbf{Y-family}. Neither arm meets \(\{m_1, m_2\}\):
the parameter ranges put every outer member beyond \([k+1]\) and every
inner member inside it, so an arm is always disjoint from the universal
pair Theorem 7.6 identifies. All three names are ours and are not
standard. The first two say where the members sit relative to \([k+1]\):
every member of \(P\) contains an element beyond \([k+1]\), reaching
\textbf{out} as far as \(r\), while every member of \(Q\) lies
\textbf{inside} \([k+1]\), cut in as deep as \(s\). So the two
parameters read off the names --- \(r\) is how far the outer arm
reaches, \(s\) is how deep the inner arm cuts.

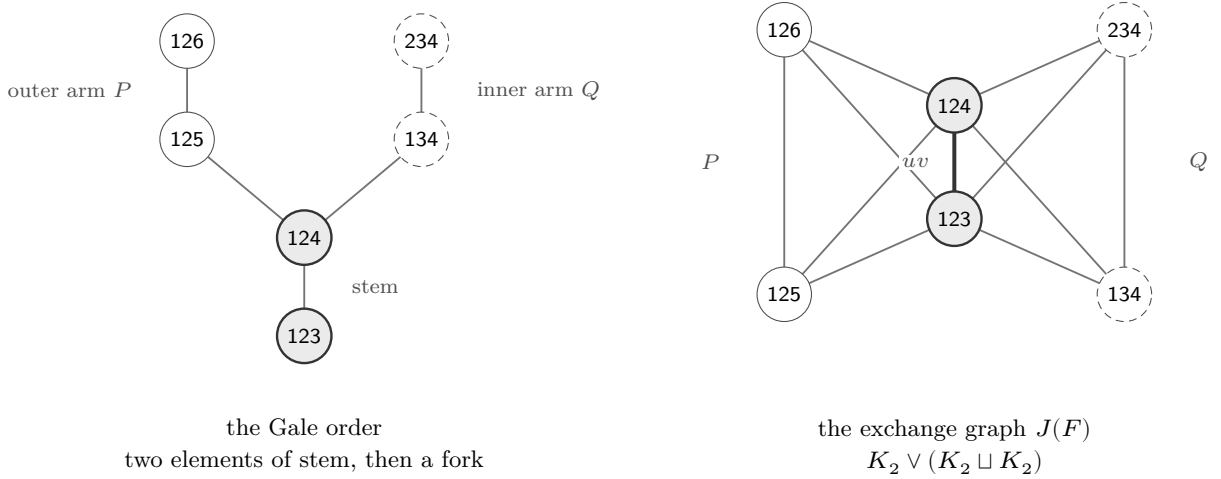
\begin{figure}[tbp]
\centering
\begin{tikzpicture}[x=1cm,y=1cm,
  ord/.style={black!55,line width=0.7pt},
  gedge/.style={black!55,line width=0.7pt},
  stem/.style={circle,draw=black!80,fill=black!8,line width=0.9pt,minimum size=7.2mm,inner sep=0pt,font=\scriptsize},
  outerarm/.style={circle,draw=black!70,fill=white,minimum size=7.2mm,inner sep=0pt,font=\scriptsize},
  innerarm/.style={circle,draw=black!70,fill=white,densely dashed,minimum size=7.2mm,inner sep=0pt,font=\scriptsize},
  armlabel/.style={font=\scriptsize,text=black!70},
  paneltitle/.style={font=\footnotesize,anchor=north,align=center},
]
\draw[ord] (0,0) -- (0,1.3);
\draw[ord] (0,1.3) -- (-1.55,2.6);
\draw[ord] (0,1.3) -- (1.55,2.6);
\draw[ord] (-1.55,2.6) -- (-1.55,3.9);
\draw[ord] (1.55,2.6) -- (1.55,3.9);
\node[stem] at (0,0) {$\mathsf{123}$};
\node[stem] at (0,1.3) {$\mathsf{124}$};
\node[outerarm] at (-1.55,2.6) {$\mathsf{125}$};
\node[outerarm] at (-1.55,3.9) {$\mathsf{126}$};
\node[innerarm] at (1.55,3.9) {$\mathsf{234}$};
\node[innerarm] at (1.55,2.6) {$\mathsf{134}$};
\node[armlabel,anchor=east] at (-2.15,3.25) {outer arm $P$};
\node[armlabel,anchor=west] at (2.15,3.25) {inner arm $Q$};
\node[armlabel,anchor=west] at (0.5,0.65) {stem};
\node[paneltitle] at (0,-0.95) {the Gale order\\[1pt]two elements of stem, then a fork};
\draw[gedge] (10.85,4.05) -- (10.85,0.55);
\draw[gedge] (8.6,3.05) -- (6.35,0.55);
\draw[gedge] (8.6,1.55) -- (8.6,3.05);
\draw[gedge] (6.35,4.05) -- (8.6,3.05);
\draw[gedge] (8.6,1.55) -- (10.85,4.05);
\draw[gedge] (10.85,4.05) -- (8.6,3.05);
\draw[gedge] (6.35,4.05) -- (8.6,1.55);
\draw[gedge] (6.35,4.05) -- (6.35,0.55);
\draw[gedge] (8.6,1.55) -- (10.85,0.55);
\draw[gedge] (8.6,3.05) -- (10.85,0.55);
\draw[gedge] (8.6,1.55) -- (6.35,0.55);
\draw[black!80,line width=1.5pt] (8.6,1.55) -- (8.6,3.05);
\node[stem] at (8.6,1.55) {$\mathsf{123}$};
\node[stem] at (8.6,3.05) {$\mathsf{124}$};
\node[outerarm] at (6.35,0.55) {$\mathsf{125}$};
\node[outerarm] at (6.35,4.05) {$\mathsf{126}$};
\node[innerarm] at (10.85,4.05) {$\mathsf{234}$};
\node[innerarm] at (10.85,0.55) {$\mathsf{134}$};
\node[armlabel,anchor=east] at (5.63,2.30) {$P$};
\node[armlabel,anchor=west] at (11.57,2.30) {$Q$};
\node[armlabel,anchor=east,fill=white,inner sep=1pt] at (8.3,2.30) {$uv$};
\node[paneltitle] at (8.6,-0.95) {the exchange graph $J(F)$\\[1pt]$K_2 \vee (K_{2} \sqcup K_{2})$};
\end{tikzpicture}
\caption{The Y-family $F(3; 6, 1)$, drawn twice. Left, its Gale order: the two Gale-least members $m_1 = \mathsf{123}$ and $m_2 = \mathsf{124}$ form a stem below everything, and above them the order forks into the outer chain $P$ and the inner chain $Q$, which are never comparable to one another. Right, the same family's exchange graph $J(F)$: the shared edge $uv = m_1m_2$ is drawn heavy, both of its endpoints are adjacent to everything, and the two arms are cliques with no edge between them. Deleting $m_1$ and $m_2$ leaves the two arms with nothing joining them, which is the clique-sum hypothesis of Lemma 7.4; that lemma then denies a Hamilton $m_1$--$m_2$ path, and that is why this family is the exception of Theorem 7.7.}
\label{fig:yfamily}
\end{figure}

\textbf{The third name records the shape of the Gale order, and it is
worth seeing.} Each arm is a \emph{chain}: raising \(a\) moves an outer
member up, lowering \(i\) moves an inner member up, and within an arm
any two members are comparable. The two Gale-least members sit below
everything, \(m_1 ≺ m_2 ≺ X\) for every other member \(X\). And the arms
are never comparable to one another --- no outer member is ever above or
below an inner member. So the order is a two-element stem that forks
into two chains which never rejoin --- the shape the name records. Drawn
as a Hasse diagram, \(m_1\) sits at the bottom with \(m_2\) directly
above it, and above \(m_2\) the order divides once and for all into the
outer chain on one side and the inner chain on the other, with no
relation across the divide. Two elements of stem, then a fork: a \(Y\).
It also explains the clique-sum of Section 7.4 exactly: each arm
together with the stem is a clique of \(J(F)\), the two cliques meet in
the stem and nowhere else, and the stem is the shared edge. The two
names are one object seen from its order and from its graph, and Figure
\ref{fig:yfamily} draws the smallest instance in which both arms are
large enough to see.

The arms are not comparable in the Gale order: no member of one is ever
above or below a member of the other. Corollary 7.9 records the instance
of this for the two arm tops, and it is this cross-arm incomparability
that makes the Y-family a clique-sum in the first place.

\textbf{Theorem 7.6.} \emph{A shifted family is a clique-sum if and only
if it is a Y-family. The universal pair of \(F(k; r, s)\) is
\(\{m_1, m_2\} = \{[k], [k−1] ∪ \{k+1\}\}\), and its two parts are \(P\)
and \(Q\).}

\emph{Proof.} Let \(F = F(k; r, s)\), so \(m_1 △ m_2 = \{k, k+1\}\).
Each member of \(P\) differs from \(m_1\) in \(\{k, a\}\) and from
\(m_2\) in \(\{k+1, a\}\); each member of \(Q\) differs from \(m_1\) in
\(\{i, k+1\}\) and from \(m_2\) in \(\{i, k\}\). So \(m_1\) and \(m_2\)
are adjacent to everything. \(P\) is a clique, its members sharing
\([k−1]\); \(Q\) is a clique, its members lying in \([k+1]\). A pair
drawn one from each arm has symmetric difference \(\{i, a, k, k+1\}\) of
size four, since \(i ≤ k−1 < k < k+1 < k+2 ≤ a\), so no edge joins the
arms. Both arms are nonempty, so no arm member is dominating, and
\(J(F) = K_2 ∨ (K_{|P|} ⊔ K_{|Q|})\). That \(F(k; r, s)\) is shifted
follows from the census: the decrements of \([k−1] ∪ \{a\}\) are
\(m_1\), \(m_2\) and the smaller outer members, and those of
\([k+1] ∖ \{i\}\) are the inner members of larger index together with
\(m_1\) and \(m_2\).

Conversely let \(F\) be a clique-sum, with universal pair
\(\{m_1, m_2\}\) by Lemma 7.5. Any other member \(X\) is adjacent to
\(m_1\), so \(X = [k] − i + a\) with \(i ≤ k < a\), and \(X\) is
adjacent to \(m_2\). Now \(X △ m_2 ⊆ \{i, a, k, k+1\}\), of size four
unless \(i = k\) or \(a = k+1\). Taking \(i = k\) gives
\(X = [k−1] ∪ \{a\}\) with \(a ≥ k+2\), an outer member; taking
\(a = k+1\) gives \(X = [k+1] ∖ \{i\}\) with \(i ≤ k−1\), an inner
member. Outer members are pairwise adjacent and inner members are
pairwise adjacent, while an outer member and an inner member are at
symmetric difference four. So each of the two sets is a clique and lies
inside a single component of \(J(F) − \{m_1, m_2\}\); that graph is
\(K_p ⊔ K_q\) with both parts nonempty, and neither set can be all of
it, since it is disconnected and they are connected. Hence both sets are
nonempty and they are exactly \(P\) and \(Q\). Down-closure makes the
\(a\) contiguous upward from \(k+2\) and the \(i\) contiguous downward
from \(k−1\), and both cliques are nonempty since \(p, q ≥ 1\); a
nonempty inner clique needs \(s ≤ k−1\), which forces \(k ≥ 2\).
\(\blacksquare\)

\hypertarget{the-main-theorem}{%
\subsection{7.6 The main theorem}\label{the-main-theorem}}

\textbf{What it says, in words.} For a shifted family the Johnson graph
is as traversable as it could be: between any two members there is a
Hamilton path. There is one family of exceptions, and in those the only
unreachable pair is \(\{m_1, m_2\}\) --- the two forced smallest members
that Section 7.3 showed every shifted family has. So the failure is not
a scattering of awkward pairs to be checked one at a time. It is a
single shape, occurring at a single pair, identifiable from three
integers. That is what makes it something a construction can be steered
around, which is what Section 8 does.

\textbf{Theorem 7.7.} \emph{Let \(F\) be shifted of rank \(k\) on
\([n]\) and let \(A ≠ B\) be members. Then \(J(F)\) has a Hamilton
\(A\)--\(B\) path unless \(F\) is a Y-family and \(\{A, B\}\) is its
universal pair, in which case it has none.}

\textbf{Corollary 7.8.} \emph{If \(F\) is shifted and is not a Y-family,
then \(J(F)\) is Hamilton-connected.}

\emph{Proof.} If \(F\) is not a Y-family then the exception of Theorem
7.7 cannot arise for any pair, so every two distinct members are joined
by a Hamilton path. \(\blacksquare\)

Shiftedness is inherited from Theorem 7.7 and cannot be dropped from
this corollary either. Without it the statement is false:
\(F = \{12, 34\}\) is not a Y-family --- a Y-family has at least four
members --- and \(J(F)\) is two isolated vertices, so it is not
Hamilton-connected.

\emph{The other pairs.} At a clique-sum every pair \textbf{other} than
the universal pair does have a Hamilton path between it. For a shifted
family this is immediate from Theorem 7.7: such a family is a Y-family
by Theorem 7.6, and Theorem 7.7 supplies a Hamilton path for every pair
that is not the Y-family's universal pair. The formalization records it
exactly this way, as a consequence of the main theorem rather than as a
separate construction.

\textbf{Corollary 7.9.} \emph{Every Y-family has two incomparable
Gale-maximal members, the tops of its two arms. A shifted family with a
Gale-greatest member is therefore never a Y-family, and its Johnson
graph is Hamilton-connected.}

\emph{Proof.} Write the arms as in Section 7.5: the outer members are
\([k−1] ∪ \{a\}\) for \(k+2 ≤ a ≤ r\), and the inner members are
\([k+1] ∖ \{i\}\) for \(s ≤ i ≤ k−1\), so \(a\) and \(i\) are the
indices running along the two arms. A larger \(a\) raises an outer
member in the Gale order and a smaller \(i\) raises an inner member, so
the tops of the two arms are \([k−1] ∪ \{r\}\) and \([k+1] ∖ \{s\}\).
Sorted, these are \((1, …, k−1, r)\) and \((1, …, s−1, s+1, …, k+1)\).
At position \(s\) the first has entry \(s\) and the second has \(s+1\);
at position \(k\) the first has \(r ≥ k+2\) and the second has \(k+1\).
So neither dominates the other and they are incomparable. Each is also
maximal in \(F\): it is above the rest of its own arm and above \(m_1\)
and \(m_2\), and the two displayed positions show it is below no member
of the other arm. A family with a Gale-greatest member \(M\) has every
member below \(M\), so \(M\) is its only maximal member and no two
incomparable maximal members exist; such a family is therefore not a
Y-family, and Corollary 7.8 applies. \(\blacksquare\)

Theorem 7.7 has two halves and they are of very different weight. Its
\textbf{obstruction half} --- at the exception there is no Hamilton path
--- is Lemma 7.4 together with Theorem 7.6, which identifies the
Y-families as exactly the clique-sums; that is the short half, and it is
already done. Its \textbf{existence half} --- that everywhere else a
Hamilton path does exist --- is the substantial one, and Section 7.8
proves it.

\hypertarget{duality}{%
\subsection{7.7 Duality}\label{duality}}

The case analysis below is cut in half by an \textbf{involution} --- a
map that is its own inverse, so applying it twice returns every set to
itself --- which we record first because it also explains why the
exception is symmetric in a way the definition does not make evident.

Fix \(n\) and \(k\) with \(1 ≤ k ≤ n−1\). For \textbf{any} subset
\(X ⊆ [n]\) put

\[
X^* = \{ n+1−x : x ∈ [n] ∖ X \},
\]

that is, complement and then reverse the linear order on the ground; a
\(k\)-subset becomes an \((n−k)\)-subset. For any family \(F\) of
subsets write \(F^* = \{X^* : X ∈ F\}\). Nothing here needs the members
to be equicardinal, which is what lets the notation reach the barred
sets \(S^e\) and \(S^{¬e}\) of Section 7.8.

\textbf{Lemma 7.10 (duality).} \emph{The map \(X ↦ X^*\) is an
involution carrying \(k\)-subsets of \([n]\) onto \((n−k)\)-subsets,
and:}

\emph{(i) \(X ≼ Y\) if and only if \(X^* ≼ Y^*\);} \emph{(ii)
\(|X^* △ Y^*| = |X △ Y|\);} \emph{(iii) for \(e ∈ [n]\) and
\(e^* = n+1−e\), \(e ∈ X\) if and only if \(e^* ∉ X^*\);} \emph{(iv)
\(X ⊆ Y\) if and only if \(Y^* ⊆ X^*\).}

\emph{Consequently, if \(F\) is shifted of rank \(k\) then \(F^*\) is
shifted of rank \(n−k\), the map is an isomorphism \(J(F) → J(F^*)\), it
carries the members of \(F\) containing \(e\) onto the members of
\(F^*\) avoiding \(e^*\) and conversely, and \(F\) is a Y-family with
universal pair \(\{u, v\}\) if and only if \(F^*\) is, with universal
pair \(\{u^*, v^*\}\) and the two arms exchanged: the inner arm of \(F\)
maps onto the outer arm of \(F^*\).}

\emph{Proof.} Write \(σ(x) = n+1−x\), so \(X^* = σ(X^c)\), where the
complement is taken in \([n]\). Both \(σ\) and complementation are
involutions and they commute, so \(X ↦ X^*\) is one, and it changes the
size from \(k\) to \(n−k\).

\begin{enumerate}
\def\labelenumi{(\roman{enumi})}
\tightlist
\item
  For a set \(Z\) and a threshold \(t ∈ [n]\) write
  \(c_Z(t) = |Z ∩ [t]|\). Since \(x_i ≤ t\) exactly when \(c_X(t) ≥ i\),
  two \(k\)-sets satisfy \(X ≼ Y\) if and only if \(c_X(t) ≥ c_Y(t)\)
  for every \(t\): the Gale order is ``at least as many elements below
  every threshold''. Now count \(X^*\) below \(t\):
\end{enumerate}

\[
\begin{gathered}
c_{X^*}(t) = |\{ σ(x) : x ∉ X, σ(x) ≤ t \}| = |\{ x ∉ X : x ≥ n+1−t \}| \\
= t − |X ∩ [n+1−t, n]| = t − (k − c_X(n−t)),
\end{gathered}
\]

so \(c_{X^*}(t) = t − k + c_X(n−t)\). The term \(t − k\) does not depend
on the set, so \(c_{X^*}(t) ≥ c_{Y^*}(t)\) is equivalent to
\(c_X(n−t) ≥ c_Y(n−t)\). As \(t\) runs over \([n]\), \(n−t\) runs over
\(\{0, 1, …, n−1\}\), and both counts agree at \(0\) and at \(n\), so
the two families of inequalities are the same family. Hence
\(X^* ≼ Y^*\) if and only if \(X ≼ Y\).

\begin{enumerate}
\def\labelenumi{(\roman{enumi})}
\setcounter{enumi}{1}
\item
  Complementation preserves symmetric differences,
  \(X^c △ Y^c = X △ Y\), and \(σ\) is a bijection, so
  \(X^* △ Y^* = σ(X △ Y)\), which has the same size.
\item
  \(e ∈ X\) if and only if \(e ∉ X^c\), if and only if
  \(σ(e) = e^* ∉ σ(X^c) = X^*\).
\item
  Complementation reverses inclusion and \(σ\) preserves it.
\end{enumerate}

For the consequences: a shifted family is a nonempty down-set of \(≼\),
and (i) says \(X ↦ X^*\) is an order isomorphism onto the
\((n−k)\)-subsets, so it carries down-sets to down-sets; (ii) makes it
an isomorphism of Johnson graphs; (iii) is the statement about \(e\).
Being a Y-family is a property of the family and its order, not of
\(J(F)\) alone --- \(\{12, 13, 14, 23\}\) and \(\{14, 23, 24, 34\}\)
have isomorphic Johnson graphs and only the first is a Y-family. What
transfers is the graph-level property: by Theorem 7.6, among
\textbf{shifted} families being a Y-family is equivalent to \(J(F)\)
being a clique-sum, and both \(F\) and \(F^*\) are shifted, so
Y-familyness and the universal pair carry along the isomorphism. For the
arms, let \(\{u,v\}\) be the universal pair. A member other than \(u\)
and \(v\) lies in the inner arm exactly when \(X ⊆ u ∪ v\); the
qualification matters, since \(u ⊆ u ∪ v\) holds trivially. By (iv) this
holds if and only if \(X^* ⊇ (u ∪ v)^*\). Complementation turns the
union into an intersection, so \((u ∪ v)^* = σ(u^c ∩ v^c) = u^* ∩ v^*\),
and \(X^* ⊇ u^* ∩ v^*\) is membership in the outer arm of \(F^*\).
\(\blacksquare\)

The hypothesis \(k ≤ n−1\) is automatic wherever we use this: two
distinct members of \(F\) force \(k < n\), since \([n]\) has only one
\(n\)-subset. The rank changes under the duality, from \(k\) to \(n−k\),
so Lemma 7.10 is not a symmetry of a single family. It is used only to
transport one completed case of the analysis below onto another, never
inside an induction, so no circularity arises: Theorem 7.11 is proved by
a direct case analysis on its parameters, and the duality permutes those
parameters.

\hypertarget{the-induction-1}{%
\subsection{7.8 The induction}\label{the-induction-1}}

The section has two halves. The first sets up the split and says what a
crossing has to avoid; the second proves that a usable one always
exists. A reader who wants only the statement can stop after Theorem
7.11.

\hypertarget{the-split-and-what-a-usable-crossing-must-avoid}{%
\subsubsection{7.8.1 The split, and what a usable crossing must
avoid}\label{the-split-and-what-a-usable-crossing-must-avoid}}

\textbf{What actually happens here, before the notation arrives.} Split
the family at an element \(e\) into the members containing it and the
members avoiding it. Both live on a smaller ground, so the induction
hypothesis gives a Hamilton path inside each; splice the two across a
single edge joining them and you have the path you wanted. Three things
can go wrong. The recursion must land on shifted families, which the
prefix count below gives --- for the avoiding layer \(F^{¬e}\) that is
the layer itself, but on the containing side it is the \textbf{link}
\(L\), the layer with \(e\) deleted from every member. The containing
layer \(F^e\) need not be closed under decrement, which is exactly why
the induction runs through \(L\) and never through \(F^e\). Each side's
own induction may fail at that side's exceptional pair, so each bars at
most one member beyond the endpoint, and the splice must avoid both. And
a joining edge must survive those bars --- that is Theorem 7.11, the
crossing lemma, and it is the one hard thing here. \textbf{Given it, the
assembly is three lines}; almost everything after it is the proof of
7.11 itself, a case analysis that duality cuts in half.

We induct on \(n\), the size of the ground --- which is what ``ground
order \(n\)'' means everywhere in this paper, Section 7.7's reversal of
the ground's linear order being the one place the phrase could be read
the other way. Fix \(e ∈ A △ B\). Exactly one of \(A\), \(B\) contains
it, and we name that one \(A\) --- a Hamilton \(A\)--\(B\) path is a
Hamilton \(B\)--\(A\) path, so swapping the two labels costs nothing.
Split

\[
F^e = \{ W ∈ F : e ∈ W \}, \quad F^{¬e} = \{ W ∈ F : e ∉ W \}, \quad L = \{ W − e : W ∈ F^e \}.
\]

\(F^{¬e}\) is shifted on \([n] ∖ \{e\}\) of rank \(k\), and \(L\) is
shifted there of rank \(k−1\); both follow from the \textbf{prefix
count} \(c_Z(t) = |Z ∩ [t]|\) of Section 7.1 --- how many members of
\(Z\) lie in the prefix \([t]\) of the ground --- which satisfies
\(c_{W ∪ \{e\}}(t) = c_W(t) + [e ≤ t]\), so inserting \(e\) shifts every
count by the same amount and Gale domination transfers unchanged. The
map \(W ↦ W − e\) is an isomorphism \(J(F^e) → J(L)\), being a bijection
that preserves symmetric differences. We call \(F^e\) and \(F^{¬e}\) the
\textbf{layers} of \(F\) at \(e\) --- the \textbf{containing layer} and
the \textbf{avoiding layer}, after whether their members contain \(e\)
--- and \(L\) the \textbf{link} --- the name is the standard one for
this operation on a simplicial complex, where the link of a face is what
remains of the sets containing it once that face is deleted, which is
exactly \(W ↦ W − e\) applied to \(F^e\). An edge of \(J(F)\) joining
the two layers is a \textbf{crossing}. The superscript names the element
the layer is taken at, and it earns its place: \(e\) is not fixed below.
The proof runs through the regimes \(e = k−1\), \(k\), \(k+1\), \(k+2\),
and Section 7.7 carries the layers of \((F, e)\) onto those of
\((F^*, e^*)\) for \(e^* = n+1−e\), so two different splitting elements
are in play at once there.

The link and the avoiding layer live on a ground of size \(n − 1\), so
the induction on \(n\) is well founded. The step is needed only where
\(J(F)\) is not already complete, and a complete graph has a Hamilton
path between any two of its vertices, so the following cases are settled
outright and stop the recursion. First \(|F| ≤ 3\), where \(J(F)\) is
complete: trivially for \(|F| ≤ 1\), and by Lemma 7.2 otherwise. No
Y-family arises there, since a Y-family has \(|F| ≥ 4\). Second
\(k ≤ 1\), where any two members meet in the empty set, so \(J(F)\) is
complete again. Every ground small enough to have no smaller ground
below it falls under one of these.

\begin{figure}[tbp]
\centering
\begin{tikzpicture}[x=1cm,y=1cm,
  mem/.style={circle,draw=black!72,fill=white,minimum size=7.4mm,inner sep=0pt,font=\scriptsize},
  endpt/.style={circle,draw=black!85,fill=black!10,line width=0.9pt,minimum size=7.4mm,inner sep=0pt,font=\scriptsize},
  linkmem/.style={circle,draw=black!55,fill=white,densely dotted,minimum size=7.0mm,inner sep=0pt,font=\scriptsize},
  pathedge/.style={black!60,line width=0.8pt},
  crossedge/.style={black!85,line width=1.5pt},
  layerbox/.style={rounded corners=3pt,draw=black!35,densely dashed},
  boxtitle/.style={font=\scriptsize,align=center,text=black!75},
  ann/.style={font=\scriptsize,text=black!70},
]
\draw[layerbox] (-2.27,-5.27) rectangle (2.27,-4.03);
\node[boxtitle,anchor=north] at (0,-5.39) {the link $L = F^{e} - e$, shifted of rank $2$\\on the smaller ground --- this is what the induction uses};
\draw[layerbox] (-0.72,-3.12) rectangle (0.72,0.62);
\node[boxtitle,anchor=south] at (0,0.72) {containing layer $F^{e}$};
\draw[layerbox] (3.88,-5.12) rectangle (5.32,-1.38);
\node[boxtitle,anchor=south] at (4.6,-1.28) {avoiding layer $F^{\neg e}$};
\draw[pathedge] (0,-0) -- (0,-1.25);
\draw[pathedge] (0,-1.25) -- (0,-2.5);
\draw[pathedge] (4.6,-2) -- (4.6,-3.25);
\draw[pathedge] (4.6,-3.25) -- (4.6,-4.5);
\draw[crossedge] (0,-2.5) -- (4.6,-2);
\draw[black!45,densely dotted,-{Stealth[length=4pt]},line width=0.7pt] (0,-3) -- (0,-3.93);
\node[ann,anchor=west] at (0.22,-3.465) {delete $e$};
\node[endpt] at (0,-0) {$\mathsf{123}$};
\node[endpt] at (4.6,-4.5) {$\mathsf{124}$};
\node[mem] at (4.6,-2) {$\mathsf{125}$};
\node[mem] at (4.6,-3.25) {$\mathsf{126}$};
\node[mem] at (0,-1.25) {$\mathsf{134}$};
\node[mem] at (0,-2.5) {$\mathsf{135}$};
\node[linkmem] at (-1.55,-4.65) {$\mathsf{12}$};
\node[linkmem] at (0,-4.65) {$\mathsf{14}$};
\node[linkmem] at (1.55,-4.65) {$\mathsf{15}$};
\node[ann,anchor=east] at (-0.62,-0) {$A$};
\node[ann,anchor=west] at (5.22,-4.5) {$B$};
\node[ann,anchor=south,fill=white,inner sep=1.5pt] at (2.3,-2.19) {the crossing};
\end{tikzpicture}
\caption{The split of Section 7.8.1, carried out on $F = \{\mathsf{123}, \mathsf{124}, \mathsf{125}, \mathsf{126}, \mathsf{134}, \mathsf{135}\}$, shifted of rank $3$ on $[6]$, at $e = 3$. The containing layer $F^{e}$ and the avoiding layer $F^{\neg e}$ partition $F$. The induction is not applied to $F^{e}$, which is not closed under decrement, but to the link $L = F^{e} - e$, shifted of rank $2$ on the ground with $e$ removed; $W \mapsto W - e$ is an isomorphism $J(F^{e}) \to J(L)$. A Hamilton path of each layer, joined across the single heavy crossing $\mathsf{135}$--$\mathsf{125}$, splices to the Hamilton $A$--$B$ path $\mathsf{123} \to \mathsf{134} \to \mathsf{135} \to \mathsf{125} \to \mathsf{126} \to \mathsf{124}$ of $J(F)$. What the induction must supply is the crossing, and that is Theorem 7.11.}
\label{fig:layersplit}
\end{figure}
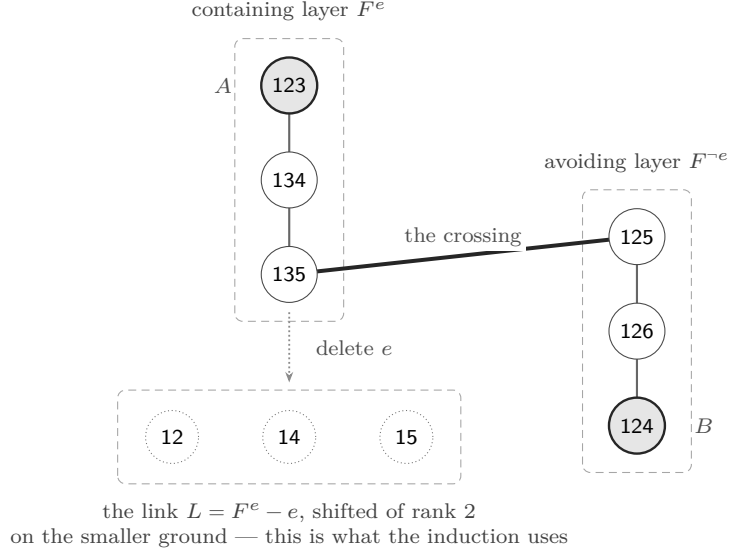

The assembly is immediate once a crossing is available. A Hamilton
\(A\)--\(W\) path of \(J(F^e)\), the edge \(WZ\), and a Hamilton
\(Z\)--\(B\) path of \(J(F^{¬e})\) concatenate to a Hamilton
\(A\)--\(B\) path of \(J(F)\), since the layers partition \(F\) and
adjacency restricts correctly to each. What the induction must supply is
the crossing, with both layer pairs kept off their own exceptions.
Figure \ref{fig:layersplit} carries the whole split out on one family.

\textbf{The barred sets.} Name the two forbidden partners. If \(L\) is a
Y-family whose universal pair contains \(A − e\), let \(W_f\) be the
member of \(F^e\) lying over the other member of that pair; otherwise
\(W_f\) is undefined. If \(F^{¬e}\) is a Y-family whose universal pair
contains \(B\), let \(Z_f\) be the other member of that pair; otherwise
\(Z_f\) is undefined. Write

\[
S^e = \{A\} \text{ or } \{A, W_f\}, \quad S^{¬e} = \{B\} \text{ or } \{B, Z_f\},
\]

according to whether the second member is defined. The superscript is
the layer it bars: \(S^e ⊆ F^e\) and \(S^{¬e} ⊆ F^{¬e}\). The first is
stated for the \textbf{link} and lifted back, not for \(F^e\), which
need not be closed under decrement: a decrement at \(e\) leaves the
layer. (It can happen to be closed; at \(e = 1\) nothing can be
decremented out. But not in general.) Since \(W ↦ W − e\) is an
isomorphism \(J(F^e) → J(L)\), the condition on \(L\) may equally be
read on \(J(F^e)\): \(S^e\) has two members exactly when
\(J(F^e) ≅ K_2 ∨ (K_p ⊔ K_q)\) with \(A\) in the dominating pair, and
\(S^e\) is then that pair. We use the second reading when applying Lemma
7.10, since it refers to \(F^e\) and not to the link.

Throughout, a barred set's superscript names the element its layer is
taken at, so on the dual side that element is \(e^*\) and never \(e\).
Only the containing-layer set \(S^{e^*}(F^*)\) goes through the dual's
own link and is re-read on its containing layer, exactly as above;
\(S^{¬e^*}(F^*)\) is defined directly on the avoiding layer, which needs
no link. Duality exchanges the two layers and moves the splitting
element, so the containing layer at \(e\) becomes the \textbf{avoiding}
layer at \(e^*\), and conversely:

\[
(F^e)^* = (F^*)^{¬e^*}, \quad (F^{¬e})^* = (F^*)^{e^*}.
\]

Lemma 7.10 carries \(J(F^e)\) isomorphically onto \(J((F^*)^{¬e^*})\),
sending \(A\) to \(A^*\), and \(J(F^{¬e})\) onto \(J((F^*)^{e^*})\),
sending \(B\) to \(B^*\). Having the clique-sum shape with a named
member dominating is preserved by an isomorphism, so

\[
S^e(F, A)^* = S^{¬e^*}(F^*, A^*) \quad \text{ and } \quad S^{¬e}(F, B)^* = S^{e^*}(F^*, B^*).
\]

This is the seam at which the duality could fail, and it is the one
place where the two readings of a barred set must be reconciled on both
sides at once.

A barred set has two members only when its layer has the clique-sum
shape with the endpoint dominating, and a Y-family has at least four
members, so

\[
\begin{gathered}
F^e ∖ S^e ≠ ∅ \quad \text{ whenever } |F^e| ≥ 2, \\
F^{¬e} ∖ S^{¬e} ≠ ∅ \quad \text{ whenever } |F^{¬e}| ≥ 2. \quad (∗)
\end{gathered}
\]

\textbf{Theorem 7.11 (crossing lemma).} \emph{Let \(F\) be shifted and
let \(A ≠ B\) be members with \((F, \{A,B\})\) not the exception of
Theorem 7.7. Then for \textbf{every} \(e ∈ A △ B\) with \(e ∈ A\) there
are \(W ∈ F^e\) and \(Z = W − e + ξ ∈ F^{¬e}\), for some \(ξ ≠ e\) with
\(ξ ∉ W\), such that \(W = A\) if \(|F^e| = 1\) and otherwise
\(W ∉ S^e\), and \(Z = B\) if \(|F^{¬e}| = 1\) and otherwise
\(Z ∉ S^{¬e}\).}

The exclusions are the content. The induction hypothesis fails inside a
layer only at that layer's own exceptional pair, so at most one member
of each layer is barred beyond the endpoint, and a barred member exists
only when that layer's recursive family --- \(F^{¬e}\) below, the link
\(L\) above --- is itself a Y-family with the endpoint in its universal
pair. Theorem 7.11 says the supply of crossings survives all four
exclusions.

Crossings come in two kinds. Downward: for \(W ∈ F^e\) and \(ξ < e\)
with \(ξ ∉ W\), down-closure puts \(W − e + ξ\) in \(F^{¬e}\). Upward:
for \(Z ∈ F^{¬e}\) and \(ζ ∈ Z\) with \(ζ > e\), down-closure puts
\(Z − ζ + e\) in \(F^e\). Every crossing is of exactly one kind, and
distinct \(ξ\) give distinct partners. Write

\[
d(W) = |[e−1] ∖ W| \quad \text{ for } W ∈ F^e, \quad u(Z) = |Z ∩ (e, n]| \quad \text{ for } Z ∈ F^{¬e},
\]

Call \(W − e + ξ\) a \textbf{downward target} of \(W\) and \(Z − ζ + e\)
an \textbf{upward source} of \(Z\), for the \(ξ\) and \(ζ\) of the two
kinds above; then \(d(W)\) counts the downward targets of \(W\) and
\(u(Z)\) the upward sources of \(Z\).

\textbf{Failure, in the form the proof uses.} Suppose no crossing as in
Theorem 7.11 exists. What that denies depends on the layer sizes, and
the dependence is easy to lose: Theorem 7.11 asks a crossing to
\textbf{meet} a singleton layer's own endpoint and to \textbf{avoid} a
larger layer's barred set. \textbf{When both layers have at least two
members} the first alternative is unavailable on either side, and
failure says that every crossing partner of a \(W ∈ F^e ∖ S^e\) lies in
\(S^{¬e}\), and every crossing partner of a \(Z ∈ F^{¬e} ∖ S^{¬e}\) lies
in \(S^e\); call this statement \((†)\). \((†)\) is the form the main
case uses, and it is \textbf{not} the right form when a layer is a
singleton --- there \(F^e ∖ S^e\) is empty and \((†)\) holds vacuously
while a crossing may well exist. The two degenerate cases are therefore
settled first and directly, below, before \((†)\) is ever invoked. In
particular, by \((†)\), the downward targets of such a \(W\) lie in
\(S^{¬e}\) and the upward sources of such a \(Z\) lie in \(S^e\), so

\[
\begin{gathered}
d(W) ≤ |S^{¬e}| ≤ 2 \quad \text{ for } W ∈ F^e ∖ S^e, \\
u(Z) ≤ |S^e| ≤ 2 \quad \text{ for } Z ∈ F^{¬e} ∖ S^{¬e}. \quad (∗∗)
\end{gathered}
\]

The two are not interchangeable and the proof below is careful about
which it uses. \((†)\) says \textbf{where} a partner lies; \((∗∗)\) says
only \textbf{how many} there can be, and is the weaker of the two, being
a consequence of \((†)\). A step that has to place a particular set
cannot be run from a bound on cardinalities.

\textbf{All of it on the family the figure already carries.} Nothing
above needs a second example: take
\(F = \{123, 124, 125, 126, 134, 135\}\) at \(e = 3\), with \(A = 123\)
and \(B = 124\), exactly as drawn. Then

\[
F^e = \{123, 134, 135\}, \quad F^{¬e} = \{124, 125, 126\}, \quad L = F^e − e = \{12, 14, 15\}
\]

on the punctured ground \(\{1,2,4,5,6\}\). Each layer has three members,
and a Y-family needs four, so neither layer is one and \textbf{both
barred sets are singletons}: \(S^e = \{A\}\) and \(S^{¬e} = \{B\}\). The
heavy crossing of the figure is \(135\)--\(125\); the two sets differ in
exactly \(\{2,3\}\), so they are adjacent, and neither is the prescribed
endpoint of its own layer, so the crossing is usable. The two counts
read off directly --- \(d(134) = d(135) = 1\) and \(d(123) = 0\), since
\(d\) measures how much of \([e−1] = \{1,2\}\) a containing member is
missing, while \(u(124) = u(125) = u(126) = 1\), each avoiding member
carrying exactly one element above \(e\). A reader who holds this
instance can check every statement of the rest of this subsection, and
of the crossing proof that follows it, against it.

\textbf{The layer bottoms.} The case analysis needs the two Gale-least
members of each layer, computed by Lemma 7.2 on the punctured ground
\([n] ∖ \{e\}\). Lemma 7.2 asks for at least two members, so a
second-bottom row is defined only for a layer that has them; a singleton
layer is handled directly in the two degenerate cases and never consults
one. Write \(m_1(G)\) and \(m_2(G)\) for the two Gale-least members of a
family \(G\) on its own ground, so that the \(m_1\) and \(m_2\) of §7.1
are \(m_1(F)\) and \(m_2(F)\). Every entry of the two tables below is
Lemma 7.2 read on the punctured ground, so the only thing needed is
where that ground's prefixes sit inside \([n]\).

\textbf{The puncture rule.} \emph{Write \(⟨j⟩\) for the set of the first
\(j\) elements of \([n] ∖ \{e\}\). Then}

\[
⟨j⟩ = [j] \quad \text{ for } j < e, \quad ⟨j⟩ = [j+1] ∖ \{e\} \quad \text{ for } j ≥ e,
\]

\emph{and the \(j\)-th element of \([n] ∖ \{e\}\) is \(j\) for \(j < e\)
and \(j+1\) for \(j ≥ e\).}

\emph{Proof.} The increasing enumerations of \([n]\) and of
\([n] ∖ \{e\}\) agree in their first \(e−1\) places, and from there on
the second is the first advanced by one. So the prefix is unchanged
while \(j\) stays below \(e\), and gains one element beyond \(e\) once
\(j\) reaches it. \(\blacksquare\)

We state it rather than use it silently because the two arms of a
Y-family straddle the puncture: an index below \(e\) is not shifted and
an index above it is, within one row. The worked regime after the tables
shows what that costs. Since \(W ↦ W − e\) is a bijection \(F^e → L\),
every member of \(L\) has a unique preimage \(X ∪ \{e\}\), and we call
it the \textbf{lift} of \(X\) into \(F^e\). Writing \(μ_1\) and \(μ_2\)
for the lifts of the two Gale-least members of \(L\):

The two Gale-least members of the avoiding layer, on the punctured
ground:

\begin{longtable}[]{@{}lll@{}}
\toprule\noalign{}
regime & \(m_1(F^{¬e})\) & \(m_2(F^{¬e})\) \\
\midrule\noalign{}
\endhead
\bottomrule\noalign{}
\endlastfoot
\(e ≤ k−1\) & \([k+1] ∖ \{e\}\) & \(([k] ∖ \{e\}) ∪ \{k+2\}\) \\
\(e = k\) & \([k+1] ∖ \{k\}\) & \([k−1] ∪ \{k+2\}\) \\
\(e = k+1\) & \([k]\) & \([k−1] ∪ \{k+2\}\) \\
\(e ≥ k+2\) & \([k]\) & \(m_2\) \\
\end{longtable}

And the two Gale-least members of the link, lifted back into \(F^e\):

\begin{longtable}[]{@{}lll@{}}
\toprule\noalign{}
regime & \(μ_1\) & \(μ_2\) \\
\midrule\noalign{}
\endhead
\bottomrule\noalign{}
\endlastfoot
\(e ≤ k−1\) & \([k]\) & \(m_2\) \\
\(e = k\) & \([k]\) & \([k+1] ∖ \{k−1\}\) \\
\(e = k+1\) & \(m_2\) & \([k+1] ∖ \{k−1\}\) \\
\(e ≥ k+2\) & \([k−1] ∪ \{e\}\) & \(([k] ∖ \{k−1\}) ∪ \{e\}\) \\
\end{longtable}

When a layer is a Y-family, Theorem 7.6 on its own ground gives its two
arms, both nonempty, since that is part of being a Y-family. \textbf{The
analysis below never needs a whole arm --- only its first member, the
one at the extreme index} --- so those are what we record, lifted into
\(F^e\) where relevant. \(s\) is not needed at all, and neither is any
outer parameter.

\begin{longtable}[]{@{}
  >{\raggedright\arraybackslash}p{(\columnwidth - 4\tabcolsep) * \real{0.3333}}
  >{\raggedright\arraybackslash}p{(\columnwidth - 4\tabcolsep) * \real{0.3333}}
  >{\raggedright\arraybackslash}p{(\columnwidth - 4\tabcolsep) * \real{0.3333}}@{}}
\toprule\noalign{}
\begin{minipage}[b]{\linewidth}\raggedright
the Y-family
\end{minipage} & \begin{minipage}[b]{\linewidth}\raggedright
first inner member
\end{minipage} & \begin{minipage}[b]{\linewidth}\raggedright
first outer member
\end{minipage} \\
\midrule\noalign{}
\endhead
\bottomrule\noalign{}
\endlastfoot
\(F^{¬e}\), at \(e = k\) & \([k−2] ∪ \{k+1, k+2\}\) & not needed \\
the link \(L\), at \(e = k−1\), lifted & \([k+1] ∖ \{k−2\}\) & not
needed \\
the link \(L\), at \(e = k\), lifted & \([k+1] ∖ \{k−2\}\) &
\([k−2] ∪ \{k, k+2\}\) \\
\end{longtable}

\textbf{One regime worked.} Take \(e = k\). The puncture rule gives
\(⟨k−2⟩ = [k−2]\), \(⟨k−1⟩ = [k−1]\), \(⟨k⟩ = [k+1] ∖ \{k\}\) and
\(⟨k+1⟩ = [k+2] ∖ \{k\}\), and makes the \(j\)-th element of the
punctured ground equal to \(j\) for \(j < k\) and to \(j+1\) for
\(j ≥ k\). One relevant entry from each of the three tables follows.

\emph{\(m_1(F^{¬e})\).} The avoiding layer has rank \(k\) on the
punctured ground, so by Lemma 7.2 its Gale-least member is that ground's
first \(k\) elements, namely \(⟨k⟩ = [k+1] ∖ \{k\}\).

\emph{\(μ_2\).} The link has rank \(k−1\), so
\(m_2(L) = ⟨k−2⟩ ∪ \{the k-th element\} = [k−2] ∪ \{k+1\}\). Lifting
restores \(e = k\), and \([k−2] ∪ \{k, k+1\} = [k+1] ∖ \{k−1\}\).

\emph{The first inner member of \(F^{¬e}\).} If the avoiding layer is a
Y-family, its inner arm is \(⟨k+1⟩\) with one element deleted, and the
first such member deletes the \((k−1)\)-th element of the punctured
ground. That index lies below the puncture and is therefore unshifted,
and \(⟨k+1⟩ = [k+2] ∖ \{k\}\), so the member is
\(([k+2] ∖ \{k\}) ∖ \{k−1\} = [k−2] ∪ \{k+1, k+2\}\).

\hypertarget{proof-of-the-crossing-lemma-and-of-theorem-7.7}{%
\subsubsection{7.8.2 Proof of the crossing lemma, and of Theorem
7.7}\label{proof-of-the-crossing-lemma-and-of-theorem-7.7}}

\textbf{The shape of the argument, before the details.} Suppose no
usable crossing exists, and read off what that forces. The proof splits
on the sizes of the two layers. There are four cardinality patterns,
which group into three cases once duality is used:

\begin{itemize}
\tightlist
\item
  both layers are singletons --- settled directly, in three lines;
\item
  exactly one layer is a singleton --- settled directly, and its mirror
  image follows by duality;
\item
  both layers have at least two members --- the main case, and the only
  long one.
\end{itemize}

In the main case a counting argument, the \textbf{range reduction},
shows that failure is possible only for four values of the splitting
element, \(k−1 ≤ e ≤ k+2\). Two of those four are treated here and the
other two follow by duality, so four regimes cost the work of two. Every
case ends in a contradiction except one configuration, up to duality,
and in that one the family is forced to be a Y-family taken at its
Gale-least pair --- the exception, which the hypothesis excludes.

Three statements from Section 7.8.1 are used throughout and are worth
having in view: \((∗)\) says each layer of size at least two has an
unbarred member; \((†)\) says that under failure every crossing partner
of an unbarred member lands in the \emph{other} layer's barred set; and
\((∗∗)\) is the counting consequence of \((†)\), that \(d(W) ≤ 2\) and
\(u(Z) ≤ 2\) for unbarred \(W\) and \(Z\).

\textbf{Proof of Theorem 7.11.}

Suppose failure, and write \(f^e = |F^e|\), \(f^{¬e} = |F^{¬e}|\), both
at least one since \(A ∈ F^e\) and \(B ∈ F^{¬e}\). By Lemma 7.10 it
suffices to treat the parameters up to duality. Under \(Z ↦ Z^*\) the
two layers exchange, so \((f^e, f^{¬e})\) becomes \((f^{¬e}, f^e)\), the
endpoints exchange, and

\[
\begin{gathered}
e^* − (n−k) = (n+1−e) − (n−k) = k + 1 − e, \\
\text{ so } \quad e − k = t \quad \text{ becomes } \quad e^* − (n−k) = 1 − t.
\end{gathered}
\]

Hence the degenerate case \(f^e = 1 < f^{¬e}\) is dual to
\(f^{¬e} = 1 < f^e\); and in the main case the regime \(e = k−1\) is
dual to \(e = k+2\), and \(e = k\) to \(e = k+1\).

Failure itself transfers, but the statement that transfers is the full
one and not \((†)\). Call a crossing \((W, Z)\) \textbf{usable} when
\(W = A\) if \(f^e = 1\) and \(W ∉ S^e\) otherwise, and \(Z = B\) if
\(f^{¬e} = 1\) and \(Z ∉ S^{¬e}\) otherwise; failure is the absence of a
usable crossing, which is the negation of Theorem 7.11's conclusion with
the degenerate layers included. Crossings being the edges of \(J(F)\)
that join the two layers, Lemma 7.10 carries crossings of \((F, e)\)
bijectively onto crossings of \((F^*, e^*)\) with the two ends
exchanged; the barred sets correspond as displayed above, and \(f^e\)
and \(f^{¬e}\) exchange with the layers, so each of the two cardinality
clauses is carried to the other. Usability is therefore preserved, and a
crossing witnessing non-failure for \((F, A, B, e)\) maps to one
witnessing non-failure for \((F^*, B^*, A^*, e^*)\), and conversely.

Two things are deliberate here. We argue from usability and not from
\((∗∗)\), which is a consequence of failure rather than a restatement of
it, and transferring a consequence would not transfer the hypothesis.
And we do not argue from \((†)\) either: duality is used below to obtain
the degenerate case \(f^{¬e} = 1 < f^e\) from \(f^e = 1 < f^{¬e}\), and
a singleton layer is precisely where \((†)\) is not the right form. The
hypothesis that \((F, \{A,B\})\) is not the exception transfers because
Y-families and universal pairs do. \textbf{We therefore prove the cases
\(f^e = f^{¬e} = 1\), \(f^e = 1 < f^{¬e}\), \(e = k−1\) and \(e = k\),
and obtain the remaining three by duality.}

\textbf{Case 1. Both layers are singletons.} If \(f^e = f^{¬e} = 1\)
then \(F = \{A, B\}\), which is \(\{m_1, m_2\}\) by Lemma 7.2(iii), and
\(A △ B = \{k, k+1\} ∋ e\). If \(e = k+1\) then \(A = m_2\) and
\(A − (k+1) + k = m_1 = B\); if \(e = k\) then \(A = m_1\) and
\(A − k + (k+1) = m_2 = B\). Either way the crossing exists, and both
layers being singletons, no exclusion applies. There is no failure.

\textbf{Case 2. Exactly one layer is a singleton.} Let \(f^e = 1\) and
\(f^{¬e} ≥ 2\), so \(F^e = \{A\}\) and \(W = A\) is forced; failure
means that every crossing partner of \(A\) lies in \(S^{¬e}\).

First, \(e ≥ k\). If \(e ≤ k−1\) then both \(m_1 = [k]\) and
\(m_2 = [k−1] ∪ \{k+1\}\) contain \(e\), and they are distinct members
of \(F\) by Lemma 7.2, since \(|F| ≥ 2\); so \(F^e\) would have two
members. Hence \(e ≥ k\).

Next, \(A ∖ \{e\} = [k−1]\). Since \(F^e = \{A\}\), the link
\(L = \{A − e\}\) is a one-member shifted family on the punctured
ground, so \(A − e\) is its Gale-least member and \(A = μ_1\). Reading
the table at \(e = k\), at \(e = k+1\) and at \(e ≥ k+2\) gives
\(μ_1 ∖ \{e\} = [k−1]\) in each of the three columns that \(e ≥ k\)
leaves available.

Let \(Π = \{ [k−1] ∪ \{ξ\} : ξ ∉ [k−1] ∪ \{e\} \} ∩ F\) be the family of
one-point extensions of \([k−1]\) that lie in \(F\) and avoid \(e\); the
excluded extension \([k−1] ∪ \{e\}\) is \(A\) itself. Every member of
\(Π\) is a \(k\)-set avoiding \(e\), so lies in \(F^{¬e}\), and equals
\(A − e + ξ\), so \(Π\) consists of crossing partners of \(A\). Reading
the table at \(e ≥ k\) shows \(m_1(F^{¬e})\) and \(m_2(F^{¬e})\) are
both of this form --- at \(e = k\) they are \([k−1] ∪ \{k+1\}\) and
\([k−1] ∪ \{k+2\}\), at \(e = k+1\) they are \([k]\) and
\([k−1] ∪ \{k+2\}\), and at \(e ≥ k+2\) they are \([k]\) and \(m_2\) ---
and both lie in \(F^{¬e}\), by Lemma 7.2 since \(f^{¬e} ≥ 2\). So
\(\{m_1(F^{¬e}), m_2(F^{¬e})\} ⊆ Π\).

If \(S^{¬e} = \{B\}\) then \(Π ∖ \{B\}\) is nonempty, since \(|Π| ≥ 2\),
and any of its members is a valid \(Z\). If \(S^{¬e}\) has two members
then \(F^{¬e}\) is a Y-family, and by Theorem 7.6 on the ground
\([n] ∖ \{e\}\) its outer arm is nonempty and consists of sets of the
form \([k−1] ∪ \{a\}\), again members of \(Π\). An arm of a Y-family is
disjoint from its universal pair, so the arm supplies a member of
\(Π ∖ S^{¬e}\), which is a valid \(Z\). Either way failure is
impossible.

The case \(f^{¬e} = 1 < f^e\) follows by duality.

\textbf{Case 3. Both layers have at least two members --- the range
reduction.} Now let \(f^e, f^{¬e} ≥ 2\), so by \((∗)\) both
\(F^e ∖ S^e\) and \(F^{¬e} ∖ S^{¬e}\) are nonempty. A member
\(W ∈ F^e ∖ S^e\) contains \(e\) and so at most \(k−1\) elements of
\([e−1]\), giving \(d(W) ≥ (e−1) − (k−1) = e − k\). A member
\(Z ∈ F^{¬e} ∖ S^{¬e}\) has at most \(e−1\) elements below \(e\), giving
\(u(Z) ≥ k − (e−1)\). With \((∗∗)\),

\[
e − k ≤ 2 \quad \text{ and } \quad k − e + 1 ≤ 2, \quad \text{ so } \quad k − 1 ≤ e ≤ k + 2.
\]

\textbf{Case 3a. The regime \(e = k−1\).} Every \(Z ∈ F^{¬e} ∖ S^{¬e}\)
has \(u(Z) ≥ 2\), so \(|S^e| = 2\) by \((∗∗)\): the link \(L\) is a
Y-family and \(S^e = \{μ_1, μ_2\} = \{[k], m_2\}\) by the table. Since
\(L\) has rank \(k−1\) and a Y-family has rank at least two, this regime
is empty unless \(k ≥ 3\), which we may therefore assume. \textbf{Fix
one such \(Z\) and keep it for the rest of the case}; by \((∗)\) there
is one. Also \(u(Z) = 2\) exactly; by \((†)\) both upward sources lie in
\(S^e\), and \(S^e\) has two members, so they are exactly \([k]\) and
\(m_2\).

Write \(Z ∩ (k−1, n] = \{ζ < ζ'\}\). Removing \(ζ\) cannot produce
\([k]\): the retained element \(ζ'\) would have to lie in \([k]\),
forcing \(ζ' = k\) and hence \(k−1 < ζ < k\), which is impossible. So
removing \(ζ'\) gives \([k]\) and removing \(ζ\) gives \(m_2\), that is

\[
Z = [k] − (k−1) + ζ' = [k−2] ∪ \{k, ζ'\} \quad \text{ and } \quad Z = m_2 − (k−1) + ζ = [k−2] ∪ \{k+1, ζ\}.
\]

Hence \(\{k, ζ'\} = \{k+1, ζ\}\), forcing \(ζ = k\) and \(ζ' = k+1\), so
the \(Z\) we fixed is \([k+1] ∖ \{k−1\}\). \textbf{It was chosen
unbarred, so \([k+1] ∖ \{k−1\} ∉ S^{¬e}\), and that is the whole of what
this case needs from the avoiding layer} --- not its size, not its
second bottom, and not \(B\).

Since \(L\) is a Y-family its lifted inner arm is nonempty, and its
\textbf{first} member is \(W_0 = [k+1] ∖ \{k−2\}\). An arm of a Y-family
is disjoint from its universal pair, and \(S^e\) is the lift of that
pair, so \(W_0 ∈ F^e ∖ S^e\) and \((†)\) applies to it. Since
\(k−2 < k−1 = e\) and \(k−2 ∉ W_0\), the set
\(W_0 − (k−1) + (k−2) = [k+1] ∖ \{k−1\}\) is a downward target of
\(W_0\), and it is the unbarred set just fixed. So \(W_0\) has a
downward target outside \(S^{¬e}\), contradicting \((†)\). So
\(e = k−1\) is impossible.

\textbf{Case 3b. The regime \(e = k\).} Here \(d(W) ≥ 0\) and
\(u(Z) ≥ 1\) carry no information by themselves, and the argument splits
on the two barred sets. There are four configurations, which we group
into three.

\textbf{(a) \(S^e = \{A\}\) and \(S^{¬e} = \{B\}\).} Each
\(Z ∈ F^{¬e} ∖ \{B\}\) has \(u(Z) = 1\) by \((∗∗)\), and by \((†)\) that
one upward source is \(A\), so \(Z = A − k + ζ\) for some \(ζ > k\)
outside \(A\). Then \(u(Z) = |A ∩ (k, n]| + 1\), so \(A ∩ (k, n] = ∅\),
and since \(e = k\) lies in \(A\) this gives \(A = [k] = m_1\) and
\(F^{¬e} ∖ \{B\} ⊆ \{ [k−1] ∪ \{ζ\} : ζ > k \}\).

On the other side each \(W ∈ F^e ∖ \{A\}\) has \(d(W) ≤ 1\). If
\(d(W) = 0\) then \([k−1] ⊆ W\) and \(k ∈ W\), so \(W = [k] = A\), which
is excluded; so \(d(W) = 1\) and \(W = ([k−1] ∖ \{ξ_W\}) ∪ \{k, γ_W\}\)
with \(ξ_W ∈ [k−1]\) and \(γ_W ≥ k+1\). Its downward target
\(W − k + ξ_W = [k−1] ∪ \{γ_W\}\) must lie in \(S^{¬e} = \{B\}\) by
\((†)\), so \(B = [k−1] ∪ \{γ\}\) with \(γ = γ_W\) the same for every
such \(W\).

Suppose \(γ ≥ k+2\). Pick any \(W ∈ F^e ∖ \{A\}\), which exists since
\(f^e ≥ 2\), and decrement \(γ\) to \(k+1\), which is legal since
\(k+1 < γ\) and \(k+1 ∉ W\). The result
\(W' = ([k−1] ∖ \{ξ_W\}) ∪ \{k, k+1\}\) lies in \(F\), contains
\(e = k\), and is not \([k] = A\), so \(W' ∈ F^e ∖ \{A\}\) with
\(γ_{W'} = k+1 ≠ γ\), a contradiction. Hence \(γ = k+1\) and
\(B = m_2\).

So \(F^e ∖ \{A\} ⊆ \{ [k+1] ∖ \{ξ\} : ξ ∈ [k−1] \}\) and
\(F^{¬e} ⊆ \{ [k−1] ∪ \{ζ\} : ζ ≥ k+1 \}\), and down-closure makes both
index ranges contiguous: \(F^{¬e} = \{ [k−1] ∪ \{ζ\} : k+1 ≤ ζ ≤ ρ \}\)
with \(ρ ≥ k+2\) since \(f^{¬e} ≥ 2\), and
\(F^e = \{m_1\} ∪ \{ [k+1] ∖ \{ξ\} : s ≤ ξ ≤ k−1 \}\) with \(s ≤ k−1\)
since \(f^e ≥ 2\). Collecting,

\[
F = \{m_1, m_2\} ∪ \{ [k−1] ∪ \{ζ\} : k+2 ≤ ζ ≤ ρ \} ∪ \{ [k+1] ∖ \{ξ\} : s ≤ ξ ≤ k−1 \} = F(k; ρ, s),
\]

with both arms nonempty. By Theorem 7.6, \(F\) is a Y-family with
universal pair \(\{m_1, m_2\} = \{A, B\}\). That is the exception, which
the hypothesis excludes. \textbf{This is the one configuration in which
failure is consistent, and it forces the exception rather than
contradicting outright.}

\textbf{(b) \(S^e = \{A\}\) and
\(S^{¬e} = \{m_1(F^{¬e}), m_2(F^{¬e})\}\).} Then \(F^{¬e}\) is a
Y-family, so its inner arm is nonempty; take its first member
\(Z = [k−2] ∪ \{k+1, k+2\}\) from the table. An arm of a Y-family is
disjoint from its universal pair, so \(Z ∈ F^{¬e} ∖ S^{¬e}\) and
\((∗∗)\) applies to it. It contains both \(k+1\) and \(k+2\) and avoids
\(k\), so \(u(Z) = 2 > 1 = |S^e|\), a contradiction.

\textbf{(c) \(S^e = \{μ_1, μ_2\} = \{[k], [k+1] ∖ \{k−1\}\}\), either
\(S^{¬e}\).} Then \(L\) is a Y-family, so both its lifted arms are
nonempty, and \textbf{one member from each is all this case needs}. Take
the first of each, from the table: \(W_1 = [k+1] ∖ \{k−2\}\) from the
inner arm and \(W_2 = [k−2] ∪ \{k, k+2\}\) from the outer. Both lie in
\(F^e ∖ S^e\), since an arm is disjoint from the universal pair whose
lift is \(S^e\). Now \(W_1\) misses exactly \(k−2\) from \([k−1]\), so
\(d(W_1) = 1\) with downward target
\(W_1 − k + (k−2) = [k+1] ∖ \{k\} = m_2\); and \(W_2\) misses exactly
\(k−1\) from \([k−1]\), so \(d(W_2) = 1\) with downward target
\([k−1] ∪ \{k+2\}\). By \((†)\) both targets lie in \(S^{¬e}\), and they
are distinct, so

\[
S^{¬e} = \{ m_2, [k−1] ∪ \{k+2\} \},
\]

which is two members: \(F^{¬e}\) is therefore a Y-family, and by the
table these are its two Gale-least members. Its inner arm is then
nonempty with first member \(Z_3 = [k−2] ∪ \{k+1, k+2\}\), which is
unbarred, being in an arm. But replacing \(k+1\) by \(k\) in \(Z_3\)
gives exactly \(W_2\), so \(Z_3\) and \(W_2\) are joined by a crossing
with both ends unbarred --- a usable crossing, contradicting failure.

\textbf{The barred-set transfer.} \emph{Write \(e^* = n+1−e\). Then}

\[
S^e(F, A)^* = S^{¬e^*}(F^*, A^*) \quad \text{ and } \quad S^{¬e}(F, B)^* = S^{e^*}(F^*, B^*).
\]

\emph{Proof.} The care is needed because \textbf{being a Y-family is not
an isomorphism invariant}: Section 7.7 exhibits \(\{12,13,14,23\}\) and
\(\{14,23,24,34\}\), isomorphic as Johnson graphs with only the first a
Y-family. So the transfer cannot be run on the definition of a barred
set. It is run on the graph condition, which \emph{is} invariant, with
Theorem 7.6 converting at each end.

\begin{enumerate}
\def\labelenumi{(\arabic{enumi})}
\item
  \emph{Read the bar as a condition on a layer's Johnson graph.} By
  Section 7.8.1, \(S^e(F, A)\) has two members exactly when
  \(J(F^e) ≅ K_2 ∨ (K_p ⊔ K_q)\) with \(A\) in the dominating pair, and
  \(S^e(F, A)\) is then that pair. That is Theorem 7.6 applied to the
  link and transported along the isomorphism \(J(F^e) → J(L)\); what it
  leaves behind mentions \(F^e\) and not \(L\), so the rank \(k−1\) has
  gone out of it. The definition of \(S^{¬e}(F, B)\) is already of this
  form on \(F^{¬e}\).
\item
  \emph{Transport.} Lemma 7.10 makes \(X ↦ X^*\) an isomorphism
  \(J(F) → J(F^*)\) carrying \(F^e\) onto \((F^*)^{¬e^*}\) and
  \(F^{¬e}\) onto \((F^*)^{e^*}\). Carrying the clique-sum shape with a
  named member in the dominating pair is what an isomorphism does, so
  \(J((F^*)^{¬e^*})\) has that shape at \(A^*\) exactly when \(J(F^e)\)
  has it at \(A\), and the dominating pairs correspond.
\item
  \emph{Read back.} Apply (1) on \(F^*\). For its avoiding layer the
  criterion is already stated on the layer and nothing more is needed;
  for its containing layer it is transported through \((F^*)\)'s own
  link, and that is the second use of Theorem 7.6. The two ranks that
  looked mismatched --- \(k−1\) through a link on one side, \(n−k\) on a
  layer on the other --- never meet, because step (1) removed both.
\item
  \emph{Conclude.} The two-membered cases correspond by (2) and (3), and
  the one-membered cases correspond because \(A^*\) and \(B^*\) are the
  images of \(A\) and \(B\). \(\blacksquare\)
\end{enumerate}

\textbf{Case 3c. The regimes \(e = k+1\) and \(e = k+2\), by duality.}
Write \(k^* = n−k\) for the rank of \(F^*\). At \(e = k+1\) we get
\(e^* = n−k = k^*\), and at \(e = k+2\) we get \(e^* = k^*−1\), so these
two regimes are carried onto the regimes \(e = k\) and \(e = k−1\) of
\(F^*\), which are Cases 3b and 3a. By Lemma 7.10, \(F^*\) is shifted of
rank \(k^*\) with \(|F^*| = |F|\), and \(X ↦ X^*\) is an isomorphism
\(J(F) → J(F^*)\) exchanging the two layers; by the transfer above it
carries \(S^e\) and \(S^{¬e}\) onto \(S^{¬e^*}\) and \(S^{e^*}\). One
relabeling is needed and is worth naming: \(e ∈ A\) gives \(e^* ∉ A^*\)
and \(e^* ∈ B^*\) by Lemma 7.10(iii), so on the dual side the member
containing the splitting element is \(B^*\), and it is \(B^*\) that
plays \(A\)'s role. That costs nothing, a Hamilton \(A\)--\(B\) path
being a Hamilton \(B\)--\(A\) path. A usable crossing for
\((F^*, B^*, A^*)\) at \(e^*\) is therefore exactly the image of one for
\((F, A, B)\) at \(e\), and \((F, \{A,B\})\) is the exception precisely
when \((F^*, \{A^*,B^*\})\) is, by the last clause of Lemma 7.10. Cases
3a and 3b applied to \(F^*\) therefore settle both regimes, the
exception arising at \(e = k+1\) exactly as configuration (a) arises at
\(e = k\).

\textbf{Assembly.} The four cases \(f^e = f^{¬e} = 1\),
\(f^e = 1 < f^{¬e}\), \(f^{¬e} = 1 < f^e\) and \(f^e, f^{¬e} ≥ 2\)
exhaust the possibilities, and in the last the range reduction leaves
only \(k−1 ≤ e ≤ k+2\). Outside configuration (a) of \(e = k\) and its
dual in \(e = k+1\), failure is impossible; in those two it forces
\((F, \{A,B\})\) to be the exception. This proves Theorem 7.11.
\(\blacksquare\)

What remains is to spend it. Theorem 7.7's obstruction half is already
done --- Lemma 7.4 with Theorem 7.6 --- so the following completes its
existence half, and with it the theorem.

The converse is worth recording, since it shows the hypothesis is
exactly right. In the exception \(A △ B = m_1 △ m_2 = \{k, k+1\}\), and
both \(e = k\) and \(e = k+1\) are genuine failures --- one for each
orientation of the pair, since Theorem 7.11 requires \(e ∈ A\). Taking
\(A = m_1\) the qualifying element is \(e = k\), and taking \(A = m_2\)
it is \(e = k+1\); the Y-family satisfies every constraint of
configuration (a) in the first case and of its dual in the second. They
must be failures, since a valid crossing would assemble a Hamilton
\(m_1\)--\(m_2\) path, which Lemma 7.4 forbids.

Theorem 7.11 gives Theorem 7.7 at once. For a non-exceptional pair,
choose \(e ∈ A ∖ B\), which is nonempty because \(A ≠ B\) and the two
have equal size; interchanging \(A\) and \(B\) costs nothing, a Hamilton
\(A\)--\(B\) path being a Hamilton \(B\)--\(A\) path. Take the crossing
Theorem 7.11 supplies. A layer with a single member needs no induction:
the path required of it is the one-vertex path. A layer with two or more
gets its far endpoint outside its barred set, which contains its near
endpoint, so the two are distinct and form a pair that is not that
layer's exception; the induction hypothesis applies there. Splicing the
two paths along the crossing gives the Hamilton \(A\)--\(B\) path. This
completes the existence half of Theorem 7.7, and with the obstruction
half --- Lemma 7.4 with Theorem 7.6 --- the theorem. \(\blacksquare\)

\hypertarget{where-this-sits-in-the-literature}{%
\subsection{7.9 Where this sits in the
literature}\label{where-this-sits-in-the-literature}}

Corollary 7.9 explains why the exception does not arise in the principal
case, which is where the literature we know of sits. Two words need
separating. A \textbf{shifted family} is any nonempty down-set of the
Gale order, the object of this section; a \textbf{shifted matroid} is
one that happens to be a matroid's basis family. Not every shifted
family is one, and that is the whole reason the exception can occur at
all. A shifted family is a matroid basis family exactly when it is
principal, a theorem of Klivans {[}19, Thm. 5.4.1{]}, whose proof of
Theorem 5.5.2 records the basis form we use: the bases of the shifted
matroid with top element \(M\) are exactly the \(k\)-sets Gale-below
\(M\).

\textbf{That class has been found repeatedly and renamed each time},
which locates what is new here and matters for searching, since a reader
who knows one name will miss most of the literature. Crapo {[}11{]}
introduced the matroids in order to show that an \(n\)-element set
carries at least \(2^n\) nonisomorphic matroids, one for each
length-\(n\) binary word; Welsh {[}35{]} proved that bound within the
transversal matroids, by way of the nested presentations that supply the
first name below. Bonin and de Mier {[}7{]} set out the history and we
follow their account.

\begin{longtable}[]{@{}
  >{\raggedright\arraybackslash}p{(\columnwidth - 2\tabcolsep) * \real{0.5000}}
  >{\raggedright\arraybackslash}p{(\columnwidth - 2\tabcolsep) * \real{0.5000}}@{}}
\toprule\noalign{}
\begin{minipage}[b]{\linewidth}\raggedright
name
\end{minipage} & \begin{minipage}[b]{\linewidth}\raggedright
introduced as such by
\end{minipage} \\
\midrule\noalign{}
\endhead
\bottomrule\noalign{}
\endlastfoot
nested matroids & Welsh {[}35{]} \\
Schubert matroids & Sohoni {[}33{]} \\
shifted matroids & Ardila {[}2{]}, as the matroids whose independence
complex is shifted \\
freedom matroids & Crapo and Schmitt {[}12{]} \\
generalized Catalan matroids & Bonin and de Mier {[}7{]} \\
PI-matroids & Billera, Jia and Reiner {[}6{]} \\
non-exchangeable matroids & Partida {[}27{]}, who adds that no such list
is exhaustive \\
\end{longtable}

We say \textbf{shifted matroids} throughout, and these names should be
read as denoting one class. They matter for more than credit, because
the class carries a Hamiltonicity literature under them: Fernandes,
Hernández-Vélez, de Pina and Ramírez Alfonsín {[}14{]} bound the number
of Hamilton cycles through any edge of the basis graph of a generalized
Catalan matroid --- a lower bound on cycles, not a statement about paths
between prescribed bases, and so not a route to anything below.

\textbf{On principal families Theorem 7.7 recovers the known conclusion;
its new content is the non-principal case.} The exception cannot occur
in a principal family, by Corollary 7.9, so a class that has been
rediscovered repeatedly is one on which the exception is empty. It is
not the only such class: \(F = \{123, 124, 125, 134, 135, 234\}\) is
shifted, is not principal --- \(135\) and \(234\) are incomparable
maxima --- and is not a Y-family either, so Corollary 7.8 applies and
\(J(F)\) is Hamilton-connected with no exceptional pair at all. What
Theorem 7.7 adds is the non-principal case, where a shifted family need
not be a matroid basis family, and we are not aware of prior work on
Hamiltonicity there.

Two proofs of the principal case are available. Naddef and Pulleyblank
{[}25{]} proved that the graph of a 0/1-polytope is Hamilton-connected
unless it is a hypercube, and for a matroid basis family \(J(F)\) is the
1-skeleton of \(\operatorname{conv}(F)\) --- the convex hull of the
indicator vectors of the members --- since the edges of a matroid
polytope are exactly the pairs of bases differing by a single exchange
{[}17{]}; a principal family with three or more members contains a
triangle by Corollary 7.3, so for \(|F| ≥ 3\) the hypercube alternative
never applies, and for \(|F| ≤ 2\) the graph is \(K_1\) or \(K_2\).
Hladík and Fink {[}18{]} give an elementary, self-contained proof.
Neither route is used below: Theorem 7.7 is proved directly and
Corollary 7.9 recovers the principal case from it.

\textbf{Two nearby results look as though they settle the question, and
neither does.}

Alspach and Liu {[}1{]} proved the base graph of a matroid
Hamilton-connected, and indeed edge-pancyclic and panconnected. All
three of their theorems assume the matroid \textbf{simple} --- no loops,
and no two elements parallel, where \(x\) and \(y\) are parallel when
neither is a loop and no basis contains both --- and their closing
example shows the hypothesis cannot be dropped. Shifted matroids fail it
in bulk: at rank one every one with two or more bases fails, since the
bases are singletons and so no basis contains a pair, making every two
nonloop elements parallel; the proportion falls as the rank rises, but
over all ranks on a fixed ground it is about three quarters. Naddef and
Pulleyblank carry no such hypothesis, which is why it is the polytope
route that gives the principal case above.

That route does not extend beyond the matroid case, because \(J(F)\) is
in general a \textbf{proper} subgraph of the skeleton of
\(\operatorname{conv}(F)\). Naddef and Pulleyblank's theorem is
unaffected --- it holds of any 0/1-polytope --- and what fails is the
identification of \(J(F)\) with the skeleton. Every pair at symmetric
difference two is an edge of \(\operatorname{conv}(F)\), since the
functional that is \(1\) on \(X ∩ Y\), \(1/2\) on \(X △ Y\) and \(0\)
elsewhere is maximized exactly on \(\{X, Y\}\); so \(J(F)\) sits inside
the skeleton, and the containment can be strict. \textbf{At every
exception it is strict in the way that matters.} A Y-family has at least
four members, so it contains a triangle by Corollary 7.3, so
\(\operatorname{conv}(F)\) is not a hypercube and its skeleton is
Hamilton-connected by Naddef and Pulleyblank --- while \(J(F)\) is not,
by Lemma 7.4. The smallest exception shows it concretely: for
\(F = \{12, 13, 14, 23\}\) the four indicator vectors are affinely
independent, so \(\operatorname{conv}(F)\) is a 3-simplex with skeleton
\(K_4\), but \(|14 △ 23| = 4\), so \(J(F) = K_4 − e\) and there is no
Hamilton \(12\)--\(13\) path. Nor is the divergence confined to the
exceptions: at \(n = 6\), \(k = 3\) the skeleton is strictly larger for
45 of the 65 shifted families, and 41 of those are not Y-families.
Hamilton-connectivity of the skeleton therefore says nothing about
\(J(F)\) where the exception lives.

Ruskey, Sawada and Williams {[}30{]} prove that listing any
\textbf{bubble language} --- a family of indicator strings closed under
swapping its first \(01\) to \(10\), which every shifted family is,
being closed under every such swap --- in \textbf{cool-lex} order yields
a Gray code in which successive strings differ by one \textbf{or two}
transpositions. The ``or two'' is the whole of the difference. Two
disjoint transpositions move four coordinates, while the edges of
\(J(F)\) are the pairs differing in two. A cool-lex listing is therefore
not in general a Hamilton path of \(J(F)\), and the gap is not an
artifact of small cases: on the full family of \(k\)-subsets, cool-lex
takes a distance-four step once for \((n,k) = (4,2)\), twice for
\((5,2)\), five times for \((6,3)\), nine times for \((7,3)\) and
nineteen times for \((8,4)\).

Finally, what makes a classification the right form of answer. Merino,
Namrata and Williams {[}23{]} show that deciding whether an arbitrary
list of \(k\)-subsets admits a single-exchange Gray code ---
equivalently, whether its Johnson graph has a Hamilton path --- is
NP-hard, so unless \(P = NP\) there is no polynomial-time criterion in
general. On shifted families the obstruction is instead confined to a
single class, recognizable by inspection.

\hypertarget{choosing-a-pivot-that-dodges-the-exception}{%
\section{8. Choosing a pivot that dodges the
exception}\label{choosing-a-pivot-that-dodges-the-exception}}

Section 7 supplies a Hamilton path in the quotient except at one
configuration, and Section 6 supplies a pivot with two properties. The
assembly of Section 10 walks the quotient between the two
\textbf{prescribed} projections \(A = N_{G_0}(v)\) and
\(B = N_{G_1}(v)\), which it does not get to choose. So it needs a pivot
at which the prescribed pair is not the exception of Theorem 7.7. That
is the content of this section.

This section needs two properties of a single pivot \(v\). They are
defined separately because the section treats them differently: one is
supplied to it, and the other is what it has to avoid.

Call \(v\) \textbf{admissible} for the pair \((G_0, G_1)\) when it
separates them and some \(v\)-fiber is non-bipartite. Admissibility is
what Section 6 delivers.

Call \(v\) \textbf{exceptional} for the pair when the quotient family at
\(v\) is a Y-family whose universal pair is
\(\{N_{G_0}(v), N_{G_1}(v)\}\). Exceptionality is what this section must
dodge, and it is the subject of everything below.

\textbf{Exceptionality factors, and the factoring is the whole
strategy.} Read the definition of \emph{exceptional} again. It asks two
things at once: that the quotient family \(I_v\) \textbf{is a Y-family},
and that the two prescribed realizations \textbf{project onto its
universal pair}. The first is a question about \(d\) and \(v\) and
nothing else. Name it:

\medskip\par\noindent

\textbf{Definition.} \(H(d)\) is the set of ground vertices \(v\) whose
quotient family \(I_v\) is a Y-family.

\par\medskip

So \(v\) is exceptional for \((G_0, G_1)\) exactly when \(v ∈ H(d)\)
\textbf{and} the prescribed pair projects onto the universal pair of
\(I_v\) --- the first condition depending on \(d\) alone, the second on
the realizations. The second is a single test, because by Theorem 7.6 a
Y-family's universal pair is always \(\{[k], [k−1] ∪ \{k+1\}\}\) read in
the degree order. So \(d\) alone says \textbf{which} two neighborhoods
would make \(v\) exceptional, before any realization is examined; what
the pair decides is only whether it happens to land on them.

That split is what makes the section possible. The assembly cannot
choose \(G_0\) and \(G_1\), which are prescribed, but it can choose
\(v\) --- and the expensive structural half of exceptionality is
computable from \(d\) in advance, leaving a cheap projection test that
depends on the pair. Everything below therefore works on \(H(d)\) rather
than on the pair.

\textbf{Two pairs appear here and they are different kinds of object.}
The \textbf{universal pair} is a pair of \emph{neighborhoods}, two
members of \(I_v\) and so two vertices of the quotient. The
\textbf{forced pair} \(g(v)\) of Lemma 8.2 is a pair of \emph{ground
vertices}, the two coordinates in which those neighborhoods differ.
Lemma 8.2 is the step that crosses between them, and it is what puts the
argument in reach of \(d\).

\textbf{Theorem 8.1.} \emph{Under the hypotheses of Theorem 6.1, every
pair of distinct realizations has an admissible pivot that is not
exceptional.}

The exception is not a formality. It occurs at admissible pivots on
degree sequences satisfying every hypothesis of Theorem 6.1; the
smallest are \(d = (4,4,3,3,3,1)\) and \(d = (4,2,2,2,1,1)\), each with
three such pivots, and at ground order eight there are thirty admissible
Y-family pivots on sequences satisfying every hypothesis of Theorem 6.1.
Theorem 8.1 is therefore strictly stronger than Theorem 6.1 and does not
follow from it.

\textbf{The shape of the argument, in three moves.} Being exceptional
looks like a property of the pair \((G_0, G_1)\), which would be
hopeless to control, since the assembly does not choose the pair. The
section shows it is not.

\begin{enumerate}
\def\labelenumi{\arabic{enumi}.}
\tightlist
\item
  \textbf{Pin it.} Section 8.1: an exceptional pivot has exactly two
  \(Δ\)-neighbors, and they are not free --- they are the \textbf{forced
  pair} \(g(v)\), computable from \(d\) and \(v\) with the realizations
  never consulted. So the pair can make \(v\) exceptional in one way
  only.
\item
  \textbf{Confine it.} Sections 8.2 and 8.3: the pivots that can be
  exceptional all sit in \textbf{one degree class}, and that class has
  exactly three members. This is the long part, and it is where the
  four-corner analysis lives.
\item
  \textbf{Count it.} Section 8.4: at most two of the three can be
  exceptional for the same pair, because the three would form a triangle
  and a two-coloring cannot alternate around an odd cycle. So one
  survives, which is Theorem 8.1.
\end{enumerate}

\textbf{Pin, confine, count.} A reader who wants only the statement can
take Theorem 8.1 and go on. What later sections do cite from inside the
section is small and specific: Section 10 uses Theorem 8.3 and the
triangle count inside the proof of Theorem 8.1 when it walks the worked
example, and Section 10.4 names Lemmas 8.3d and 8.3e as the two places
the Erdős--Gallai criterion enters. Nothing else here is quoted
elsewhere.

\hypertarget{the-forced-pair-is-determined-by-d-alone}{%
\subsection{\texorpdfstring{8.1 The forced pair is determined by \(d\)
alone}{8.1 The forced pair is determined by d alone}}\label{the-forced-pair-is-determined-by-d-alone}}

Write the ground vertices other than \(v\) as
\(u_1 ≻ u_2 ≻ ⋯ ≻ u_{n−1}\), ordered by non-increasing degree with ties
broken by label, and put \(k = d(v)\). Let \(Δ = E(G_0) △ E(G_1)\), so
that \(\deg_Δ(w) = |N_{G_0}(w) △ N_{G_1}(w)|\).

\textbf{Lemma 8.2.} \emph{If \(v\) is exceptional then}

\[
N_{G_0}(v) ∩ N_{G_1}(v) = \{u_1, …, u_{k−1}\} \quad \text{ and } \quad N_{G_0}(v) △ N_{G_1}(v) = \{u_k, u_{k+1}\}.
\]

\emph{Proof.} The quotient family at \(v\) is a Y-family, and its
universal pair is \(\{A, B\}\). By Theorem 7.6 that pair is
\(\{[k], [k−1] ∪ \{k+1\}\}\) read in the order \(u_1 ≻ ⋯ ≻ u_{n−1}\).
Intersect and difference. \(\blacksquare\)

Lemma 8.2 is the reason this section works. It says that an exceptional
pivot has \(\deg_Δ(v) = 2\), and that its two \(Δ\)-neighbors are
\(\{u_k, u_{k+1}\}\), a pair determined by \(d\) and \(v\) with no
reference to \(G_0\) or \(G_1\). We write \(g(v)\) for that pair and
call it the \textbf{forced pair}. In particular a vertex with
\(\deg_Δ ≥ 4\) is never exceptional.

Two inputs deserve to be named here, because Lemma 8.2 reads a statement
about a quotient as a statement about a Johnson graph. The quotient
family at \(v\) is shifted, and the quotient adjacency is the Johnson
adjacency; both are established in Section 5.

\hypertarget{the-y-family-pivots-form-one-degree-class}{%
\subsection{8.2 The Y-family pivots form one degree
class}\label{the-y-family-pivots-form-one-degree-class}}

Everything below is about \(H(d)\), defined in Section 8.1: the pivots
whose quotient family is a Y-family, which is the half of exceptionality
that \(d\) decides on its own.

\textbf{The idea, before the notation.} Suppose the quotient at \(v\) is
a Y-family. Theorem 7.6 says its two arms are nonempty and contiguous,
and that is a statement about which neighborhoods of \(v\) are
realizable: four particular sets are, two are not. By Lemma 5.2
realizability is graphicality of a residual, so the clique-sum shape
hands us four residual degree functions that are graphical and two that
are not. Six residuals differing only in two coordinates cannot be
graphical and non-graphical in that pattern unless the degrees at those
coordinates are pinned. Pinning them is exactly what the theorem below
concludes, and everything after it turns the pattern into the pinning,
in the currency of Erdős--Gallai slacks.

The theorem therefore confines every such pivot to a single degree
class, and pins that class to exactly three vertices --- so the final
argument has three pivots to count rather than a sequence to search.

\textbf{Theorem 8.3.} \emph{If \(v ∈ H(d)\), then with \(k = d(v)\),
\(a = u_k\) and \(b = u_{k+1}\) we have \(d(a) = d(b) = k\), and
\(\{v, a, b\}\) is exactly the degree-\(k\) class of \(d\).}

\textbf{Read on the sorted sequence, Theorem 8.3 says this.} Write \(d\)
non-increasing as \(d_1 ≥ d_2 ≥ ⋯ ≥ d_n\). If \(v ∈ H(d)\) and
\(k = d(v)\), then

\[
d_k = d_{k+1} = d_{k+2} = k :
\]

the value \(k\) occurs exactly three times in \(d\), and its run begins
at position \(k\). Theorem 8.3 makes the degree-\(k\) class exactly
\(\{v, u_k, u_{k+1}\}\), so \(k\) occurs three times; equal degrees are
consecutive when sorted, so the three sit at positions \(j, j+1, j+2\)
for some \(j\); and since the order \(u_1 ≻ u_2 ≻ ⋯\) deletes \(v\),
which is one of the three, exactly \(j−1\) larger degrees precede the
other two and place them at \(j\) and \(j+1\), which Theorem 8.3 puts at
\(k\) and \(k+1\).

The proof of Theorem 8.1 does not require this sorted formulation, and
it is recorded because it is the form a reader can apply: no order on
the ground, no deletion of \(v\), no family, just a look at position
\(k\) for three copies of \(k\). The two smallest instances read
\(d = (3,2,2,2,1)\) with \(d_2 = d_3 = d_4 = 2\) and
\(d = (4,4,3,3,3,1)\) with \(d_3 = d_4 = d_5 = 3\). It is necessary and
far from sufficient --- at ground order eight it admits 174 pivots of
which 102 lie in \(H(d)\) --- but it disposes of the rest at a glance.

The six lemmas below divide as follows, and it is worth having the shape
in mind before the notation arrives. Lemmas 8.3a to 8.3d are bookkeeping
about Erdős--Gallai slacks: 8.3a fixes the two strict outer gaps, 8.3b
and 8.3c confine where a violation can sit, and 8.3d --- much the
longest --- propagates any violation below rank \(k−1\) up to rank
\(k−1\). What that bookkeeping delivers is a single inequality, and
Lemma 8.3e converts it into the one-sided bound \(d(b) ≥ k\). Lemma 8.3f
is the same argument run on the complement, giving \(d(a) ≤ k\). With
\(d(b) ≤ d(a)\) from sortedness the two bounds close as
\(k ≤ d(b) ≤ d(a) ≤ k\), so \(d(a) = d(b) = k\), which is Theorem 8.3.

Write \(C_0 = \{u_1, …, u_{k−2}\}\) and name the four vertices that
follow it by their position in the order:

\[
c = u_{k−1}, \quad a = u_k, \quad b = u_{k+1}, \quad x = u_{k+2},
\]

so that sortedness of the order reads

\[
d(c) ≥ d(a) ≥ d(b) ≥ d(x).
\]

\textbf{Every degree in this section is written \(d(·)\) of the vertex
it belongs to}, so no separate letters have to be matched to vertices.
The four keep the names they arrive with: \(a\) and \(b\) are the forced
pair of Lemma 8.2 --- the neighbor \(v\) takes on and the one it gives
up, the same \(a\) and \(b\) Section 5 swaps across an interface ---
while \(c\) and \(x\) are the order-neighbors bracketing them.

\textbf{The letters on the smallest live instance.} Take
\(d = (4,2,2,2,1,1)\) on \(\{1,…,6\}\), the sequence Section 10 also
runs. It has nine realizations, and the pivots whose quotient is a
Y-family are \(2\), \(3\) and \(4\). Take \(v = 2\). Then
\(k = d(v) = 2\), the order on the remaining ground is
\(1 ≻ 3 ≻ 4 ≻ 5 ≻ 6\) by non-increasing degree, and the four names land
on

\[
C_0 = ∅, \quad c = 1, \quad a = 3, \quad b = 4, \quad x = 5.
\]

The realizable neighborhoods of \(v\) are \(\{1,3\}\), \(\{1,4\}\),
\(\{1,5\}\), \(\{1,6\}\) and \(\{3,4\}\). Reading the six corners
against that list:

\begin{longtable}[]{@{}lll@{}}
\toprule\noalign{}
corner & the set & realizable \\
\midrule\noalign{}
\endhead
\bottomrule\noalign{}
\endlastfoot
\(\{c,a\}\) & \(\{1,3\}\) & yes \\
\(\{c,b\}\) & \(\{1,4\}\) & yes \\
\(\{c,x\}\) & \(\{1,5\}\) & yes \\
\(\{a,b\}\) & \(\{3,4\}\) & yes \\
\(\{a,x\}\) & \(\{3,5\}\) & \textbf{no} \\
\(\{b,x\}\) & \(\{4,5\}\) & \textbf{no} \\
\end{longtable}

Four realizable and two not, which is the configuration the rest of this
subsection is about. In residuals, \(D = f_{a,x}\) is \((4,1,2,0,1)\),
which is not graphical --- its largest entry is \(4\) while the others
sum to only \(4\) --- whereas \(U = f_{c,x} = (3,2,2,0,1)\) and
\(V = f_{a,b} = (4,1,1,1,1)\), the star, both are. And the conclusion to
be proved is visible in the sequence itself: the degree-\(2\) class is
\(\{2,3,4\}\), exactly three vertices, and sorted \(d\) reads
\((4,2,2,2,1,1)\) with \(d_2 = d_3 = d_4 = 2\).

Two standing assumptions are implicit in these names and we record them:
\(k ≥ 2\), or \(c = u_{k−1}\) does not exist, and \(k + 2 ≤ n − 1\), or
\(x = u_{k+2}\) does not exist. Both hold whenever \(v\) is a Y-family
pivot, since a Y-family has at least four members.

Because both Y-family arms are nonempty and contiguous (Theorem 7.6),
the four sets \(C_0 ∪ \{c,a\}\), \(C_0 ∪ \{c,b\}\), \(C_0 ∪ \{c,x\}\)
and \(C_0 ∪ \{a,b\}\) are realizable neighborhoods of \(v\), while
\(C_0 ∪ \{a,x\}\) and \(C_0 ∪ \{b,x\}\) are not, lying in neither arm.
This configuration is the \textbf{four-corner hypothesis}. All six, and
where each comes from:

\begin{longtable}[]{@{}
  >{\raggedright\arraybackslash}p{(\columnwidth - 4\tabcolsep) * \real{0.3333}}
  >{\raggedright\arraybackslash}p{(\columnwidth - 4\tabcolsep) * \real{0.3333}}
  >{\raggedright\arraybackslash}p{(\columnwidth - 4\tabcolsep) * \real{0.3333}}@{}}
\toprule\noalign{}
\begin{minipage}[b]{\linewidth}\raggedright
neighborhood of \(v\)
\end{minipage} & \begin{minipage}[b]{\linewidth}\raggedright
realizable
\end{minipage} & \begin{minipage}[b]{\linewidth}\raggedright
how we know
\end{minipage} \\
\midrule\noalign{}
\endhead
\bottomrule\noalign{}
\endlastfoot
\(C_0 ∪ \{c,a\}\) & yes & shiftedness, from \(C_0 ∪ \{a,b\}\) \\
\(C_0 ∪ \{c,b\}\) & yes & shiftedness, from \(C_0 ∪ \{a,b\}\) \\
\(C_0 ∪ \{c,x\}\) & yes & assumed \\
\(C_0 ∪ \{a,b\}\) & yes & assumed \\
\(C_0 ∪ \{a,x\}\) & no & assumed \\
\(C_0 ∪ \{b,x\}\) & no & shiftedness, or the transposition when
\(d(a) = d(b)\) \\
\end{longtable}

\textbf{Three of those six conditions suffice, and the formal
development assumes exactly those three} --- the rows the table marks as
assumed. The other three follow, and the paragraphs below derive them.
The name counts the realizable corners, not the ones that do work.
Everything from Lemma 8.3a onward runs on three residuals ---
\(f_{a,x}\), \(f_{c,x}\) and \(f_{a,b}\), named \(D\), \(U\) and \(V\)
below --- which are exactly the three corners the hypothesis assumes.
The other two realizable corners, \(C_0 ∪ \{c,a\}\) and
\(C_0 ∪ \{c,b\}\), are derived here and then not used again; the table
displays all six because its purpose is to show the pattern is closed,
which is what makes the complement of Lemma 8.3f available.

Two of them follow from shiftedness of \(I_v\) (Theorem 5.3). Here
\(c = u_{k−1}\) precedes \(a = u_k\) and \(b = u_{k+1}\) by
construction, so shiftedness exchanges either for \(c\) in a realizable
neighborhood: from \(C_0 ∪ \{a,b\}\) it delivers \(C_0 ∪ \{c,a\}\) and
\(C_0 ∪ \{c,b\}\) at once. The sixth runs the same way with a second
mechanism available when the degrees tie. Suppose \(C_0 ∪ \{b,x\}\) were
realizable. If \(d(b) < d(a)\) then shiftedness replaces \(b\) by \(a\);
if \(d(b) = d(a)\) then transposing the two labels fixes \(d\), hence
carries realizable neighborhoods to realizable neighborhoods, and
carries \(C_0 ∪ \{b,x\}\) to \(C_0 ∪ \{a,x\}\). Either way
\(C_0 ∪ \{a,x\}\) would be realizable, contrary to the one
non-realizability we assume.

\textbf{One step of that is available here and not below, and the
difference matters.} The derivation just given places \(c\) before \(a\)
and \(b\) by their \emph{indices} in the pivot setup. The statements
from Lemma 8.3a onward are deliberately about the configuration and not
about the pivot --- see the remark before 8.3a --- and at that level all
that is assumed is \(d(c) ≥ d(a) ≥ d(b) ≥ d(x)\), which is \textbf{not
strict} and so does not by itself place \(c\) first. The formalization
reaches the same conclusion from \textbf{Lemma 8.3a}, which gives
\(d(c) > d(a)\) and hence \(d(b) ≤ d(a) < d(c)\); and 8.3a consumes only
the three assumed conditions, so nothing is circular.

We state all six in the display because that is how the Y-family arms
deliver them. All six are carried across the complement --- it sends a
pair to its complement within \(\{c,a,b,x\}\) and renames by
\((c,a,b,x) ↦ (x,b,a,c)\), so \(\{c,x\}\) and \(\{a,b\}\) exchange while
\(\{c,a\}\), \(\{c,b\}\), \(\{a,x\}\) and \(\{b,x\}\) are each fixed ---
which is what lets Lemma 8.3f run the whole argument on the complement.
Assuming three rather than six does not create that closure; it makes
the assumed set \textbf{minimally} closed, so nothing is carried through
the complement except what a later lemma reads.

Write \(F_Z(t)\) for the Erdős--Gallai {[}13{]} slack of a sequence
\(Z\) at rank \(t\), the right side minus the left side:

\[
F_Z(t) = t(t−1) + Σ_{i>t} \min(z_i, t) − Σ_{i≤t} z_i,
\]

where \(z_1 ≥ z_2 ≥ ⋯\) are the sorted entries of \(Z\). Thus \(Z\) is
graphical if and only if it is non-negative, has even sum, and
\(F_Z(t) ≥ 0\) for every \(t\). Only the \(⟸\) half is taken from the
citation. The \(⟹\) half --- that a graphical sequence has non-negative
slack at every rank --- is proved below by a double count, inside Step 3
of Lemma 8.3d, and it is used before it is proved there; the double
count depends on nothing in this section, so there is no circularity.
The three residuals the hypothesis names are the three that carry the
argument, and we name them once:

\[
D = f_{a,x}, \quad U = f_{c,x}, \quad V = f_{a,b},
\]

where \(D\) is the forbidden residual and the other two are graphical.
Set \(σ = F_D(k−1)\). Two choices are made in that line and both are
forced. The slack is read at the forbidden corner because that is the
only one the hypothesis makes non-graphical, and non-graphicality is
what will later supply a violating rank. And it is read at rank \(k−1\)
because, by Lemma 8.3b below, that is the rank up to which all six
residuals fill their sorted prefix from the same entries --- the \(C_0\)
values together with \(c\)'s --- so that only \(c\)'s contribution
varies from corner to corner, and it varies by one. That is what makes
the slack there computable in closed form and, more to the point,
comparable across corners. Everything after this turns on \(σ\): both
closed-form identities, the conclusion of Lemma 8.3d, and the top-level
case split of Lemma 8.3e. In \(D\) the vertices \(a\) and \(x\) are
decremented, so its four middle values are
\((d(c), d(a)−1, d(b), d(x)−1)\) at \((c, a, b, x)\); in \(U\) they are
\((d(c)−1, d(a), d(b), d(x)−1)\) and in \(V\) they are
\((d(c), d(a)−1, d(b)−1, d(x))\).

The next six statements are about the four-corner hypothesis itself, not
about the pivot \(v\) we started from. We state them that way
deliberately: Lemma 8.3f applies Lemma 8.3e to a \textbf{different}
degree sequence, the complement, which need not be the quotient of a
Y-family pivot of \(d\). Reading them as statements about \(v\) would be
a quantifier slip of exactly the kind that has broken an argument in
this line of work before.

\textbf{Lemma 8.3a (strict outer gaps).} \emph{In every four-corner
configuration, \(d(c) > d(a)\) and \(d(b) > d(x)\).}

\emph{Proof.} If \(d(c) = d(a)\), transposing \(c\) and \(a\) fixes
\(d\) and sends the graphical corner \(C_0 ∪ \{c,x\}\) to the
non-graphical \(C_0 ∪ \{a,x\}\); graphicality depends only on the
multiset of residual values, so this is a contradiction. If
\(d(b) = d(x)\), transposing \(b\) and \(x\) sends the graphical
\(C_0 ∪ \{a,b\}\) to the non-graphical \(C_0 ∪ \{a,x\}\).
\(\blacksquare\)

\textbf{Lemma 8.3b (stability at \(k−1\)).} \emph{In every four-corner
configuration, for all six residuals the top \(k−1\) sorted entries are
exactly the \(C_0\) values together with \(c\)'s value.}

\emph{Proof.} Each \(C_0\) entry is at least \(d(c) − 1\), and \(c\)'s
value is \(d(c)\) or \(d(c) − 1\), while every entry at position \(k\)
or later is at most \(d(a) ≤ d(c) − 1\) by Lemma 8.3a. A tie at the
boundary can exchange labels but not values, and only the multiset of
values is asserted. \(\blacksquare\)

Stability makes the slack at \(k−1\) computable in closed form. Every
corner \(P\) contains \(C_0\), contains exactly two of \(c, a, b, x\),
and contains no \(u_i\) with \(i ≥ k+3\), so its residual assigns
\(e_i − 1\) to each \(u_i ∈ C_0\), the value \(d(c) − [c ∈ P]\) to
\(c\), and likewise for \(a, b, x\), leaving \(e_i\) untouched for
\(i ≥ k+3\). By Lemma 8.3b the first \(k−1\) sorted entries sum to
\(L_0 + d(c) − [c ∈ P]\) where \(L_0 = Σ_{i ≤ k−2} (e_i − 1)\), and
everything else lies in the tail. Writing \(m(y) = \min(y, k−1)\) and
\(L_1 = Σ_{i ≥ k+3} m(e_i)\), the tail sum beside the prefix sum
\(L_0\), the definition of \(F\) gives, for every corner \(P\),

\[
\begin{gathered}
F_P(k−1) = (k−1)(k−2) + m(d(a) − [a ∈ P]) + m(d(b) − [b ∈ P]) \\
+ m(d(x) − [x ∈ P]) + L_1 − L_0 − (d(c) − [c ∈ P]).
\end{gathered}
\]

Only the four bracketed indicators change from one corner to another, so
subtracting the expression for \(\{a,x\}\) from those for \(\{c,x\}\)
and \(\{a,b\}\), and using \(m(y) − m(y−1) = [y ≤ k−1]\), gives the two
identities

\[
F_{c,x}(k−1) = σ + 1 + [d(a) ≤ k−1], \quad F_{a,b}(k−1) = σ − [d(b) ≤ k−1] + [d(x) ≤ k−1].
\]

\textbf{Lemma 8.3c (agreement above the pivot rank).} \emph{In every
four-corner configuration \(F_D(j) = F_U(j)\) for every \(j ≥ k+1\); and
if \(d(a) > d(b)\) then also \(F_D(k) = F_U(k)\).}

\emph{Proof.} Passing from \(D\) to \(U\) replaces the pair of values
\((d(c), d(a)−1)\) by \((d(c)−1, d(a))\), a unit transfer, and changes
nothing else. In \(D\) the only entries that can strictly exceed
\(d(a)−1\) are the \(k−2\) entries of \(C_0\), the value \(d(c)\) at
\(c\), and --- only when \(d(b) = d(a)\) --- the value \(d(b)\) at
\(b\); so the changed \(a\)-value has sorted rank at most \(k+1\). In
\(U\) both changed values likewise have rank at most \(k+1\). Hence by
rank \(k+1\) both ends of the transfer lie inside the sorted prefix, so
the two prefixes have equal sums and the two tails have equal multisets,
and every term of \(F\) agrees. If \(d(a) > d(b)\) then \(b\) does not
exceed \(d(a)−1\), the ranks drop by one, and the same argument applies
already at rank \(k\). \(\blacksquare\)

We record what Lemma 8.3c does \textbf{not} say. It is not true that all
six slacks agree above rank \(k+1\): for \(d = (4,2,2,2,1,1)\) with
\(k = 2\) the rank-4 slacks of the six corners are \(6, 6, 4, 6, 4, 4\).
Only the \(D\)-to-\(U\) comparison is claimed, and only that is used.

\textbf{Lemma 8.3d (lower-witness propagation).} \emph{In every
four-corner configuration, if \(F_D(t) < 0\) for some \(t ≤ k−2\), then
\(F_D(k−1) < 0\).}

\emph{Proof.} Fix such a \(t\).

\emph{Step 1: the lower transfer has defect exactly one.} First, the
top-\(t\) sums of \(D\) and \(V\) agree. The \(k−1\) entries carried by
\(C_0 ∪ \{c\}\) are identical in \(D\) and in \(V\): every member of
\(C_0\) is decremented by both corners, since \(C_0\) lies in all six of
them, and neither corner's distinguishing pair meets \(C_0 ∪ \{c\}\), so
the block is decremented the same way on both sides. All of them are at
least \(d(c) − 1\), since a \(C_0\) entry is \(e_i − 1 ≥ d(c) − 1\) and
\(c\)'s entry is \(d(c)\). Every remaining entry of either sequence is
at most \(d(c) − 1\): in \(D\) these are \(d(a) − 1 ≤ d(c) − 2\),
\(d(b) ≤ d(a) ≤ d(c) − 1\) and \(d(x) − 1\), in \(V\) they are
\(d(a) − 1\), \(d(b) − 1\) and \(d(x) ≤ d(b) ≤ d(c) − 1\), and each tail
entry is at most \(d(x)\). Since \(t ≤ k − 2 < k − 1\), the common block
\(C_0 ∪ \{c\}\) already supplies more than \(t\) entries at the top, and
an outside entry can enter the top \(t\) only by tying at \(d(c) − 1\),
which changes labels but not values. So both top-\(t\) multisets equal
the top-\(t\) multiset of \(C_0 ∪ \{c\}\), and in particular the two
sums agree.

We do \textbf{not} claim the top \(t\) entries are
\(d(c), e_1 − 1, …, e_{t−1} − 1\). That description is false in general
--- for \(d = (6,4,3,3,3,2,1)\) with \(k = 3\) and \(t = 1\) it offers
\(d(c) = 4\) where the true top entry is \(e_1 − 1 = 5\) --- and it
becomes correct only once \(e_t = d(c)\), which is Step 2's conclusion.
Only the equality of sums is used here, and it holds unconditionally.

The only values differing between \(D\) and \(V\) are \((d(b), d(x)−1)\)
against \((d(b)−1, d(x))\), so in the tail

\[
F_V(t) − F_D(t) = [d(x) ≤ t] − [d(b) ≤ t].
\]

\(V\) is graphical and \(D\) violates at \(t\), so the left side is at
least one; the right side is at most one. Hence both equal one, and

\[
d(x) ≤ t < d(b) ≤ d(a) < d(c), \quad F_D(t) = −1, \quad F_V(t) = 0. \quad (8.1)
\]

\emph{Step 2: the upper transfer forces an equality in \(U\).} By (8.1)
all four of \(d(c), d(c)−1, d(a), d(a)−1\) are at least \(t\), so each
caps to \(t\). Capping every entry at \(t\) therefore sends \(D\) and
\(U\) to the same multiset: they differ only in the two entries at \(c\)
and \(a\), and all four values occurring there cap to \(t\). Equal
capped multisets have equal sorted sequences, so the capped tail sum
\(Σ_{i>t} \min(z_i, t)\) is the same for \(D\) and \(U\). The transfer
\(D → U\) replaces \((d(c), d(a)−1)\) at \((c,a)\) by
\((d(c)−1, d(a))\), so the two sequences differ only there. Let

\[
g = \#\{ i ≤ k−2 : e_i > d(c) \},
\]

the number of \(C_0\) entries strictly above \(d(c) − 1\); it is the
count that decides whether the top \(t\) of \(D\) and of \(U\) is drawn
entirely from untouched \(C_0\) entries, which is the case split of the
next paragraph. Since \(d(a) − 1 ≤ d(c) − 2\) and \(d(b) ≤ d(c) − 1\),
the entries strictly above \(d(c) − 1\) are exactly those \(g\) together
with \(c\)'s entry in \(D\), and exactly those \(g\) in \(U\).

If \(t ≤ g\) the top \(t\) of both sequences is drawn from the same
\(g\) untouched \(C_0\) entries, so the sums agree; and \(e_t > d(c)\),
since \(e_t\) is one of them. If \(t > g\), the top \(t\) of \(D\) is
those \(g\) entries, then \(c\)'s value \(d(c)\), then \(t − g − 1\)
entries equal to \(d(c) − 1\); the top \(t\) of \(U\) is those \(g\)
entries and then \(t − g\) entries equal to \(d(c) − 1\). \(U\) has
enough of them: the \(k − 2 − g\) \(C_0\) entries equal to \(d(c) − 1\)
together with \(c\)'s new value \(d(c) − 1\) give \(k − 1 − g ≥ t − g\),
using \(t ≤ k − 1\). So the sums differ by exactly
\(d(c) − (d(c)−1) = 1\); and \(e_t ≤ d(c)\), hence \(e_t = d(c)\), since
\(e_t ≥ e_{k−1} = d(c)\). Writing \(ε = [e_t = d(c)]\), in both cases

\[
F_U(t) − F_D(t) = ε,
\]

and the difference is exactly \(ε\), not merely at least it; this
matters, because Step 3 needs Erdős--Gallai \emph{equality} in \(U\) and
not just non-negativity. As \(U\) is graphical and \(F_D(t) = −1\), this
forces \(ε = 1\) and

\[
e_t = e_{t+1} = ⋯ = e_{k−1} = d(c), \quad F_U(t) = 0. \quad (8.2)
\]

\emph{Step 3: what equality in Erdős--Gallai forces.} Let
\(y_1 ≥ y_2 ≥ ⋯\) be the sorted entries of \(U\). By (8.2) and the
top-prefix description, \(y_t = d(c) − 1\).

We first record what equality in Erdős--Gallai says about a realization.
Take any graph, and let \(T\) be any set of \(t\) of its vertices.
Counting the edges leaving \(T\),

\[
\begin{gathered}
Σ_{u ∈ T} \deg(u) = Σ_{u ∈ T} |N(u) ∩ T| + Σ_{u ∉ T} |N(u) ∩ T| \\
≤ t(t−1) + Σ_{u ∉ T} \min(\deg(u), t),
\end{gathered}
\]

where both estimates are termwise: \(N(u) ∩ T ⊆ T ∖ \{u\}\) for
\(u ∈ T\), since \(u ∉ N(u)\), and \(N(u) ∩ T\) is contained in both
\(N(u)\) and \(T\) for \(u ∉ T\). Taking \(T\) to consist of \(t\)
vertices of largest degree makes the two sides the two sides of
Erdős--Gallai at rank \(t\), so the slack \(F(t)\) is the sum of these
termwise slacks. \textbf{Equality therefore forces every one of them to
be tight}: \(T\) is a clique, and every vertex outside \(T\) meets \(T\)
in exactly \(\min(\deg, t)\) vertices, so every outside vertex of degree
at least \(t\) is adjacent to all of \(T\).

\emph{Step 4: a union of tight sets.} Write \(H = C_0 ∪ \{c\}\), so
\(|H| = k−1\), and fix a realization \(R\) of \(U\), graphical by
hypothesis. Let

\[
K = \{ u ∈ H : \deg_R(u) > d(c) − 1 \}, \quad g = |K|.
\]

\emph{First, \(g < t\).} The \(t\)-th largest entry of \(U\) is
\(y_t = d(c) − 1\), so at most \(t−1\) entries exceed \(d(c) − 1\).
Hence \(g ≤ t−1 < t\), and every vertex of \(H ∖ K\) has degree exactly
\(d(c) − 1\).

\emph{A family of tight sets.} For each \((t−g)\)-element subset
\(J ⊆ H ∖ K\), put \(T_J = K ∪ J\), so \(|T_J| = t\). Each \(T_J\)
consists of \(t\) vertices of largest degree in \(R\): it contains every
vertex above \(d(c) − 1\), and all its remaining vertices sit at
\(d(c) − 1\), which no vertex outside \(H\) exceeds. By (8.2) the slack
\(F_U(t)\) is zero, so \textbf{Step 3 applies to every \(T_J\)}: each is
a clique, and every vertex outside it of degree at least \(t\) is
adjacent to all of it.

\emph{Their union is tight.} Since \(1 ≤ t−g ≤ |H ∖ K|\), every element
of \(H ∖ K\) lies in some \(J\), so the \(T_J\) have union \(H\); and
every vertex of \(H\) has degree at least \(d(c) − 1 ≥ d(a) > t\) by
(8.1). That is all it takes. Given \(h ≠ h'\) in \(H\), choose \(J\)
with \(h ∈ T_J\): either \(h' ∈ T_J\), which is a clique, or
\(h' ∉ T_J\) has degree at least \(t = |T_J|\) and is therefore adjacent
to all of \(T_J\). Either way \(h h'\) is an edge, so \(H\) is a clique.
A vertex outside \(H\) of degree at least \(t\) is adjacent to every
\(T_J\) and hence to all of \(H\); one of degree below \(t\) has all its
neighbors inside any single \(T_J\), hence inside \(H\), so it meets
\(H\) in \(\min(\deg, |H|)\) vertices either way. Those are exactly the
equality conditions of Step 3, so \(H\) is tight, and \(|H| = k−1\)
gives

\[
F_U(k−1) = 0. \quad (8.3)
\]

\emph{Conclusion.} Lemma 8.3b at rank \(k−1\) reads
\(F_U(k−1) = F_D(k−1) + 1 + [d(a) ≤ k−1]\). With (8.3) this gives
\(F_D(k−1) = −1 − [d(a) ≤ k−1] < 0\). \(\blacksquare\)

\textbf{What this argument replaces, and why the set form is the right
one.} Earlier drafts propagated the violation upward through a plateau
of tied sorted values, by a recurrence on consecutive slacks. The
saturation above is a statement about a \emph{set of vertices} rather
than about a sequence of numerical slacks, and equality at a tied prefix
can therefore be read at \textbf{every} choice of tied labels at once.
Taking the union of those choices jumps straight from rank \(t\) to rank
\(k−1\), with no recurrence and no degree count in between.

\textbf{Lemma 8.3e (one-sided bound).} \emph{Every four-corner
configuration satisfies \(d(b) ≥ k\).}

\emph{Proof.} Suppose first \(σ < 0\). If \(d(a) ≤ k−1\) then
\(d(b) ≤ d(a) ≤ k−1\) and \(d(x) ≤ k−1\), so \(F_{a,b}(k−1) = σ < 0\),
contradicting that \(\{a,b\}\) is a realizable corner; so \(d(a) ≥ k\).
Then \(F_{c,x}(k−1) = σ + 1 ≥ 0\) forces \(σ = −1\), and
\(F_{a,b}(k−1) = −1 − [d(b) ≤ k−1] + [d(x) ≤ k−1] ≥ 0\) forces
\(d(b) ≥ k\).

Suppose instead \(σ ≥ 0\). First, \(D\) is non-negative and has even
sum: its decremented \(a\)-entry also occurs in the graphical residual
\(V\), its decremented \(x\)-entry occurs in \(U\), and it has the same
sum as \(U\). So its non-graphicality supplies an Erdős--Gallai
violating rank rather than a parity or negative-entry obstruction, which
is what the argument needs. Every rank at least \(k+1\) is excluded by
Lemma 8.3c, since \(U\) is graphical; rank \(k−1\) is excluded by the
supposition \(σ ≥ 0\); and every rank at most \(k−2\) is excluded by
Lemma 8.3d, again using \(σ ≥ 0\). So the violating rank is \(k\), that
is \(F_D(k) < 0\). Were \(d(a) > d(b)\), Lemma 8.3c would give
\(F_D(k) = F_U(k) ≥ 0\); so \(d(a) = d(b)\).

When \(d(a) = d(b)\), the sorted \(k\)-th value of \(D\) is exactly
\(d(a)\), since the \(k−2\) entries of \(C_0\), the value \(d(c)\), and
the undecremented \(b\)-value supply its top \(k\) while every other
entry is at most \(d(a) − 1\). If \(d(a) ≤ k−1\), consecutive-slack
subtraction gives \(F_D(k) − F_D(k−1) = 2(k−1−d(a)) ≥ 0\), contradicting
\(F_D(k) < 0 ≤ σ\). So \(d(a) = d(b) ≥ k\). \(\blacksquare\)

\textbf{Lemma 8.3f (complementation).} \emph{Every four-corner
configuration satisfies \(d(a) ≤ k\).}

\emph{Proof.} Put \(N = n − 1\). Complementing every neighborhood sends
\(k\) to \(N − k\), reverses the four positional roles as
\((c,a,b,x) ↦ (x,b,a,c)\), and replaces each residual value by \(N − 1\)
minus it. Graphicality is preserved and each corner maps to the
corresponding corner, so the four-corner hypothesis is preserved,
forbidden corner included. Applying Lemma 8.3e to the complement gives
\(N − d(a) ≥ N − k\). \(\blacksquare\)

\emph{Proof of Theorem 8.3.} Lemmas 8.3e and 8.3f give
\(d(a) = d(b) = k\), which is \(d(a) = d(b) = k\). Lemma 8.3a then gives
\(d(c) > k > d(x)\), so every vertex before \(c\) has degree above \(k\)
and every vertex after \(x\) has degree below \(k\). With \(v\) itself
of degree \(k\), the class is exactly \(\{v, a, b\}\). \(\blacksquare\)

We note that this argument uses complementation on the \textbf{degree
bound} only, and not on any correspondence between individual violating
thresholds.

\textbf{Corollary 8.4.} \emph{If \(v ∈ H(d)\) with \(k = d(v)\), then
\(H(d) ⊆ \{v, a, b\}\), the degree-\(k\) class of \(d\); and the forced
pair \(g(v)\) is \(\{a, b\}\).}

\emph{Proof.} Let \(v ∈ H(d)\) with \(k = d(v)\). Theorem 8.3 puts the
whole degree-\(k\) class at \(\{v, a, b\}\), so \(H(d)\) contains at
most those three vertices of degree \(k\). It contains no vertex of
another degree: if \(w ∈ H(d)\) had \(d(w) = k' > k\), then applying
Theorem 8.3 at \(w\) would put exactly \(k'−1\) vertices above degree
\(k'\) and three at it, while applying it at \(v\) bounds the number of
vertices of degree at least \(k'\) by \(k−1\), giving \(k'+2 ≤ k−1\) and
so \(k' ≤ k−3\), a contradiction; the case \(k' < k\) is the same count
read downward. So \(H(d) ⊆ \{v, a, b\}\). That \(g(v) = \{a, b\}\) is
the definition of the forced pair, \(a\) and \(b\) being \(u_k\) and
\(u_{k+1}\). \(\blacksquare\)

We state the corollary as an \textbf{inclusion} on purpose. The reverse,
that \(a\) and \(b\) are themselves Y-family pivots and so \(H(d)\) is
empty or has exactly three elements, is true and follows from Lemma
8.5's relabeling argument applied to the quotient rather than the fiber,
but it is not proved here and nothing below uses it.

\textbf{Where an obstruction can live, read off \(d\) alone.} Section
8.4 proves Theorem 8.1 without this corollary, working only inside the
one degree class it is handed. The corollary is kept because of what it
tells a reader who is holding a degree sequence and nothing else. Put it
beside the sorted form of Theorem 8.3, which says an obstruction of
degree \(k\) forces \(d_k = d_{k+1} = d_{k+2} = k\). At most one index
can meet even the first of those three conditions: if \(d_k = k\) and
\(d_j = j\) with \(k < j\), then \(d_k ≥ d_j\) gives \(k ≥ j\).
\textbf{A non-increasing sequence therefore names its own candidate
degree, and names at most one.} The corollary supplies the other half,
that the pivots which can be exceptional do not spread across two
degrees, so there is nowhere else to look.

The test is a glance. Scan \(d\) for the single index with \(d_k = k\).
If there is none, or if \(d_{k+1}\) and \(d_{k+2}\) are not also \(k\),
then \(H(d)\) is empty: no pivot of \(d\) is exceptional at all, and
nothing in this section has anything to do. Otherwise every exceptional
pivot of \(d\) lies among the three vertices at positions
\(k, k+1, k+2\), and Section 8.4 shows that for a given pair at most two
of those three can be exceptional at once. A degree-\(2\) obstruction
must be three \(2\)s beginning at position \(2\), a degree-\(3\)
obstruction three \(3\)s beginning at position \(3\), and so on up; on
\(d = (3,2,2,2,1)\) the scan returns \(k = 2\), and on
\(d = (4,4,3,3,3,1)\) it returns \(k = 3\) with the run at positions
\(3, 4, 5\).

\hypertarget{admissibility-is-constant-on-a-degree-class}{%
\subsection{8.3 Admissibility is constant on a degree
class}\label{admissibility-is-constant-on-a-degree-class}}

\textbf{Lemma 8.5.} \emph{If \(d(u) = d(w)\) then \(u\) has a
non-bipartite fiber if and only if \(w\) does.}

\emph{Proof.} The transposition \(τ = (u\,w)\) satisfies \(d ∘ τ = d\),
so it permutes the realizations of \(d\). A 2-switch is equivariant
under relabeling, so \(τ\) induces an automorphism of \(G(d)\). It
carries a triangle to a triangle and the ground support of a triangle to
the image of that support. A triangle inside a \(u\)-fiber is one whose
support avoids \(u\); its image avoids \(w\) and lies inside a
\(w\)-fiber. \(\blacksquare\)

Lemma 8.5 does real work rather than holding vacuously. Y-family pivots
without a non-bipartite fiber exist: all three at \(d = (3,2,2,2,1)\).

\hypertarget{proof-of-theorem-8.1}{%
\subsection{8.4 Proof of Theorem 8.1}\label{proof-of-theorem-8.1}}

By Theorem 6.1 the pair has an admissible pivot \(v\). If \(v\) is not
exceptional we are done, so suppose it is, and write \(k = d(v)\). By
Theorem 8.3 the degree-\(k\) class is exactly \(\{v, a, b\}\), and by
Lemma 8.2 the two \(Δ\)-neighbors of \(v\) are \(a\) and \(b\), so both
separate the pair. Since \(v\) is admissible it has a non-bipartite
fiber, and \(a\) and \(b\) share its degree, so by Lemma 8.5 they have
non-bipartite fibers too: both are admissible.

Suppose both were exceptional as well. Theorem 8.3 and Lemma 8.2 then
apply at \(a\) and at \(b\), and the degree class each returns is the
same \(\{v, a, b\}\), because \(d(a) = d(b) = k\). So
\(N_Δ(a) = \{v, b\}\) and \(N_Δ(b) = \{v, a\}\), which together with
\(N_Δ(v) = \{a, b\}\) makes the \(Δ\)-component containing \(v\) exactly
the triangle on \(\{v, a, b\}\), with no longer trail passing through.

Now count colors. Each vertex meets as many \(G_0\)-only as \(G_1\)-only
edges of \(Δ\), because \(\deg_{G_0} = \deg_{G_1} = d\). A vertex of
\(Δ\)-degree two therefore meets exactly one edge of each color. Going
around the triangle the colors must alternate, which is impossible on a
cycle of odd length. So \(a\) and \(b\) are not both exceptional, and at
least one of them is an admissible non-exceptional pivot. One survives.
\(\blacksquare\)

The color count is used only at vertices of \(Δ\)-degree two. It does
not extend to larger degrees: taking \(G_0 = C_5\) and \(G_1\) its
complement on the same five vertices, both realize \((2,2,2,2,2)\) and
\(Δ = K_5\), which is balanced and full of triangles.

The argument is \textbf{local}: it never leaves the degree-\(k\) class
\(\{v, a, b\}\), so it needs nothing about exceptional pivots of any
other degree. It uses no hypothesis on \(d\) beyond what Theorem 6.1
needs to produce the first admissible pivot, and what it establishes is
slightly stronger than the statement: if any admissible pivot exists, an
admissible non-exceptional one exists.

\hypertarget{the-interface-extensions}{%
\section{9. The interface extensions}\label{the-interface-extensions}}

The assembly of Section 10 crosses from one fiber to the next along the
quotient path. The edges available for that are the \textbf{connectors}
of Section 5.4, and this section keeps that name: an edge of \(G(d)\)
with one endpoint in each of two adjacent fibers is exactly a connector,
since by the classification at the head of Section 5 a 2-switch that
leaves a fiber must move the pivot's neighborhood by one element. The
companion paper calls the same object a \emph{crossing edge} or an
\emph{interface edge}, and Section 7.8 uses \emph{crossing} for the
analogous edge one level down, between the two layers of a shifted
family; the two are the same idea in different categories and never
appear in one argument.

At each crossing the assembly needs a spanning walk of the current fiber
whose last realization carries a connector into the next fiber, and at
one crossing on each side it needs that connector to \textbf{avoid a
single prescribed target}. This section proves both.

\textbf{Why a target has to be forbidden.} The construction of Section
10 sends two chains inward, one from each prescribed realization, and
they meet at the buffer. If both chains were allowed to enter the buffer
at the same realization the two halves would collide there and no
spanning walk of the buffer could join them. So the second chain must
arrive somewhere the first did not, which is what Theorem 9.2 delivers:
a connector whose target differs from one named in advance. Section 10.2
is where the two arrival points are finally named.

\hypertarget{the-interface-1}{%
\subsection{9.1 The interface}\label{the-interface-1}}

Fix a pivot \(v\) and adjacent quotient vertices \(S\) and \(T\),
oriented as \(T = S − b + a\) with \(b ∈ S\) and \(a ∉ S\). Write
\(f = f_S\) for the residual degree function of the \(S\)-fiber, so that
\(f_T = f − e_a + e_b\), and set

\[
δ = f(a) − f(b).
\]

Connectors, their witnesses, the count \(w_{ab}(G)\) and the connector
graph between \(Φ_S\) and \(Φ_T\) are as defined in Section 5.4: a
\textbf{connector} from \(G ∈ Φ_S\) is the 2-switch
\(\{vb, ac\} → \{va, bc\}\) for a witness \(c ∉ \{v,a,b\}\) with
\(ac ∈ G\) and \(bc ∉ G\), and \(w_{ab}(G)\) counts them on the ground
realization \(F = G − v\). Every connector is a literal edge of
\(G(d)\), not merely a quotient adjacency, and \(G\) has exactly
\(w_{ab}(G)\) connector partners. We call \(G\) a \textbf{connector
source} when it has one. \textbf{Both counts are taken on the ground},
for the reason Section 5.4 gives: reading \(w_{ba}\) in \(G\) would
count the pivot and shift the identity below by one.

Three facts proved in Section 5 are used throughout what follows, and we
display them together because this is where they are consumed; Section
9.2 draws on Theorem 5.4 and Section 9.3 on Theorem 5.7 besides. The
signed identity of Theorem 5.6,

\[
w_{ab}(G) − w_{ba}(G) = f(a) − f(b),
\]

which in the target orientation reads
\(w_{ba}(H) − w_{ab}(H) = f_T(b) − f_T(a) = 2 − δ\) for \(H ∈ Φ_T\). The
exact color expansion of Theorem 5.8: if \(Φ_S\) is bipartite with color
classes \(C_0\) and \(C_1\), then for \(k = \max(1, δ)\) every
\(Λ ⊆ C_i\) satisfies \(|N(Λ)| = k|Λ|\) in the connector graph. Two
letters here carry a different sense in Section 8, and the sections
never meet: there \(k\) is the pivot degree \(d(v)\) and \(C_0\) is a
set of ground vertices, while throughout this section and in Theorems
5.7 and 5.8 \(k\) is \(\max(1, δ)\) and \(C_0\), \(C_1\) are color
classes of a fiber. And the difference taken \emph{across} a connector
rather than at a single realization,

\[
w_{ab}(G) − w_{ba}(H) = δ − 1,
\]

the same value for every connector of every such pair. Section 9.3 uses
the last at \(δ = 1\), where the two one-sided bounds give only \(≥ 1\)
and the argument needs \(≥ 2\).

Section 5.4's warning about the sign applies here unchanged: \(δ ≥ 0\)
is not a general fact about interfaces, it is asserted only for a
\textbf{triangle-free} source, and that is the only form any argument
below uses.

\hypertarget{ordinary-propagation}{%
\subsection{9.2 Ordinary propagation}\label{ordinary-propagation}}

\textbf{Theorem 9.1.} \emph{Let \(Φ_S\) be maximally Hamiltonian. For
every \(G ∈ Φ_S\) there is a connector source \(H ∈ Φ_S\) together with
a spanning \(G\)--\(H\) walk of \(Φ_S\).}

\emph{Proof.} Three cases.

If \(|Φ_S| = 1\), quotient adjacency gives a literal 2-switch between
the two fibers, so its only member is a connector source and the
one-vertex walk works.

If \(Φ_S\) is bipartite with \(|Φ_S| ≥ 2\), it is balanced and both
color classes are nonempty, by the remark on maximal Hamiltonicity in
Section 2. Let \(C\) be the class opposite \(G\). Exact color expansion
gives \(|N(C)| = k|C| ≥ 1\), so some \(H ∈ C\) is a connector source,
and Hamilton-laceability supplies the spanning \(G\)--\(H\) walk.

If \(Φ_S\) is non-bipartite, it contains an odd cycle, so \(|Φ_S| ≥ 3\),
and Hamilton-connectedness supplies a spanning walk to any \(H ≠ G\). It
is enough to find a connector source other than \(G\), and we split on
the sign of \(δ\). If \(δ ≥ 1\), the signed identity gives
\(w_{ab}(G') = w_{ba}(G') + δ ≥ 1\) for every \(G' ∈ Φ_S\), so all of
them are connector sources and at least two differ from \(G\). If
\(δ ≤ 0\), pick any \(H ∈ Φ_T\); the target identity gives
\(w_{ba}(H) = w_{ab}(H) + (2 − δ) ≥ 2\), and distinct witnesses give
distinct partners, so \(H\) alone has two distinct connector sources in
\(Φ_S\) and at least one differs from \(G\). \(\blacksquare\)

The split on \(δ\) is what lets Theorem 9.1 stand without any hypothesis
on the target fiber. That matters, because ordinary crossings occur at
every step of the quotient path and their targets are arbitrary.

\hypertarget{avoiding-one-target}{%
\subsection{9.3 Avoiding one target}\label{avoiding-one-target}}

\textbf{The idea, and where the difficulty sits.} Theorem 9.1 delivers
\emph{some} connector; here one target is forbidden, for the reason
given at the head of this section, and we need a second. The forbidden
one is where the other chain will enter the buffer, so a connector
landing there would make the two halves collide. Where the target fiber
has three or more members and the source has room to choose its exit,
that is easy --- pick a different exit and a different target. The whole
subsection is the accounting for when there is not much room: small
\(|Φ_S|\), and the interface parameter at \(δ = 1\), where each of the
two one-sided bounds gives only one connector and the argument needs
two. The constant difference across the interface, recorded at the end
of Section 5.4, is what supplies the second.

\textbf{Theorem 9.2.} \emph{Let \(Φ_S\) be maximally Hamiltonian and
\(Φ_T\) non-bipartite. For every \(G ∈ Φ_S\) and every \(z_0 ∈ Φ_T\)
there are \(H\), \(z\), \(c\) and a spanning \(G\)--\(H\) walk of
\(Φ_S\) such that \(c\) is a witness for \(H\), the connector carries
\(H\) to \(z ∈ Φ_T\), and \(z ≠ z_0\).}

\emph{Proof.} Since \(Φ_T\) is non-bipartite, \(t = |Φ_T| ≥ 3\). Put
\(r = |Φ_S|\).

Suppose \(Φ_S\) is non-bipartite, so \(r ≥ 3\) and
Hamilton-connectedness makes every \(H ≠ G\) a legal endpoint. It is
enough to produce a connector from some \(H ≠ G\) to a target other than
\(z_0\), and as in Theorem 9.1 we split on the sign of \(δ\). Fix two
distinct \(H_1, H_2 ∈ Φ_S ∖ \{G\}\), which \(r ≥ 3\) supplies.

If \(δ ≥ 2\), the signed identity gives
\(w_{ab}(H_1) = w_{ba}(H_1) + δ ≥ 2\), so \(H_1\) has two distinct
connector partners and at least one of them is not \(z_0\).

If \(δ ≤ 0\), pick any \(z ∈ Φ_T\) other than \(z_0\), which \(t ≥ 3\)
supplies. The target identity gives
\(w_{ba}(z) = w_{ab}(z) + (2 − δ) ≥ 2\), so \(z\) has two distinct
connector sources in \(Φ_S\) and at least one of them is not \(G\).

There remains \(δ = 1\), where both identities give only \(≥ 1\) and
neither side has a vertex known to carry two connectors. Here the order
of the choices matters: \(z_0\) is fixed first, and the endpoint is
chosen after. Each of \(H_1\) and \(H_2\) has \(w_{ab} ≥ 1\) and so has
at least one connector partner. If either partner is other than \(z_0\),
that pair is the one we want. Otherwise \(z_0\) is a connector partner
of both, so \(z_0\) has two distinct connector sources and \(w_{ba}\) at
\(z_0\) is at least \(2\). Differencing along the connector from \(H_1\)
to \(z_0\) then gives \(w_{ab}(H_1)\) equals \(w_{ba}\) at \(z_0\) plus
\(δ − 1\), hence equals it, so \(w_{ab}(H_1) ≥ 2\) after all, and
\(H_1\) has a second partner, necessarily other than \(z_0\).

Suppose \(Φ_S\) is bipartite and \(r ≥ 3\), the case \(r ≤ 2\) being
taken up below. It is balanced by the Section 2 remark, so \(r\) is even
and \(r ≥ 3\) gives \(r ≥ 4\). Let \(C\) be the class opposite \(G\), so
\(|C| = r/2 ≥ 2\). Exact color expansion gives \(|N(C)| = k|C| ≥ 2\), so
some connector from \(C\) has target other than \(z_0\), and laceability
supplies the walk.

There remains \(r ≤ 2\). The source is then triangle-free, so the
cross-degree theorem applies and gives both \(δ ≥ 0\) and that every
source has exactly \(k = \max(1, δ)\) connectors. (\(δ ≥ 0\) is not
automatic from graphicality of \(f_T\); see Section 5.4.) Suppose
\(k = 1\). Then \(δ ≤ 1\) and the interface has \(kr = r ≤ 2\) edges in
total. For every \(H ∈ Φ_T\) the target identity gives
\(w_{ba}(H) = w_{ab}(H) + 2 − δ ≥ 1\), so every member of \(Φ_T\) meets
an interface edge and \(t ≤ 2\), contradicting \(t ≥ 3\). So \(k ≥ 2\).
If \(r = 1\) the singleton walk at \(G\) has \(k ≥ 2\) distinct targets,
one of them other than \(z_0\); if \(r = 2\) the two members are
adjacent, so the fiber edge \(GH\) is a spanning walk to the other
member \(H\), and \(H\) again has \(k ≥ 2\) distinct targets.
\(\blacksquare\)

Theorem 9.2 is stated with the exclusion inside the same existential as
the walk and the witness. The order matters, and the selection-first
statement it would otherwise be confused with is false. Take a
triangle-free source with \(δ ≤ 1\): Theorem 5.7 then gives every member
of \(Φ_S\) exactly \(k = \max(1, δ) = 1\) connector, so whichever
endpoint is fixed before \(z_0\) is seen has a single target, and
setting \(z_0\) to that target defeats the choice.

\hypertarget{the-one-pass-assembly}{%
\section{10. The one-pass assembly}\label{the-one-pass-assembly}}

We now prove the main theorem. The pieces are in place: Section 4
reduces to an indecomposable, non-bipartite sequence that is not the
\(K_3\) base; Section 8 chooses a pivot; Section 7 walks the quotient;
Section 9 crosses between fibers.

\textbf{Theorem 10.1.} \emph{Every realization graph is
Hamilton-laceable when bipartite on more than one vertex, and
Hamilton-connected otherwise.}

\textbf{Corollary 10.2.} \emph{Every realization graph is homogeneously
traceable: for every realization \(G\) there is a Hamilton path of
\(G(d)\) beginning at \(G\).}

\emph{Proof.} If \(G(d)\) has one vertex, that vertex is itself a
Hamilton path beginning at \(G\). Otherwise \(G(d)\) is connected and
has order at least two. If it is not bipartite it is Hamilton-connected,
so take any \(H ≠ G\) and a Hamilton \(G\)--\(H\) path. If it is
bipartite, connectedness and order at least two make its 2-coloring
proper and surjective, so the class opposite \(G\) is nonempty; take
\(H\) in it and a Hamilton \(G\)--\(H\) path, which laceability
supplies. \(\blacksquare\)

\hypertarget{the-induction-2}{%
\subsection{10.1 The induction}\label{the-induction-2}}

We induct on the number of ground vertices. Let \(d\) be graphical and
let \(G_0 ≠ G_1\) be realizations. One of four branches applies; this is
the dispatch of Section 4.3, carried out here with the conclusion each
branch delivers. The fourth is the work, and it is restated as an
explicit construction in Section 10.2 and run on a worked example at the
end of the section.

If \(G(d)\) is bipartite, Theorem 3.1 gives the conclusion by
classification.

If \(d\) is Tyshkevich-decomposable, Corollary 4.3 propagates the
conclusion from the factors, each of which has fewer ground vertices.

If \(d\) has exactly three realizations then, since \(G(d)\) is not
bipartite, no two of them can be non-adjacent --- coloring a
non-adjacent pair alike and the third oppositely would properly 2-color
\(G(d)\) --- so \(G(d)\) is \(K_3\) and is Hamilton-connected directly.
The branch is stated by the count because the count is what the
induction can test.

Otherwise \(d\) is indecomposable and non-bipartite, and \(G(d)\) is not
\(K_3\). This is the case the rest of the section treats. Activity is
not assumed: under these hypotheses it follows, since the only
indecomposable graph with an inactive vertex is the one-vertex graph and
its realization graph is bipartite.

\textbf{The measure, in one place.} Two of the four branches recurse,
and each recurses on a strictly smaller ground:

\begin{longtable}[]{@{}
  >{\raggedright\arraybackslash}p{(\columnwidth - 4\tabcolsep) * \real{0.3333}}
  >{\raggedright\arraybackslash}p{(\columnwidth - 4\tabcolsep) * \real{0.3333}}
  >{\raggedright\arraybackslash}p{(\columnwidth - 4\tabcolsep) * \real{0.3333}}@{}}
\toprule\noalign{}
\begin{minipage}[b]{\linewidth}\raggedright
branch
\end{minipage} & \begin{minipage}[b]{\linewidth}\raggedright
what it recurses on
\end{minipage} & \begin{minipage}[b]{\linewidth}\raggedright
supplied by
\end{minipage} \\
\midrule\noalign{}
\endhead
\bottomrule\noalign{}
\endlastfoot
Tyshkevich factor & \(d_1\) and \(d_2\), each with fewer & Section
4.2 \\
main case & each fiber \(Φ_S ≅ G(f_S)\), one fewer ground vertex & Lemma
5.1 \\
\end{longtable}

The other two terminate outright: the bipartite branch by Theorem 3.1
and the \(K_3\) branch directly. So the number of ground vertices is a
well-founded measure for all four branches at once. Section 4.3 records
the first; the second is easy to miss, because it happens inside the
construction of Section 10.2 rather than in the dispatch, and it is the
reason Lemma 5.1 is stated where it is.

The branches overlap, and Section 3 is not confined to the one named for
it: the decomposable branch reaches its factors through Corollary 4.3,
which consumes the paired two-disjoint-path property that Corollary 3.3
supplies. Every bipartite factor is carried by Section 3 whichever
branch reaches it.

\hypertarget{one-pass-along-the-quotient}{%
\subsection{10.2 One pass along the
quotient}\label{one-pass-along-the-quotient}}

Apply Theorem 8.1 to \(G_0\) and \(G_1\). It supplies a pivot \(v\) such
that the projections \(A = N_{G_0}(v)\) and \(B = N_{G_1}(v)\) are
distinct, some \(v\)-fiber \(Φ_T\) is non-bipartite, and the pair
\(\{A, B\}\) is not the exception of Theorem 7.7 in the quotient family
at \(v\).

By Section 5 that family is shifted and the quotient is its Johnson
graph. Theorem 7.7 therefore gives a Hamilton \(A\)--\(B\) path

\[
A = q_0, q_1, …, q_m = B
\]

of the quotient. Theorem 8.1 supplies a pivot at which \emph{some} fiber
is non-bipartite, and it need not be the only one; fix such a fiber once
and call it the \textbf{buffer}. The quotient path visits every quotient
vertex exactly once, so the buffer we fixed is \(Φ_{q_j}\) for exactly
one \(j\). Each fiber \(Φ_{q_i}\) is a realization graph on a ground
with one fewer vertex --- this is Lemma 5.1, and it is what makes the
induction well founded --- so the induction hypothesis makes it
maximally Hamiltonian, and Theorem 9.1 applies to every crossing.
Theorem 9.2 needs more --- a non-bipartite target fiber --- and is used
at only the one or two crossings whose target is the buffer, which is
exactly where the construction below invokes it.

Now deliver two distinct entries into the buffer \(Φ_{q_j}\). Call them
\(x\) and \(y\); the quotient vertices are \(q_0, …, q_m\) and nothing
below reuses those letters for a realization.

If \(0 < j < m\), propagate from \(G_0\) along \(q_0, …, q_{j-1}\) by
Theorem 9.1, arriving at some \(x ∈ Φ_{q_j}\); then propagate from
\(G_1\) backwards along \(q_m, …, q_{j+1}\), using Theorem 9.2 on the
last hop with the prescribed target \(x\), arriving at some \(y ≠ x\).

The two end positions are not a degenerate remark, because the
concatenation has a different shape in each and the endpoint claim of
Section 10.3 has to hold in all three. If \(j = 0\) the buffer is the
first fiber, so it already contains \(G_0\): there is \textbf{no left
chain}, \(x = G_0\), and only the right-hand chain is built, its final
hop using Theorem 9.2 with the prescribed target \(G_0\) so that
\(y ≠ G_0\). The concatenation is then the buffer walk from \(G_0\) to
\(y\), followed by the reversed right-hand chain, and its first entry is
\(G_0\) because the buffer walk starts there rather than because a left
chain delivered it. If \(j = m\) the mirror holds, and it needs the same
care: no right chain, \(y = G_1\), and the left-hand chain's final hop
uses Theorem 9.2 with the prescribed target \(G_1\), without which
ordinary propagation could enter the buffer at \(G_1\) and leave the two
entries equal. The buffer walk then \emph{ends} at \(G_1\). In all three
positions the two entries into the buffer are distinct, which is what
the spanning walk needs.

\textbf{Why the buffer must be non-bipartite, which is the reason
Sections 6 and 8 exist.} At every fiber but one the entry is forced
while the exit may still be chosen, and Theorem 9.1's bipartite case
chooses it in the color class opposite the entry, so a laceable fiber
suffices. At the buffer both ends are forced: \(x\) by the chain
arriving from \(G_0\) and \(y\) by the chain arriving from \(G_1\).
Theorem 9.2 delivers \(y ≠ x\) and nothing about their colors. A
bipartite fiber is only Hamilton-laceable, so a same-colored \(x\) and
\(y\) would have no spanning walk at all, and the construction would
fail at the one fiber it cannot reroute. A non-bipartite fiber is
Hamilton-connected, which is indifferent to color. That is the whole of
it, and it is why Theorem 8.1 must deliver a pivot with a non-bipartite
fiber rather than merely a separating one.

\textbf{What the second prescribed endpoint costs.} Exactly this one
fiber. Were only \(G_0\) prescribed there would be no chain arriving
from \(G_1\), every fiber would have its entry forced and its exit still
open, and Theorem 9.1 alone would carry the pass. What the second
endpoint does not pay for is the rest of Sections 6 to 9: the induction
hypothesis must be maximal Hamiltonicity whatever the outermost
statement asks for, since every fiber below the top is entered and left
at vertices fixed from outside it. The stronger conclusion is bought for
one fiber's worth of care, not for a section.

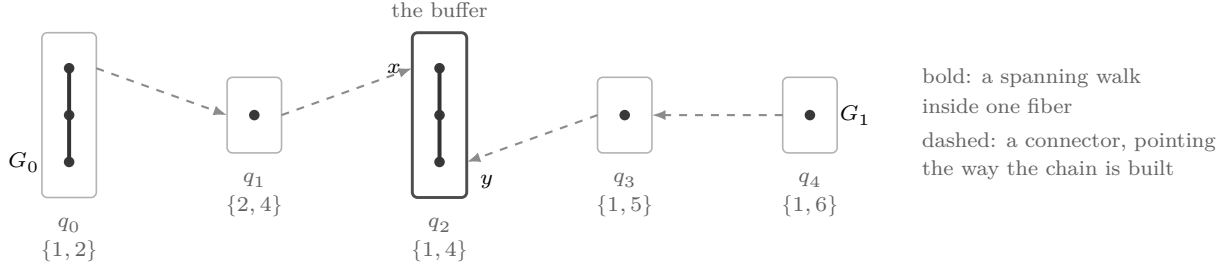
\begin{figure}[tbp]
\centering
\begin{tikzpicture}[x=1cm,y=1cm,
  dot/.style={circle,fill=black!78,inner sep=1.4pt},
  run/.style={draw=black!80,line width=1.5pt},
  conn/.style={draw=black!45,line width=0.8pt,dashed},
  box/.style={draw=black!30,rounded corners=2pt,line width=0.6pt},
  buf/.style={draw=black!70,rounded corners=2pt,line width=1.1pt},
  ann/.style={font=\scriptsize,text=black!62},
  lab/.style={font=\scriptsize}]
\node[dot] at (0,-0.62) {};
\node[dot] at (0,0) {};
\node[dot] at (0,0.62) {};
\draw[run] (0,-0.62) -- (0,0);
\draw[run] (0,0) -- (0,0.62);
\draw[box] (-0.36,-1.09) rectangle (0.36,1.09);
\node[ann] at (0,-1.45) {$q_{0}$};
\node[ann] at (0,-1.79) {$\{1,2\}$};
\node[dot] at (2.45,0) {};
\draw[box] (2.09,-0.5) rectangle (2.81,0.5);
\node[ann] at (2.45,-0.86) {$q_{1}$};
\node[ann] at (2.45,-1.2) {$\{2,4\}$};
\node[dot] at (4.9,-0.62) {};
\node[dot] at (4.9,0) {};
\node[dot] at (4.9,0.62) {};
\draw[run] (4.9,-0.62) -- (4.9,0);
\draw[run] (4.9,0) -- (4.9,0.62);
\draw[buf] (4.54,-1.09) rectangle (5.26,1.09);
\node[ann] at (4.9,-1.45) {$q_{2}$};
\node[ann] at (4.9,-1.79) {$\{1,4\}$};
\node[dot] at (7.35,0) {};
\draw[box] (6.99,-0.5) rectangle (7.71,0.5);
\node[ann] at (7.35,-0.86) {$q_{3}$};
\node[ann] at (7.35,-1.2) {$\{1,5\}$};
\node[dot] at (9.8,0) {};
\draw[box] (9.44,-0.5) rectangle (10.16,0.5);
\node[ann] at (9.8,-0.86) {$q_{4}$};
\node[ann] at (9.8,-1.2) {$\{1,6\}$};
\draw[conn,-{Latex[length=1.8mm]}] (0.36,0.62) -- (2.09,0);
\draw[conn,-{Latex[length=1.8mm]}] (2.81,0) -- (4.54,0.62);
\draw[conn,-{Latex[length=1.8mm]}] (6.99,0) -- (5.26,-0.62);
\draw[conn,-{Latex[length=1.8mm]}] (9.44,0) -- (7.71,0);
\node[lab] at (-0.6,-0.62) {$G_0$};
\node[lab] at (10.4,0) {$G_1$};
\node[lab] at (4.3,0.62) {$x$};
\node[lab] at (5.52,-0.88) {$y$};
\node[ann] at (4.9,1.4) {the buffer};
\node[ann,anchor=west] at (11.15,0.5) {bold: a spanning walk};
\node[ann,anchor=west] at (11.15,0.12) {inside one fiber};
\node[ann,anchor=west] at (11.15,-0.35) {dashed: a connector, pointing};
\node[ann,anchor=west] at (11.15,-0.72) {the way the chain is built};
\end{tikzpicture}
\caption{The concatenated path, on the instance of the worked example closing this section: $d=(4,2,2,2,1,1)$ at the pivot $v=3$. Each box is a fiber over one vertex of the quotient Hamilton path $q_0,\dots,q_4$, drawn with its realizations and the spanning walk through them; the dashed edges are the connectors that cross between consecutive fibers. Reading left to right is not how the path is built. Two chains are built inwards -- one from $G_0$ by Theorem~9.1, one back from $G_1$ -- and they enter the buffer $q_2$ at the two distinct realizations $x$ and $y$, the second of which Theorem~9.2 is invoked to keep distinct from the first. The buffer is the only fiber whose two ends are both forced, and it is the only one required to be non-bipartite: a laceable fiber would have no walk at all between an $x$ and a $y$ of the same color. Splicing the buffer walk between the two chains, with the right-hand one reversed, gives the Hamilton $G_0$--$G_1$ path listed at the end of this section.}
\label{fig:assembly}
\end{figure}

In the worked example closing this section the buffer is a triangle and
the two chains reach it from opposite sides, and the reader may find it
useful to have that instance alongside what follows. Figure
\ref{fig:assembly} draws the finished concatenation on it.

The buffer is a realization graph on a smaller ground and is
non-bipartite, so the induction hypothesis gives a spanning \(x\)--\(y\)
walk inside \(Φ_{q_j}\). Reverse the right-hand chain and concatenate:
the left fiber walks, their connectors, the buffer walk, then the right
connectors and the reversed right fiber walks.

\hypertarget{the-concatenation-is-a-hamilton-path}{%
\subsection{10.3 The concatenation is a Hamilton
path}\label{the-concatenation-is-a-hamilton-path}}

Every consecutive pair in the concatenated list is adjacent in \(G(d)\)
for one of two reasons. It is an internal step of a spanning fiber walk,
which is a 2-switch inside a fiber; or it is the explicit switch
\(\{vb, ac\} → \{va, bc\}\) attached to a quotient hop. Reversing a walk
or a connector preserves adjacency, because \(G(d)\) is undirected.

The list has no repeated realization and covers every realization. The
quotient path visits each quotient vertex once, distinct fibers are
disjoint, and each selected fiber walk contains its whole fiber exactly
once. Its first and last entries are \(G_0\) and \(G_1\).

So the list is a Hamilton \(G_0\)--\(G_1\) path, checked step by step
rather than inferred from a search. This discharges the main case. With
the Tyshkevich reduction above, to strictly smaller grounds, and the two
branches that terminate outright --- Theorem 3.1 for the bipartite case
and the direct argument for \(K_3\) --- the induction closes and Theorem
10.1 follows. \(\blacksquare\)

\hypertarget{a-worked-example}{%
\paragraph{A worked example}\label{a-worked-example}}

\textbf{Why not the running example.} Section 5's \(d = (3,2,2,2,1)\)
cannot close the paper. Keeping that section's ground \(\{1,…,5\}\), its
pivots carrying a non-bipartite fiber are \(1\) and \(5\) --- the
vertices of degree three and one --- while its pivots whose quotient
family is a Y-family are the three of degree two, \(2\), \(3\) and
\(4\), and those two sets are disjoint, so every pivot Theorem 6.1 could
offer there is free of the exception already and the choice Section 8
exists to make never arises. We change sequences for that reason and for
no other.

Take \(d = (4,2,2,2,1,1)\) on the ground \(\{1,…,6\}\), and write a
realization by its edge set, an edge \(uv\) as \(uv\). It is graphical
and Tyshkevich-indecomposable, \(G(d)\) is not bipartite, and \(G(d)\)
has nine realizations rather than three, so all three of Section 10.1's
other branches are refused and the induction is in the fourth. Prescribe

\[
G_0 = \{12,13,14,15,23,46\}, \quad \quad G_1 = \{12,13,14,15,24,36\}.
\]

\emph{The pivot, and why Section 8 is needed.} Four ground vertices
separate the pair and carry a non-bipartite fiber, so Theorem 6.1 has
four admissible pivots to offer: \(2, 3, 4\) and \(6\). If it offers
\(v = 2\), that pivot is \textbf{exceptional} --- its quotient family is
a Y-family and \(\{N_{G_0}(2), N_{G_1}(2)\}\) is exactly its universal
pair --- so by Theorem 7.7 the quotient has \textbf{no} \(A\)--\(B\)
Hamilton path at all, not merely none that the theorem supplies. This is
the failure Section 8 exists to avoid. Theorem 8.3 says the exceptional
pivot's degree class is \(\{2,3,4\}\), and the triangle count in the
proof of Theorem 8.1 says the three cannot all be exceptional; here
neither \(3\) nor \(4\) is. Take \(v = 3\). This sequence is one of the
two smallest \textbf{satisfying every hypothesis of Theorem 6.1} at
which the exception is live, which is why Section 8 names it. The
qualifier is doing work: the Y-family shape appears one ground vertex
earlier, on the five-vertex sequence of Section 5.5, whose three pivots
all carry Y-families.

\emph{The quotient.} The realizable neighborhoods of \(v = 3\) are

\[
I_3 = \{ \{1,2\}, \{1,4\}, \{1,5\}, \{1,6\}, \{2,4\} \},
\]

a down-set of the Gale order on the ground \(\{1,2,4,5,6\}\), so shifted
as Theorem 5.3 requires. It is \textbf{not} principal: \(\{1,6\}\) and
\(\{2,4\}\) are incomparable and both maximal, so there is no
Gale-greatest member; this family is the Y-family \(F(2;5,1)\) of
Section 7.5, with universal pair \(\{\{1,2\},\{1,4\}\}\). The
exceptional shape is therefore present at the pivot we use as well; what
makes \(v = 3\) serviceable is that the prescribed pair
\(\{A, B\} = \{\{1,2\},\{1,6\}\}\) is not that universal pair, which is
exactly the hypothesis Theorem 7.7 asks for. Their fibers have sizes
\(3, 3, 1, 1, 1\), and the quotient is the Johnson graph on \(I_3\) ---
five vertices and eight edges, so not complete. The prescribed
projections are \(A = \{1,2\}\) and \(B = \{1,6\}\), and Theorem 7.7
supplies a Hamilton \(A\)--\(B\) path of it:

\[
\{1,2\}, \quad \{2,4\}, \quad \{1,4\}, \quad \{1,5\}, \quad \{1,6\}.
\]

Figure \ref{fig:quotient10} draws it. The two absent edges are the ones
that make the choice of Section 8 do work rather than be free: a
complete quotient would have an \(A\)--\(B\) Hamilton path for every
prescribed pair.

\begin{figure}[tbp]
\centering
\begin{tikzpicture}[x=1cm,y=1cm,
  qdot/.style={circle,fill=black!78,inner sep=0pt,minimum size=2.4mm},
  bufring/.style={circle,draw=black!80,line width=1.1pt,inner sep=0pt,minimum size=5.2mm},
  pathedge/.style={black!85,line width=1.6pt},
  chord/.style={black!45,line width=0.7pt},
  nonedge/.style={black!45,densely dotted,line width=0.8pt},
  qlab/.style={font=\scriptsize,align=center,anchor=south},
  ann/.style={font=\scriptsize,text=black!62},
  lab/.style={font=\scriptsize}]
\draw[chord] (0,0) .. controls (1.204,-1.57) and (3.096,-1.57) .. (4.3,0);
\draw[chord] (0,0) .. controls (1.806,-2.19) and (4.644,-2.19) .. (6.45,0);
\draw[chord] (0,0) .. controls (2.408,-2.81) and (6.192,-2.81) .. (8.6,0);
\draw[chord] (4.3,0) .. controls (5.504,-1.57) and (7.396,-1.57) .. (8.6,0);
\draw[nonedge] (2.15,0) .. controls (3.354,2.45) and (5.246,2.45) .. (6.45,0);
\draw[nonedge] (2.15,0) .. controls (3.956,3.2) and (6.794,3.2) .. (8.6,0);
\draw[pathedge] (0,0) -- (2.15,0);
\draw[pathedge] (2.15,0) -- (4.3,0);
\draw[pathedge] (4.3,0) -- (6.45,0);
\draw[pathedge] (6.45,0) -- (8.6,0);
\node[qdot] at (0,0) {};
\node[qlab] at (0,0.34) {$q_{0}$\\$ \{1,2\} $\\[-1pt]{\tiny fiber of 3}};
\node[qdot] at (2.15,0) {};
\node[qlab] at (2.15,0.34) {$q_{1}$\\$ \{2,4\} $\\[-1pt]{\tiny fiber of 1}};
\node[bufring] at (4.3,0) {};
\node[qdot] at (4.3,0) {};
\node[qlab] at (4.3,0.34) {$q_{2}$\\$ \{1,4\} $\\[-1pt]{\tiny fiber of 3}\\[-1pt]{\tiny the buffer}};
\node[qdot] at (6.45,0) {};
\node[qlab] at (6.45,0.34) {$q_{3}$\\$ \{1,5\} $\\[-1pt]{\tiny fiber of 1}};
\node[qdot] at (8.6,0) {};
\node[qlab] at (8.6,0.34) {$q_{4}$\\$ \{1,6\} $\\[-1pt]{\tiny fiber of 1}};
\node[lab] at (-0.78,0) {$A$};
\node[lab] at (9.38,0) {$B$};
\node[ann,anchor=north] at (4.3,-3.36) {heavy: the Hamilton $A$--$B$ path of Theorem 7.7\quad light: the quotient's other edges};
\node[ann,anchor=north] at (4.3,-3.76) {dotted: the two non-adjacent pairs, so the quotient is not complete};
\end{tikzpicture}
\caption{The quotient of the worked example, on $d=(4,2,2,2,1,1)$ at the pivot $v=3$: the Johnson graph of $I_{3}$, one vertex per fiber. It has 5 vertices and 8 edges, two short of complete -- $\{2,4\}$ is adjacent to neither $\{1,5\}$ nor $\{1,6\}$, since those pairs differ in two elements rather than one. The prescribed projections $A$ and $B$ are the two ends, and the heavy edges are the Hamilton $A$--$B$ path Theorem 7.7 supplies; the family is a Y-family, so that path exists only because $\{A,B\}$ is not its universal pair, which is the choice Section 8 makes. The ringed vertex $q_{2}$ is the buffer: its fiber is a triangle, hence not bipartite, and it is the one fiber entered from both sides. Figure~\ref{fig:assembly} draws what happens inside these five fibers.}
\label{fig:quotient10}
\end{figure}

\emph{The buffer.} Two fibers are non-bipartite, those over \(\{1,2\}\)
and \(\{1,4\}\), each a triangle. Fix the one over \(\{1,4\}\); it is
\(q_2\), interior to the path, so the two-chain case applies.
Propagating from \(G_0\) across \(q_0\) and \(q_1\) by Theorem 9.1
enters it at \(x = \{12,13,14,15,26,34\}\); propagating back from
\(G_1\) across \(q_4\) and \(q_3\), with Theorem 9.2 on the last hop
excluding the prescribed \(x\), enters it at
\(y = \{12,13,15,16,24,34\}\). The buffer is Hamilton-connected because
it is not bipartite, which is what makes a walk from \(x\) to \(y\)
available whatever colors they would have carried, and the walk is
\(x\), \(\{12,13,14,16,25,34\}\), \(y\).

\emph{The concatenation.} Reversing the right-hand chain and splicing
gives

\[
\begin{gathered}
\{12,13,14,15,23,46\}, \quad \quad \{13,14,15,16,23,24\}, \quad \quad \{12,13,14,16,23,45\}, \\
\{12,14,15,16,23,34\}, \\
\{12,13,14,15,26,34\}, \quad \quad \{12,13,14,16,25,34\}, \quad \quad \{12,13,15,16,24,34\}, \\
\{12,13,14,16,24,35\}, \\
\{12,13,14,15,24,36\}.
\end{gathered}
\]

Nine entries, no repetition, consecutive entries differing in four
edges, first \(G_0\) and last \(G_1\): a Hamilton \(G_0\)--\(G_1\) path
of \(G(d)\), which is the test of Section 10.3 applied to this instance.
The construction was run and checked by machine, and the script refuses
to emit unless every claim above holds on it
(\texttt{compute/sec10\_walkthrough\_build.py}); the choice of sequence
was itself computed, as one of the two smallest instances satisfying
every hypothesis of Theorem 6.1 whose Section 8 exception is live
(\texttt{compute/sec10\_walkthrough\_instance\_sweep.py}).

\hypertarget{what-the-proof-rests-on}{%
\subsection{10.4 What the proof rests
on}\label{what-the-proof-rests-on}}

The argument above rests on results from outside this paper, and they
divide into three kinds. \textbf{Seven are taken on trust}, and they are
listed below because the formal development enters exactly these as
assumptions and nothing else: naming them is what makes Section 11.1's
claim checkable. \textbf{Five more are external but proved rather than
assumed} in an accompanying formalization --- \textbf{C1} to \textbf{C3}
in Section 3.2, \textbf{C4} in Section 4.2, \textbf{C5} in Section 6.4.
\textbf{Two classical facts about split graphs are proved here rather
than quoted.}

\textbf{The seven taken on trust.} Two are directions of Barrus's
classification of triangle-free realization graphs {[}5{]}, which
supports Section 3 and the fiber criterion of Section 6: that a
realization graph is bipartite exactly when it is triangle-free, and
that a bipartite one is a Cartesian product of transposition graphs and
at most one copy of \(K_{6,6} − 6K_2\). Connectivity of the realization
graph is a consequence of a theorem of Fulkerson, Hoffman and McAndrew
{[}16{]}. Section 8 uses the sufficiency half of the Erdős--Gallai
criterion {[}13{]}, in Lemma 8.3e, to turn the non-realizability of a
corner into a violating rank; the necessity half is proved rather than
assumed, inside Step 3 of Lemma 8.3d. Section 6 uses Ryser's theorem on
invariant positions in a class of zero-one matrices with fixed margins,
in the form given by Brualdi {[}9, Theorem 3.4.1{]}: an invariant
position forces every matrix of the class to decompose into an all-ones
block and an all-zeros block, and that is what carries Section 6.4 from
indecomposability of the degree sequence to primeness of the interchange
class. The last two arrive with the companion paper's buffer lemma
rather than on their own account, and we list them separately because
the formal trace charges them separately: Ryser's theorem that any two
matrices of a nonempty class with fixed margins are joined by a sequence
of interchanges {[}31, Theorem 3.1{]}, and Brualdi's theorem that a
non-bipartite interchange graph with no invariant positions contains a
triangle {[}9, Theorem 6.3.4{]}.

\textbf{The five proved elsewhere.} \textbf{C1} to \textbf{C3} are
stated in Section 3.2, \textbf{C4} in Section 4.2 and \textbf{C5} in
Section 6.4. \textbf{C1} and \textbf{C2} are Theorem 1.5 and Proposition
1.1(c) of Coleman, Fischberg, Gong, Harrington and Wong {[}10{]}, and
carry the crown factor through Section 3. \textbf{C3}, \textbf{C4} and
\textbf{C5} are the companion paper's {[}4{]}: the paired
two-disjoint-path coverability of products of complete transposition
graphs and the Cartesian-product lift, used in Sections 3 and 4, and its
buffer lemma, used in Section 6.4. The one case \textbf{C1} does not
reach, a factor \(CT_2\) of rank two, is proved in Appendix A rather
than imported. None of the five is charged by the formal trace: the
companion's development proves its own three, and it proves \textbf{C1}
and \textbf{C2} from the foundations rather than assuming them.

\textbf{The two proved here.} Section 6.4 uses the forbidden-subgraph
characterization of split graphs, in the form that a graph with no
induced \(2K_2\) and no induced \(C_4\) is split or contains an induced
\(C_5\), and the fact that every realization of a split degree sequence
is split at the same partition. The first is due to Földes and Hammer
{[}15{]}. Both are classical, both are proved in the formal development,
and neither is on the list of seven.

Section 4 cites Tyshkevich's decomposition theorem {[}34{]}, and the
formal development does not need it: the induction asks only whether
some nontrivial composition exists and applies Theorem 4.2 to whichever
one it is handed, so neither canonicity nor uniqueness enters. The
citation stays because it explains why the decomposition is a property
of the sequence rather than of a realization.

One route is no longer needed. Earlier approaches reached the quotient
through matroid basis graphs and the theorem of Naddef and Pulleyblank
{[}25,26{]}; Theorem 7.7 covers every shifted family and Section 7.9
sets out why the polytope route cannot be repaired rather than merely
avoided. The axiom trace bears it out: the machine-checked main theorem
does not charge the Naddef--Pulleyblank exchange property.

\hypertarget{mechanization-and-computation}{%
\section{11. Mechanization and
computation}\label{mechanization-and-computation}}

Computation was how this theorem was found. It is not how it is
certified, and this section is ordered to say so: what the kernel checks
comes first, and the computational record that follows reports only what
no theorem here covers. This is the accounting, not the reproduction
manual. The development itself documents the build, the module sizes,
the figures' provenance, and every divergence between the formal
statements and the printed ones.

\hypertarget{what-is-machine-checked}{%
\subsection{11.1 What is
machine-checked}\label{what-is-machine-checked}}

\textbf{Theorem 10.1 is machine-checked.} The main theorem is proved in
Lean, unconditionally, with no unproved step anywhere beneath it and
nothing outside the kernel: the audited surface carries no
\texttt{sorry} and no compiled-evaluation certificate. Its axiom trace
is the ambient foundations \texttt{propext}, \texttt{Classical.choice}
and \texttt{Quot.sound}, together with seven cited results, and nothing
else.

\textbf{The three kinds of outside result are set out in Section 10.4},
with their locators and the sections that use them. That list is the
authoritative one and is not repeated here. What belongs here is the
correspondence with the kernel: the \textbf{seven taken on trust} are
exactly what the formal development charges as assumptions, so the axiom
trace and Section 10.4's list can be read against each other item by
item. Nothing else is charged --- in particular neither the five results
\textbf{C1} through \textbf{C5} listed in Section 10.4, all of which the
accompanying development proves, nor the two classical facts about split
graphs behind Section 6.4, which are proved there as well. The one case
\textbf{C1} does not reach is proved here as Lemma A.1.

\textbf{Section 7 is machine-checked from the foundations alone}: the
quotient theorem, its crossing lemma, the classification of clique-sums
as Y-families, the duality of Section 7.7 and both corollaries are
proved in a module of 284 declarations that charges no cited axiom of
any kind.

\textbf{Every numbered result of Sections 3 through 10 is
machine-checked, and the argument given for it here is the argument the
kernel checks}, with no exceptions. Where an earlier draft argued a
numbered result differently from the formalization, the text was brought
to the formal route rather than the other way round.

The counts here are produced by a gate rather than by hand: \textbf{115
audited declarations, no \texttt{sorry} among them, and exactly the
seven cited axioms above}. The gate re-derives every one of those traces
from the kernel, and the six conditions it fails on --- an unexpected
cited axiom among them --- are listed with the development.

One convention needs stating, because the formal hypothesis is not word
for word the printed one. The formalization carries Tyshkevich
indecomposability as \emph{no realization of \(d\) decomposes}, where
Section 4 states it through the canonical decomposition;
\textbf{Tyshkevich's theorem identifies the two and we have not
formalized it}. The formal reading makes indecomposability harder to
satisfy, so it cannot let any statement below claim more than it should.

\hypertarget{what-computation-alone-supports}{%
\subsection{11.2 What computation alone
supports}\label{what-computation-alone-supports}}

Five statements in this paper are not theorems and cannot be. Each is a
count over a finite range, each is load-bearing somewhere, and none
follows from anything proved above. They are collected here as the
computational claims the main line leans on. Other counts reported in
passing elsewhere are also computational; their ranges are given where
they are reported, and Section 11.3 describes the separate checks of the
formal statements against enumerated realization graphs.

\textbf{No ground vertex need have a matroid family.} Used in Sections 1
and 5.5, where it is the reason the matroid route cannot be taken at
all. Exhaustive through ground order eleven: every degree sequence on at
most ten ground vertices has a ground vertex whose realizable family is
the basis family of a matroid, and at eleven exactly \textbf{six} do
not. Eleven is therefore the smallest order at which the matroid reading
fails, and the failure is not an artifact of where we stopped looking.
The six are

\[
\begin{gathered}
(8,8,8,4,4,4,4,2,2,2,2) \quad (8,8,8,5,5,4,4,2,2,2,2) \quad (8,8,8,5,5,5,5,2,2,2,2) \\
(8,8,8,8,5,5,5,5,2,2,2) \quad (8,8,8,8,6,6,5,5,2,2,2) \quad (8,8,8,8,6,6,6,6,2,2,2)
\end{gathered}
\]

and any one of them can be checked against the Erdős--Gallai criterion
by hand. A second implementation reproduced both the count and these
six, testing basis exchange on the family directly rather than through
the Gale-greatest-member criterion of Section 5.5, so the reproduction
does not rest on Theorem 5.3.

\textbf{The running example is the smallest of its kind.} Used in
Section 5.1, which introduces \(d = (3,2,2,2,1)\) as the smallest
sequence whose quotient exhibits the exception and says there is no
other of its order. Exhaustive over all graphical sequences through
ground order six, counting those with a pivot whose realizable family is
a Y-family: \textbf{none} at orders four and below, exactly \textbf{one}
at order five, and \textbf{four} at order six.

\textbf{The skeleton of \(\operatorname{conv}(F)\) strictly contains
\(J(F)\).} Used in Section 10.4, where it is the reason the
Naddef--Pulleyblank route cannot be repaired rather than merely avoided.
Exhaustive through ground order six, among the shifted families with at
least three members: \textbf{88 of 172}.

\textbf{\(δ < 0\) occurs.} Used in Sections 5 and 9.1 as a warning
against over-reading Theorem 5.7, whose inequality is asserted only for
a triangle-free source. Exhaustive through ground order seven: negative
values at \textbf{8,876 of 88,518} interfaces.

\textbf{The exception of Theorem 7.7 occurs at the predicted rate}, and
this is the one that does the most work. Two enumerations, deliberately
on different ranges: on the ground \([8]\) alone, at rank two,
\textbf{five} non-Hamiltonian pairs among 13,888; cumulatively over the
grounds \([2]\) through \([8]\), at all ranks \(1 ≤ k < n\), \textbf{70}
clique-sums among 15,477 family instances, every one of them a Y-family.

\textbf{That the exception is inhabited at all is proved, not computed},
and we separate the two because an earlier draft did not: Section 7.5
exhibits the Y-families and the formalization proves the exhibited
families really are shifted, Theorem 7.6 shows each is a clique-sum, and
Lemma 7.4 shows a clique-sum admits no Hamilton path between the members
of its universal pair. What the enumeration adds is the \textbf{rate},
and the rate is what makes the classification tight rather than merely
uncontradicted.

\hypertarget{the-formal-statements-against-the-objects}{%
\subsection{11.3 The formal statements against the
objects}\label{the-formal-statements-against-the-objects}}

A kernel checks a formalization, not a paper. If our definitions do not
capture realization graphs then Lean proves a true statement about the
wrong objects, and no amount of kernel-checking notices. That is the one
gap machine verification cannot close by itself, and the only instrument
against it is to compute with the objects themselves and compare.

So maximal Hamiltonicity was verified directly, on realization graphs
built by enumeration. At ground order seven the run reached 340 of the
342 graphical degree sequences and decided 331 of them, with no failures
and 9 left undecided by a search cap. At ground order eight the
indecomposable census decided 504 of 644 through that order, again with
no failures, skipping 133 by a vertex cap and leaving 7 undecided.
\textbf{Every skipped instance is reported as skipped and none as
passed.} The same reading applies to the certificate of Section 3.3,
where all 3,600 admissible demands on the crown were certified
exhaustively in the kernel.

Its force is limited and worth stating. The enumeration shares the
formalization's adjacency convention, edge sets at Hamming distance
four, so it cannot detect a mis-encoding of that convention --- which is
instead checked internally, by a theorem that a four-edge difference is
an alternating four-cycle. What the enumeration can detect, and a kernel
cannot, is a formal statement true of some other family of objects.

\hypertarget{disclosure}{%
\subsection{11.4 Disclosure}\label{disclosure}}

Where a result was load-bearing, its computation was reproduced by a
second, independently written program rather than by a second run of the
first. The biconditional form of Theorem 7.11 is the case in point: at
ground order seven alone the enumeration finds 409,416 non-exceptional
and 40 exceptional configurations, and cumulatively over grounds one
through seven it finds 427,156 configurations, of which 427,086 are
non-exceptional and 70 exceptional. The second implementation,
enumerating them a different way, reproduced the cumulative figures
exactly, with no violation either time.

The scripts that produced the numbers quoted here, and the Lean
development, are in the public companion repository named under Data
availability, with one exception: the order-seven verification of
Section 11.3 was run on separate hardware and only its summary survives.

The use of AI-assisted tools in this work is described in the
Acknowledgments and in the declaration that follows them.

\hypertarget{what-follows-and-what-does-not}{%
\section{12. What follows, and what does
not}\label{what-follows-and-what-does-not}}

\hypertarget{the-cycle-form}{%
\subsection{12.1 The cycle form}\label{the-cycle-form}}

Brualdi asked his question in the cycle form first. Quoted as the
companion paper quotes it {[}8, Problem 3.7{]}:

\begin{quote}
\textbf{Problem (Brualdi, 1980, Problem 3.7).} \emph{Does the
interchange graph \(G(R,S)\) have a Hamilton cycle? If it has no
Hamilton cycle, does it have a Hamilton path? If it has no Hamilton
path, is there a matrix \(A ∈ A(R,S)\) and a family of paths from \(A\)
covering every other matrix exactly once?}
\end{quote}

Here \(A(R,S)\) is the class of zero-one matrices with row and column
sums \(R\) and \(S\), and \(G(R,S)\) its interchange graph; by the
correspondence of Arikati and Peled {[}3{]} that is the realization
graph of a split degree sequence. Mütze's survey is less committal,
asking of the realization graph \emph{``whether \(G(d)\) admits a
Hamilton path or cycle''} {[}24, Problem P59{]}. Theorem 1.1 answers the
path form in its strongest version. The cycle form follows, and we
record it because nothing above says it.

\textbf{Corollary 12.1.} \emph{The realization graph of a graphical
degree sequence has a Hamilton cycle whenever it has at least three
vertices. The only realization graphs without one are \(K_1\) and
\(K_2\).}

\emph{Proof.} \(G(d)\) is connected, so it has an edge \(xy\). If
\(G(d)\) is non-bipartite it is Hamilton-connected by Theorem 1.1; if it
is bipartite it is Hamilton-laceable, and \(x\) and \(y\) lie in
opposite color classes because they are adjacent. Either way there is a
Hamilton \(x\)--\(y\) path, and closing it with the edge \(xy\) gives a
cycle through every vertex, of length at least three by hypothesis.
\(K_1\) and \(K_2\) have no cycle at all. \(\blacksquare\)

\hypertarget{questions}{%
\subsection{12.2 Questions}\label{questions}}

\textbf{The exception, and whether it must be dodged.} The exception
here is Theorem 7.7's --- a Y-family \(F(k;r,s)\) at its universal pair
--- and not the crown of Section 3, which is unrelated. Section 8 exists
to choose a pivot that does not present it, and it could in principle be
crossed instead. In a Y-family the two arms are cliques with no edge
between them while the universal pair is adjacent to everything, so that
pair is a two-element cut, and Lemma 7.4 is the observation that
deleting the endpoints of a Hamilton path between them would leave one
path covering a disconnected set. A walk that repeats a vertex does
exist, out through one arm and back through the other; but every
crossing of the cut uses an endpoint, so the assembly would have to
traverse one \textbf{endpoint} fiber twice, covering it by two disjoint
paths on four terminals, one of them the pinned realization. That is the
companion's doubled-layer device asked of a non-bipartite fiber, and
Hamilton-connectedness does not supply it --- the same gap Section 2
records in the bipartite case, where laceability is not enough and
Corollary 3.3 supplies the paired property instead.

\textbf{The obvious strengthening is false, which is worth knowing
before anyone attempts it.} A non-bipartite realization graph need not
admit such a pair for \emph{every} four prescribed terminals. Take
\(d = (2,2,2,1,1)\): its realization graph is \(K_{3,3}\) together with
one vertex \(u\) adjacent to all six, and the demand fails exactly when
the four terminals are one color class of the \(K_{3,3}\) together with
\(u\). Call that class \(A\) and the other \(B\). The obstruction is a
count --- a path joining two \(A\)-vertices without using \(u\)
alternates and carries one more \(A\) than \(B\), leaving the second
path two \(B\)-vertices to join and no edge between them. What is open
is the flexible form, in which only the pinned realization is prescribed
and the construction draws the other three terminals from menus. It
survives every instance we have reached: over the non-bipartite
realization graphs of order five to thirty on grounds of at most eight
vertices, \textbf{462 instances attempted}: 459 passed, none failed, and
three exhausted a per-instance time budget. That is against 22 failures
of the rigid property on the smaller range where both were measured. A
further 717 classes exceeded the order bound; those and the three budget
skips are reported as untested rather than as passes. A proof would
remove Section 8 at the cost of a stronger Section 9.

\textbf{Which non-shifted families?} Theorem 7.7 classifies the failures
of Hamilton-connectedness among shifted families. Section 7.9 records
that the corresponding decision problem for an arbitrary family of
\(k\)-subsets is NP-hard {[}23{]}, so no classification can cover all
families; the question is which structural hypotheses weaker than
shiftedness still admit one. Shiftedness is used in Section 7 through
Lemma 7.1's decrement counts and through the layer split of Section 7.8,
and both are more specific than the theorem seems to need.

\textbf{Pancyclicity.} Hamilton-connectedness gives cycles through every
vertex; it says nothing about cycles of intermediate length. For the
quotient in the principal case something is known --- Alspach and Liu's
edge-pancyclicity {[}1{]}, on the simple matroids among them --- but for
the realization graph itself nothing is. Is every non-bipartite
realization graph on \(N\) vertices pancyclic --- does it contain a
cycle of every length from three to \(N\) --- and is every bipartite one
bipancyclic in the corresponding even-length sense? An exhaustive census
finds no counterexample to either. At ground order eight it covers 1,037
classes and 602,320 realizations, reaching realization graphs of 19,355
vertices, every non-bipartite one pancyclic and every bipartite one
bipancyclic; 365 were settled by exact search and the rest by a
one-sided finder whose cycles were revalidated against an independently
rebuilt adjacency, none left unconfirmed. That computation belongs to a
companion project and is not among the ones Section 11 reports for this
paper. Nothing here bears on the question. The induction produces one
spanning walk at a time, and a shorter cycle would have to omit whole
fibers, which the quotient walk gives no way to do. This is the natural
strengthening of Theorem 1.1 along an axis the present argument does not
see.

\textbf{An algorithm.} A Hamilton path in a flip graph is a Gray code: a
listing of every object in which consecutive entries differ by one local
move. Ruskey {[}29{]} poses the companion paper's question in exactly
that setting. Exercise 53 of his Section 5.14, rated a research problem,
defines \(G(R,C)\) on the \((0,1)\)-matrices with prescribed margins,
adjacent when they differ by a \(2 × 2\) submatrix and its bit-flip ---
our 2-switch --- records that Ryser proved it connected, and asks the
reader to investigate its Hamiltonicity. The exercise after it asks why
\(G(R,C)\) is Hamiltonian in the regular case and why a Hamilton cycle
then yields a constant amortized time generation algorithm, so the link
between traversal and generation is his rather than ours. Exercise 20 of
his Section 4.13, also rated a research problem, asks for such an
algorithm for all \((0,1)\)-matrices with prescribed margins, and we
know of no solution.

Every step of the induction is an explicit construction: a Hamilton path
of the quotient, a spanning walk inside each fiber, and a connector
between consecutive fibers. So the proof should yield a procedure that
outputs a Hamilton path between two prescribed realizations. We have not
analyzed its running time, and the standard a Gray code is held to ---
amortized constant time per step, with the objects generated rather than
stored --- is a further question that the recursion through fibers does
not obviously meet. The comparison to beat is Merino and Mütze {[}22{]},
who traverse the skeleton of any 0/1-polytope with a delay only a
logarithmic factor above the cost of one linear optimization over its
vertices. That does not apply here, for the reason Section 7.9 gives:
off the matroid case the skeleton of \(\operatorname{conv}(F)\) strictly
contains \(J(F)\), so a skeleton traversal takes steps that are not
2-switches.

\hypertarget{acknowledgments}{%
\section{Acknowledgments}\label{acknowledgments}}

Petr Hladík and Jiří Fink responded to the point raised in Section 1.3
about their Lemma 3.2, and gave a sharper diagnosis of it than the one
they were sent: the symmetric difference of two realizations decomposes
into edge-disjoint alternating cycles when the underlying graph is
bipartite, but in general it may decompose into two odd cycles sharing
the vertex whose neighborhood is at issue. That account is theirs, not
ours.

The computational search and verification programs and the Lean 4
formalization of Section 11 were developed with AI-assisted tools under
author direction: Anthropic Claude (Opus 5 and Fable 5.1, and earlier
models of both lines) and OpenAI Codex (GPT-6 Astra and GPT-5.6 Sol, and
earlier models). All mathematical content was verified independently of
any model's assertions. The argument is machine-checked in Lean 4, as
Section 11.1 sets out, and the finite claims are the ones collected in
Section 11.2 with their ranges. The author takes full responsibility for
the content.

\hypertarget{funding}{%
\section{Funding}\label{funding}}

This research did not receive any specific grant from funding agencies
in the public, commercial, or not-for-profit sectors.

\hypertarget{data-availability}{%
\section{Data availability}\label{data-availability}}

No datasets were generated or analyzed for this article. The Lean 4
formalization and the verification artifacts supporting its findings are
openly available in the public companion repository at
\url{https://github.com/jbaggett/realization_graphs_lean}.

\hypertarget{declaration-of-generative-ai-and-ai-assisted-technologies-in-the-writing-process}{%
\section{Declaration of generative AI and AI-assisted technologies in
the writing
process}\label{declaration-of-generative-ai-and-ai-assisted-technologies-in-the-writing-process}}

During the preparation of this work the author used Anthropic Claude and
OpenAI Codex in order to draft and edit prose. After using these tools,
the author reviewed and edited the content as needed and takes full
responsibility for the content of the publication. The role of these
tools in developing the mathematics itself --- proof search, the
computational programs, and the Lean 4 formalization --- is described in
the Acknowledgments.

\hypertarget{appendix-a.-the-prism-lemma}{%
\section{Appendix A. The prism
lemma}\label{appendix-a.-the-prism-lemma}}

Section 3 factors a bipartite realization graph into complete
transposition graphs and at most one crown, and lifts paired
two-disjoint-path coverability through the product. When no crown
occurs, \textbf{C3} covers every case, the pure hypercube included. When
a crown occurs, the welding theorem \textbf{C1} carries any factor of
rank at least three, but not a factor \(CT_2\). Multiplying by
\(CT_2 = K_2\) is taking the prism over the rest of the product, and
this appendix supplies that one case.

The statement is about bipartite graphs and paired two-disjoint-path
coverability alone; no realization graph appears in it. It is proved in
the Lean development accompanying this paper, and its axiom trace is the
three ambient foundations and nothing else.

\hypertarget{a.1-statement}{%
\subsection{A.1 Statement}\label{a.1-statement}}

Throughout, \(H\) is balanced bipartite with coloring \(χ_H\), and
\emph{paired two-disjoint-path coverable} has the meaning fixed in
Section 2.

Write \(K_2 □ H\) for the prism over \(H\), with layers \(H^0\) and
\(H^1\), and give it the coloring

\[
χ(u^i) = χ_H(u) ⊕ i.
\]

Call the edge \(u^0u^1\) the \textbf{rung} at \(u\). The prism is again
balanced bipartite under \(χ\).

\textbf{Lemma A.1 (prism).} \emph{If \(H\) is balanced bipartite,
\(|V(H)| ≥ 4\), and \(H\) is paired two-disjoint-path coverable, then so
is \(K_2 □ H\).}

\hypertarget{a.2-proof}{%
\subsection{A.2 Proof}\label{a.2-proof}}

We first record that such an \(H\) is Hamilton-laceable, since the
argument below uses it repeatedly. Let \(x\) and \(y\) be oppositely
colored. As \(|V(H)| ≥ 4\) and \(H\) is balanced, some opposite-colored
pair \((a, b)\) avoids them; apply coverability to \((x, y)\) and
\((a, b)\). The second path joins two distinct vertices, so it has an
edge; name its endpoints \(u\) and \(v\) with \(χ_H(v) = χ_H(y)\). The
two paths are disjoint, so \(u\) and \(v\) differ from \(x\) and \(y\),
and the demands \((x, v)\) and \((u, y)\) are each opposite-colored and
involve four distinct vertices. Apply coverability again to those,
obtaining \(Q_1\) from \(x\) to \(v\) and \(Q_2\) from \(u\) to \(y\),
disjoint and spanning. Then \(Q_1 + vu + Q_2\) is a Hamilton
\(x\)--\(y\) path.

Take four distinct terminals of the prism in two opposite-colored
prescribed pairs. Exchanging the layer names is an automorphism, and it
\textbf{complements} the coloring rather than preserving it, since
\(χ(u^{1−i}) = χ_H(u) ⊕ (1−i)\). That is harmless here: a pair is
opposite-colored exactly when its two colors differ, and complementation
preserves that. So we may assume \(H^0\) holds no more terminals than
\(H^1\). Let \(r\) be the number of terminals in \(H^0\), so
\(r ∈ \{0, 1, 2\}\). Projections of terminals in different layers may
coincide, and the argument allows it throughout.

\textbf{Case \(r = 0\).} All four terminals lie in \(H^1\). Each
prescribed pair is opposite-colored in \(χ\), so by the display above it
is opposite-colored in \(χ_H\), and coverability of \(H\) gives two
disjoint paths covering \(H^1\) and joining the prescribed pairs. At
least one of them contains an edge \(xy\), since its endpoints are
distinct and oppositely colored. Replace that edge by

\[
x^1x^0 + R[x^0, y^0] + y^0y^1,
\]

where \(R[x^0, y^0]\) is a Hamilton \(x^0\)--\(y^0\) path of \(H^0\).
The result is still a path: the old one met only \(H^1\), the inserted
interior meets only \(H^0\), and the only layer changes are the two
rungs displayed. It keeps its prescribed endpoints and now covers
\(H^0\) as well. With the untouched second path this is the required
cover.

\textbf{Case \(r = 1\).} Write the pairs as \((u^0, v^1)\) and
\((a^1, b^1)\). Opposition in the prism gives

\[
χ_H(u) = χ_H(v), \quad χ_H(a) ≠ χ_H(b).
\]

Choose \(z\) in the color class opposite \(u\), avoiding whichever of
\(a, b\) lies in that class. Such a \(z\) exists because each color
class of \(H\) has at least two vertices. It differs from \(u\) and
\(v\) by color, and from \(a\) and \(b\) by construction, so
\(z, v, a, b\) are four distinct vertices of \(H\).

Take a Hamilton \(u^0\)--\(z^0\) path of \(H^0\), and apply coverability
of \(H\) in the upper layer to the pairs \((z^1, v^1)\) and
\((a^1, b^1)\). Join the lower Hamilton path to the upper
\(z^1\)--\(v^1\) path through the rung \(z^0z^1\). That produces a
\(u^0\)--\(v^1\) path; the other upper path already joins \(a^1\) to
\(b^1\). The lower path covers \(H^0\), the upper cover covers \(H^1\),
and the two are disjoint apart from the deliberate join at the two
distinct copies of \(z\).

\textbf{Case \(r = 2\).} Let \(c\) be the number of prescribed pairs
crossing between layers and \(m\) the number confined to the lower
layer. A crossing pair contributes one lower terminal and a
lower-confined pair contributes two, so

\[
c + 2m = 2,
\]

and \(c\) is \(0\) or \(2\). The intermediate distribution, exactly one
crossing pair, cannot occur alongside a two-two split between the
layers.

If \(c = 0\), each layer holds one whole prescribed pair, whose
endpoints are opposite in \(χ_H\). A Hamilton path for the lower pair
and a Hamilton path for the upper pair are disjoint and together cover
the prism.

Let \(c = 2\), with pairs \((u^0, v^1)\) and \((a^0, b^1)\), so that

\[
χ_H(u) = χ_H(v), \quad χ_H(a) = χ_H(b).
\]

Suppose first that each color class of \(H\) has at least three
vertices. Choose connectors \(x\) and \(y\) with \(χ_H(x) ≠ χ_H(u)\) and
\(χ_H(y) ≠ χ_H(a)\), as follows. If \(χ_H(u) = χ_H(a)\), take \(x\) and
\(y\) distinct in the other class; they differ from all four of
\(u, a, v, b\) by color, whatever coincidences hold between lower and
upper projections. If \(χ_H(u) ≠ χ_H(a)\), take \(x\) in the class of
\(a\) and \(b\) avoiding \(\{a, b\}\), and \(y\) in the class of \(u\)
and \(v\) avoiding \(\{u, v\}\); each forbidden set has at most two
elements, since a lower projection may equal an upper one, and each
class has at least three, so both choices exist. And \(x ≠ y\) because
they have different colors.

Either way the four demands

\[
\text{ in } H^0: \quad (u, x), (a, y) \quad \text{ in } H^1: \quad (x, v), (y, b)
\]

are opposite-colored and, within each layer, involve four distinct
vertices. Take a cover of each layer and join the lower \(u\)--\(x\)
path to the upper \(x\)--\(v\) path through the rung \(x^0x^1\), and
likewise the two paths at \(y\). The layer covers were spanning and
disjoint and \(x ≠ y\), so the two resulting paths are disjoint, span
the prism, and join the prescribed pairs.

Finally suppose each color class has two vertices, so \(|V(H)| = 4\).
Then \(H = K_{2,2} = C_4\). To see it, take any \(x\) in one class and
\(y\) in the other, prescribe \((x, y)\) and the remaining pair
\((x', y')\), and apply coverability. The two paths together have four
vertices and two edges, and each joins opposite colors, so each has odd
positive length; hence both have length one, and \(xy\) is an edge. As
\(x\) and \(y\) were arbitrary, \(H\) is complete bipartite.

So the prism is \(C_4 □ K_2 = Q_3\). Since \(CT_2 = K_2\) and
\(C_4 = CT_2 □ CT_2\), the cube is a Cartesian product of complete
transposition graphs whose order is eight, and \textbf{C3} applies to it
directly: \(Q_3\) is paired two-disjoint-path coverable. \(Q_3\) is
connected and bipartite, so its proper 2-coloring is unique up to
complementation and complementation preserves opposite-coloredness; the
conclusion therefore holds for the coloring at hand. \(\blacksquare\)

\hypertarget{a.3-two-remarks}{%
\subsection{A.3 Two remarks}\label{a.3-two-remarks}}

\textbf{The opposite-colored hypothesis is not decorative.} Weaken it to
ask only that the four terminals split two-and-two between the color
classes, allowing a prescribed pair to be monochromatic, and the lemma
fails already at the cube. Prescribe \((000, 110)\) and \((100, 010)\),
each pair monochromatic. Two paths covering eight vertices have six
edges between them; each of these pairs has distinct endpoints of the
same color, so each path has positive even length, and the lengths must
be \(2\) and \(4\) in some order. The only common neighbors of \(000\)
and \(110\) are \(100\) and \(010\), which are the other path's
endpoints, so the first path cannot have length two; the only common
neighbors of \(100\) and \(010\) are \(000\) and \(110\), so neither can
the second. So neither Lemma A.1 nor its cube base survives that
weakening.

\textbf{Where it is used, and that the case occurs.} Section 3 only, and
there only for a factor \(CT_2\) multiplying a product that contains the
crown. A product of complete transposition graphs with no crown is
covered by the companion paper's product theorem, hypercubes included,
and never reaches this lemma.

The case is not vacuous, and the smallest instance sits at ground order
nine. Composing the path \(P_4\), whose degree sequence is \((2,2,1,1)\)
and whose realization graph is \(K_2 = CT_2\), with the five-cycle,
whose degree sequence is \((2,2,2,2,2)\) and whose realization graph is
the crown, gives

\[
d = (7,7,4,4,4,4,4,1,1),
\]

and \(G(d)\) is the prism over the crown: 24 vertices, 6-regular,
bipartite. No smaller ground order produces one, since every nontrivial
bipartite realization graph on at most eight ground vertices has 2, 4,
6, 12 or 24 vertices, and the only 24 among them is \(CT_4\), which the
product theorem covers. The two are distinguished by their four-cycle
counts, 210 for the prism over the crown against 162 for \(CT_4\).

\hypertarget{refs}{}
\begin{CSLReferences}{0}{0}
\leavevmode\vadjust pre{\hypertarget{ref-AL1989}{}}%
\CSLLeftMargin{{[}1{]} }%
\CSLRightInline{B. Alspach, G. Liu,
\href{https://doi.org/10.1007/BF01788672}{Paths and cycles in matroid
base graphs}, Graphs and Combinatorics 5 (1989) 207--211.}

\leavevmode\vadjust pre{\hypertarget{ref-Ardila2003}{}}%
\CSLLeftMargin{{[}2{]} }%
\CSLRightInline{F. Ardila,
\href{https://doi.org/10.1016/s0097-3165(03)00121-3}{The {C}atalan
matroid}, Journal of Combinatorial Theory, Series A 104 (2003) 49--62.}

\leavevmode\vadjust pre{\hypertarget{ref-AP1999}{}}%
\CSLLeftMargin{{[}3{]} }%
\CSLRightInline{S.R. Arikati, U.N. Peled,
\href{https://doi.org/10.1016/s0024-3795(98)10229-x}{The realization
graph of a degree sequence with majorization gap 1 is {H}amiltonian},
Linear Algebra and Its Applications 290 (1999) 213--235.}

\leavevmode\vadjust pre{\hypertarget{ref-BY2026}{}}%
\CSLLeftMargin{{[}4{]} }%
\CSLRightInline{J.S. Baggett, H. Yan,
\href{https://doi.org/10.48550/arxiv.2607.13165}{Interchange graphs of
(0,1)-matrices are maximally {H}amiltonian}, (2026).}

\leavevmode\vadjust pre{\hypertarget{ref-Barrus2016}{}}%
\CSLLeftMargin{{[}5{]} }%
\CSLRightInline{M.D. Barrus,
\href{https://doi.org/10.1016/j.disc.2016.03.012}{On realization graphs
of degree sequences}, Discrete Mathematics 339 (2016) 2146--2152.}

\leavevmode\vadjust pre{\hypertarget{ref-BJR2009}{}}%
\CSLLeftMargin{{[}6{]} }%
\CSLRightInline{L.J. Billera, N. Jia, V. Reiner,
\href{https://doi.org/10.1016/j.ejc.2008.12.007}{A quasisymmetric
function for matroids}, European Journal of Combinatorics 30 (2009)
1727--1757.}

\leavevmode\vadjust pre{\hypertarget{ref-BM2006}{}}%
\CSLLeftMargin{{[}7{]} }%
\CSLRightInline{J.E. Bonin, A. de Mier,
\href{https://doi.org/10.1016/j.ejc.2005.01.008}{Lattice path matroids:
{S}tructural properties}, European Journal of Combinatorics 27 (2006)
701--738.}

\leavevmode\vadjust pre{\hypertarget{ref-brualdiMatricesZerosOnes1980}{}}%
\CSLLeftMargin{{[}8{]} }%
\CSLRightInline{R.A. Brualdi,
\href{https://doi.org/10.1016/0024-3795(80)90105-6}{Matrices of zeros
and ones with fixed row and column sum vectors}, Linear Algebra and Its
Applications 33 (1980) 159--231.}

\leavevmode\vadjust pre{\hypertarget{ref-brualdiCombinatorialMatrixClasses2006}{}}%
\CSLLeftMargin{{[}9{]} }%
\CSLRightInline{R.A. Brualdi,
\href{https://doi.org/10.1017/CBO9780511721182}{Combinatorial matrix
classes}, Cambridge University Press, Cambridge, UK ; New York, 2006.}

\leavevmode\vadjust pre{\hypertarget{ref-colemanPairedn1ton1Disjoint2025}{}}%
\CSLLeftMargin{{[}10{]} }%
\CSLRightInline{A. Coleman, G. Fischberg, C. Gong, J. Harrington, T.W.H.
Wong, \href{https://doi.org/10.1016/j.dam.2025.07.038}{Paired
(n-1)-to-(n-1) disjoint path covers in bipartite transposition-like
graphs}, Discrete Applied Mathematics 376 (2025) 449--461.}

\leavevmode\vadjust pre{\hypertarget{ref-Crapo1965}{}}%
\CSLLeftMargin{{[}11{]} }%
\CSLRightInline{H.H. Crapo,
\href{https://doi.org/10.6028/jres.069B.003}{Single-element extensions
of matroids}, Journal of Research of the National Bureau of Standards
Section B Mathematics and Mathematical Physics 69B (1965) 55--65.}

\leavevmode\vadjust pre{\hypertarget{ref-CS2005}{}}%
\CSLLeftMargin{{[}12{]} }%
\CSLRightInline{H. Crapo, W. Schmitt,
\href{https://doi.org/10.1016/j.ejc.2004.05.006}{A free subalgebra of
the algebra of matroids}, European Journal of Combinatorics 26 (2005)
1066--1085.}

\leavevmode\vadjust pre{\hypertarget{ref-EG1960}{}}%
\CSLLeftMargin{{[}13{]} }%
\CSLRightInline{P. Erdős, T. Gallai, Gráfok előírt fokú pontokkal,
Matematikai Lapok 11 (1960) 264--274.}

\leavevmode\vadjust pre{\hypertarget{ref-FHPR2019}{}}%
\CSLLeftMargin{{[}14{]} }%
\CSLRightInline{C.G. Fernandes, C. Hernández-Vélez, J.C. de Pina, J.L.
Ramírez Alfonsín,
\href{https://doi.org/10.1007/s00373-019-02011-8}{Counting {H}amiltonian
{C}ycles in the {M}atroid {B}asis {G}raph}, Graphs and Combinatorics 35
(2019) 539--550.}

\leavevmode\vadjust pre{\hypertarget{ref-FH1977}{}}%
\CSLLeftMargin{{[}15{]} }%
\CSLRightInline{S. Földes, P.L. Hammer, Split graphs, in: Proceedings of
the Eighth Southeastern Conference on Combinatorics, Graph Theory and
Computing, Baton Rouge, 1977: pp. 311--315.}

\leavevmode\vadjust pre{\hypertarget{ref-FHM1965}{}}%
\CSLLeftMargin{{[}16{]} }%
\CSLRightInline{D.R. Fulkerson, A.J. Hoffman, M.H. McAndrew,
\href{https://doi.org/10.4153/cjm-1965-016-2}{Some {P}roperties of
{G}raphs with {M}ultiple {E}dges}, Canadian Journal of Mathematics 17
(1965) 166--177.}

\leavevmode\vadjust pre{\hypertarget{ref-GGMS1987}{}}%
\CSLLeftMargin{{[}17{]} }%
\CSLRightInline{I.M. Gelfand, R.M. Goresky, R.D. MacPherson, V.V.
Serganova,
\href{https://doi.org/10.1016/0001-8708(87)90059-4}{Combinatorial
geometries, convex polyhedra, and {S}chubert cells}, Advances in
Mathematics 63 (1987) 301--316.}

\leavevmode\vadjust pre{\hypertarget{ref-hladikRealizationGraphEvery2026}{}}%
\CSLLeftMargin{{[}18{]} }%
\CSLRightInline{P. Hladík, J. Fink,
\href{https://doi.org/10.48550/arXiv.2607.18146}{The realization graph
of every degree sequence has a {H}amilton path}, (2026).}

\leavevmode\vadjust pre{\hypertarget{ref-Klivans2003}{}}%
\CSLLeftMargin{{[}19{]} }%
\CSLRightInline{C.J. Klivans, Combinatorial properties of shifted
complexes, Ph.D. thesis, Massachusetts Institute of Technology, 2003.}

\leavevmode\vadjust pre{\hypertarget{ref-LZ1994}{}}%
\CSLLeftMargin{{[}20{]} }%
\CSLRightInline{X. Li, F. Zhang,
\href{https://doi.org/10.1016/0166-218x(94)90099-x}{Hamiltonicity of a
type of interchange graphs}, Discrete Applied Mathematics 51 (1994)
107--111.}

\leavevmode\vadjust pre{\hypertarget{ref-MP1995}{}}%
\CSLLeftMargin{{[}21{]} }%
\CSLRightInline{N.V.R. Mahadev, U.N. Peled,
\href{https://doi.org/10.1016/s0167-5060(13)71063-x}{Threshold {G}raphs
and {R}elated {T}opics}, North-Holland (Elsevier), 1995.}

\leavevmode\vadjust pre{\hypertarget{ref-MM2024}{}}%
\CSLLeftMargin{{[}22{]} }%
\CSLRightInline{A. Merino, T. Mütze,
\href{https://doi.org/10.1137/23m1612019}{Traversing {C}ombinatorial
0/1-{P}olytopes via {O}ptimization}, SIAM Journal on Computing 53 (2024)
1257--1292.}

\leavevmode\vadjust pre{\hypertarget{ref-MNW2024}{}}%
\CSLLeftMargin{{[}23{]} }%
\CSLRightInline{A. Merino, Namrata, A. Williams,
\href{https://doi.org/10.1007/978-981-97-0566-5_9}{On the~{H}ardness
of~{G}ray {C}ode {P}roblems for~{C}ombinatorial {O}bjects}, in: WALCOM:
Algorithms and Computation: 18th International Conference and Workshops
on Algorithms and Computation, WALCOM 2024, Kanazawa, Japan, March
18--20, 2024, Proceedings, Springer Nature Singapore, 2024: pp.
103--117.}

\leavevmode\vadjust pre{\hypertarget{ref-Mutze2023}{}}%
\CSLLeftMargin{{[}24{]} }%
\CSLRightInline{T. Mütze,
\href{https://doi.org/10.37236/11023}{Combinatorial {G}ray {C}odes ---
an {U}pdated {S}urvey}, The Electronic Journal of Combinatorics DS26
(2023).}

\leavevmode\vadjust pre{\hypertarget{ref-NP1984}{}}%
\CSLLeftMargin{{[}25{]} }%
\CSLRightInline{D.J. Naddef, W.R. Pulleyblank,
\href{https://doi.org/10.1016/0095-8956(84)90043-1}{Hamiltonicity in
(0--1)-polyhedra}, Journal of Combinatorial Theory, Series B 37 (1984)
41--52.}

\leavevmode\vadjust pre{\hypertarget{ref-NP1981}{}}%
\CSLLeftMargin{{[}26{]} }%
\CSLRightInline{D. Naddef, W.R. Pulleyblank,
\href{https://doi.org/10.1016/0095-8956(81)90032-0}{Hamiltonicity and
combinatorial polyhedra}, Journal of Combinatorial Theory, Series B 31
(1981) 297--312.}

\leavevmode\vadjust pre{\hypertarget{ref-Partida2026}{}}%
\CSLLeftMargin{{[}27{]} }%
\CSLRightInline{E. Partida,
\href{https://doi.org/10.1016/j.dam.2025.08.038}{Shifted and threshold
matroids}, Discrete Applied Mathematics 378 (2026) 522--537.}

\leavevmode\vadjust pre{\hypertarget{ref-RG1979}{}}%
\CSLLeftMargin{{[}28{]} }%
\CSLRightInline{E. Ruch, I. Gutman, The branching extent of graphs,
Journal of Combinatorics, Information and System Sciences 4 (1979)
285--295.}

\leavevmode\vadjust pre{\hypertarget{ref-Ruskey2003}{}}%
\CSLLeftMargin{{[}29{]} }%
\CSLRightInline{F. Ruskey, Combinatorial {G}eneration, (2003).}

\leavevmode\vadjust pre{\hypertarget{ref-RSW2012}{}}%
\CSLLeftMargin{{[}30{]} }%
\CSLRightInline{F. Ruskey, J. Sawada, A. Williams,
\href{https://doi.org/10.1016/j.jcta.2011.07.005}{Binary bubble
languages and cool-lex order}, Journal of Combinatorial Theory, Series A
119 (2012) 155--169.}

\leavevmode\vadjust pre{\hypertarget{ref-ryserCombinatorialPropertiesMatrices1957}{}}%
\CSLLeftMargin{{[}31{]} }%
\CSLRightInline{H.J. Ryser,
\href{https://doi.org/10.4153/CJM-1957-044-3}{Combinatorial {P}roperties
of {M}atrices of {Z}eros and {O}nes}, Canadian Journal of Mathematics 9
(1957) 371--377.}

\leavevmode\vadjust pre{\hypertarget{ref-SP2025}{}}%
\CSLLeftMargin{{[}32{]} }%
\CSLRightInline{V.N. Schvöllner, A. Pastine,
\href{https://doi.org/10.48550/arxiv.2511.23327}{The 2-switch-degree of
a graph}, (2025).}

\leavevmode\vadjust pre{\hypertarget{ref-Sohoni1999}{}}%
\CSLLeftMargin{{[}33{]} }%
\CSLRightInline{M. Sohoni,
\href{https://doi.org/10.1007/s003730050032}{Rapid {M}ixing of {S}ome
{L}inear {M}atroids and {O}ther {C}ombinatorial {O}bjects}, Graphs and
Combinatorics 15 (1999) 93--107.}

\leavevmode\vadjust pre{\hypertarget{ref-Tyshkevich2000}{}}%
\CSLLeftMargin{{[}34{]} }%
\CSLRightInline{R. Tyshkevich,
\href{https://doi.org/10.1016/s0012-365x(99)00381-7}{Decomposition of
graphical sequences and unigraphs}, Discrete Mathematics 220 (2000)
201--238.}

\leavevmode\vadjust pre{\hypertarget{ref-Welsh1969}{}}%
\CSLLeftMargin{{[}35{]} }%
\CSLRightInline{D.J.A. Welsh,
\href{https://doi.org/10.1016/s0021-9800(69)80094-3}{A bound for the
number of matroids}, Journal of Combinatorial Theory 6 (1969) 313--316.}

\end{CSLReferences}

\end{document}